\documentclass[11pt]{article}
\usepackage[a4paper,margin=2.5cm]{geometry}
\usepackage{indentfirst}
\usepackage{amsmath,amssymb,amsthm,mathtools}
\usepackage{mathrsfs}
\usepackage{enumitem,booktabs,microtype}
\usepackage{tikz,placeins}
\usepackage{xcolor,needspace}
\usepackage[normalem]{ulem}

\usepackage[
  colorlinks=true,
  allcolors=black,
  bookmarksnumbered=true,
  pagebackref=true
]{hyperref}

\renewcommand*{\backref}[1]{}
\renewcommand*{\backrefalt}[4]{%
  \ifcase#1\relax
  \or
    (Cited on page~#2.)%
  \else
    (Cited on pages~#2.)%
  \fi
}

\usepackage[nottoc]{tocbibind}
\usepackage{array}

\newcommand{\soc}{\operatorname{soc}}
\newcommand{\rad}{\operatorname{rad}}
\renewcommand{\top}{\operatorname{top}}
\newcommand{\Hom}{\operatorname{Hom}}
\newcommand{\End}{\operatorname{End}}
\newcommand{\Ext}{\operatorname{Ext}}

\newcommand{\add}{\operatorname{add}}
\renewcommand{\Im}{\operatorname{Im}}
\renewcommand{\mod}{\operatorname{mod}}
\newcommand{\smod}{\operatorname{\underline{mod}}}
\newcommand{\sHom}{\operatorname{\underline{Hom}}}
\newcommand{\ismod}{\operatorname{\overline{mod}}}
\newcommand{\isHom}{\operatorname{\overline{Hom}}}
\newcommand{\isEnd}{\operatorname{\overline{End}}}
\newcommand{\id}{\operatorname{id}}
\newcommand{\ra}{\rightarrow}

\newcommand{\vv}{\mathsf{v}}
\newcommand{\ww}{\mathsf{w}}

\renewcommand{\cal}[1]{\mathcal #1}   % 只对第一个字母有效; \cal AB=\cal{AB}=\mathcal{A}B

\providecommand{\varLambda}{\Lambda}
\providecommand{\varGamma}{\Gamma}
\newtheorem{theorem}{Theorem}[section]
\newtheorem{proposition}[theorem]{Proposition}
\newtheorem{lemma}[theorem]{Lemma}
\newtheorem{corollary}[theorem]{Corollary}
\theoremstyle{definition}
\newtheorem{definition}[theorem]{Definition}
\newtheorem{remark}[theorem]{Remark}
\newtheorem{example}[theorem]{Example}
\newtheorem{question}[theorem]{Question}
\setlist[enumerate]{label=\textup{(\arabic*)},itemsep=3pt,topsep=4pt}
\title{Stable equivalence of quadratic monomial algebras}
\author{Pengyun Chen$^{a}$ and Bohan Xing$^{a,*}$}
\date{}

\begin{document}
\maketitle

\begingroup
\renewcommand{\thefootnote}{}

\makeatletter
\renewcommand{\@makefntext}[1]{%
  \setlength{\parindent}{2em}%
  \noindent\hspace*{\parindent}#1}
\makeatother

\begin{NoHyper}
\footnotetext{%
	The version of \today.
  \par
  $^{a}$Pengyun Chen and Bohan Xing,
  School of Mathematical Sciences,
  Laboratory of Mathematics and Complex Systems,
  Beijing Normal University,
  Beijing 100875, P. R. China.
  \par
  E-mail addresses:
  \texttt{pychen@mail.bnu.edu.cn} (P. Chen);
  \texttt{bhxing@mail.bnu.edu.cn} (B. Xing).
  \par
  $^{*}$Corresponding author.
}
\end{NoHyper}
\endgroup

\begin{abstract}
It is well known that gentle algebras are closed under derived equivalence.
However, it is easy to verify that they are not closed under stable equivalence.
We prove that the classes of quadratic monomial, almost gentle and string almost gentle algebras are closed under stable equivalence.
In fact, in each of these cases, stable equivalence is characterized by isomorphism after deleting all nodes and removing semisimple direct summands.
The underlying isomorphism theorem gives a positive answer to Mart\'inez-Villa's question on selfinjective-free algebras for these classes.
We also prove that the nonprojective simple modules form the unique weakly simple-minded system of a quadratic monomial algebra.
Consequently, every stable equivalence of Morita type involving such an algebra lifts to a Morita equivalence when neither algebra has semisimple direct summands.
Finally, we describe the effect of node deletion on stable categories of Gorenstein-projective modules and singularity categories, and give combinatorial descriptions of the equivalences for almost gentle algebras using marked configurations.
\end{abstract}

\tableofcontents
%\clearpage

\section{Introduction}

	Gentle algebras arose in the study of iterated tilted algebras of types $A_n$ and $\widetilde{A}_n$ \cite{AH81,AS87}. They form a class of quadratic monomial special biserial algebras \cite{SW83,WW85,Sch15} and are closed under derived equivalence \cite{SZ03}. In recent years, gentle algebras have attracted increasing attention through the development of surface models for their derived categories \cite{OPS} and connections with partially wrapped Fukaya categories \cite{HKK17,LP20}. These models give concrete geometric descriptions of indecomposable objects and morphisms in their derived categories and provide tools for studying derived equivalences. For an overview of these developments, we refer to Schroll's survey \cite{Sch26}.

	In parallel, other triangulated categories associated with gentle algebras and their generalizations have also received considerable attention. These include the stable categories of Gorenstein-projective modules \cite{Buc87,Hap91} and singularity categories \cite{Buc87,Orl04}, which are triangle equivalent for Gorenstein algebras. Kalck \cite{Kal15} gave an explicit description for gentle algebras, while the work of Chen, Shen and Zhou \cite{CSZ18} and Chen \cite{Che18} provided descriptions for quadratic monomial algebras using relation quivers. More generally, Honma and Usui \cite{HU25} described the stable categories of Gorenstein-projective modules over arbitrary monomial algebras in terms of stable module categories of selfinjective Nakayama algebras. These results connect the structure of the associated triangulated categories with the combinatorics of paths and relations, making their study more concrete.

	However, equivalences between the ordinary stable module categories of these algebras remain less well understood. We note that gentle algebras are not closed under stable equivalence: node deletion gives a stable equivalence between the Kronecker algebra $k(1\rightrightarrows 2)$ and $k[x,y]/(x,y)^2$ \cite{Mar80}, and the latter is string almost gentle \cite{GS19,FGR22} but not gentle. This motivates us to consider almost gentle algebras \cite{GS19}, which are quadratic monomial special multiserial algebras \cite{HW84,GS16MS}, and more generally arbitrary quadratic monomial algebras. Our results complete the picture of closure under equivalences summarized in Table~\ref{table of closure properties}.

	\begin{table}[!htbp]
		\centering
		\small
		\renewcommand{\arraystretch}{1.25}
		\begin{tabular}{@{}lcccc@{}}
			\toprule
			& Derived & Stable & \shortstack{Stable of\\Morita type} & Singular \\
			\midrule
			Gentle & $\checkmark$ & $\times$ & $\checkmark$ & $\times$ \\
			String almost gentle & $\times$ & $\checkmark$ & $\checkmark$ & $\times$ \\
			Almost gentle & $\times$ & $\checkmark$ & $\checkmark$ & $\times$ \\
			Quadratic monomial & $\times$ & $\checkmark$ & $\checkmark$ & $\times$ \\
			\bottomrule
		\end{tabular}
		\caption{Closure under equivalences among finite-dimensional basic algebras.}
		\label{table of closure properties}
	\end{table}

	The following examples account for the negative entries in Table~\ref{table of closure properties}.

	\begin{example}\label{example of closure counterexamples}
		Consider the following quivers:
		
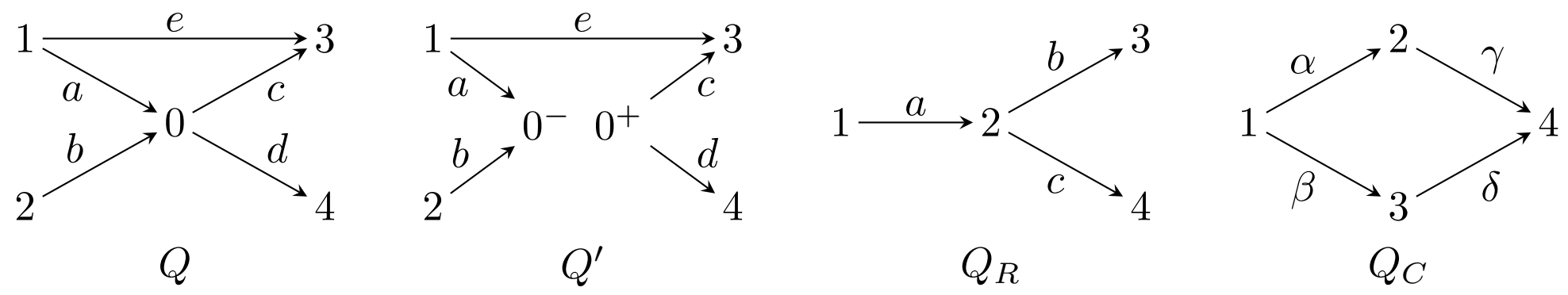
\[
			\begin{tikzpicture}[
				x=1cm,y=1cm,font=\small,>=stealth,
				every node/.style={inner sep=1.5pt}
				]
				\begin{scope}
					\node (u1) at (0,.7) {$1$};
					\node (u2) at (0,-.7) {$2$};
					\node (v) at (1.25,0) {$0$};
					\node (w1) at (2.5,.7) {$3$};
					\node (w2) at (2.5,-.7) {$4$};
					\draw[->] (u1)--node[pos=.42,below left] {$a$}(v);
					\draw[->] (u2)--node[pos=.42,above left] {$b$}(v);
					\draw[->] (v)--node[pos=.58,below right] {$c$}(w1);
					\draw[->] (v)--node[pos=.58,above right] {$d$}(w2);
					\draw[->] (u1)--node[above] {$e$}(w1);
					\node at (1.25,-1.2) {$Q$};
				\end{scope}

				\begin{scope}[xshift=3.4cm]
					\node (u1) at (0,.7) {$1$};
					\node (u2) at (0,-.7) {$2$};
					\node (vm) at (.95,0) {$0^-$};
					\node (vp) at (1.55,0) {$0^+$};
					\node (w1) at (2.5,.7) {$3$};
					\node (w2) at (2.5,-.7) {$4$};
					\draw[->] (u1)--node[pos=.4,below left] {$a$}(vm);
					\draw[->] (u2)--node[pos=.4,above left] {$b$}(vm);
					\draw[->] (vp)--node[pos=.6,below right] {$c$}(w1);
					\draw[->] (vp)--node[pos=.6,above right] {$d$}(w2);
					\draw[->] (u1)--node[above] {$e$}(w1);
					\node at (1.25,-1.2) {$Q'$};
				\end{scope}

				\begin{scope}[xshift=6.8cm]
					\node (u) at (0,0) {$1$};
					\node (v) at (1.25,0) {$2$};
					\node (w1) at (2.5,.7) {$3$};
					\node (w2) at (2.5,-.7) {$4$};
					\draw[->] (u)--node[above] {$a$}(v);
					\draw[->] (v)--node[pos=.55,above left] {$b$}(w1);
					\draw[->] (v)--node[pos=.55,below left] {$c$}(w2);
					\node at (1.25,-1.2) {$Q_R$};
				\end{scope}

				\begin{scope}[xshift=10.2cm]
					\node (u) at (0,0) {$1$};
					\node (v1) at (1.25,.7) {$2$};
					\node (v2) at (1.25,-.7) {$3$};
					\node (w) at (2.5,0) {$4$};
					\draw[->] (u)--node[above left] {$\alpha$}(v1);
					\draw[->] (u)--node[below left] {$\beta$}(v2);
					\draw[->] (v1)--node[above right] {$\gamma$}(w);
					\draw[->] (v2)--node[below right] {$\delta$}(w);
					\node at (1.25,-1.2) {$Q_C$};
				\end{scope}
			\end{tikzpicture}
		\]
		Put
		\[
			\begin{aligned}
				A&=kQ/(ac,ad,bc,bd),&
				G&=kQ',\\
				R&=kQ_R/(ab,ac),&
				C&=kQ_C/(\alpha\gamma-\beta\delta).
			\end{aligned}
		\]
		The algebras $A$ and $R$ are string almost gentle but not gentle, $G$ is gentle, and $C$ is not monomial.
		There are equivalences
		\[
			\smod A\simeq\smod G,\qquad
			D^b(\mod R)\simeq D^b(\mod C),\qquad
			D_{\mathrm{sg}}(C)\simeq D_{\mathrm{sg}}(k)=0;
		\]
		cf.~\cite{Mar80,Kel07}.
		These pairs give all negative entries of Table~\ref{table of closure properties}.
	\end{example}

	In general, stable equivalences remain less well understood than derived equivalences, which admit a characterization in terms of tilting complexes \cite{Ric89,Kel94}.
	Although stable equivalences can arise from various constructions, including node deletion \cite{Mar80}, certain socle equivalences of selfinjective algebras \cite{Ohn02}, and tensor functors inducing stable equivalences of Morita type \cite{Bro94}, there is no comparable general description of stable equivalence.

	For quadratic monomial algebras, we show that, up to algebra isomorphism and semisimple direct summands, stable equivalence is generated by node deletion and its inverse. This gives a concrete description of their stable equivalence classes. We also describe how node deletion affects stable categories of Gorenstein-projective modules and singularity categories in Subsection~\ref{subsection of triangulated node deletion}. Our main result is the following.

	\begin{theorem}[see Theorems~\ref{theorem of quadratic monomial stable closure},
		\ref{theorem of AG stable closure}, and~\ref{theorem of SAG stable closure},
		and Corollary~\ref{corollary of stable Gorenstein quotient invariant}]
		\label{theorem of main results}
		Let $\varLambda$ and $\varGamma$ be finite-dimensional basic algebras over an algebraically closed field $k$, and suppose that $\varLambda$ is quadratic monomial. Then the following are equivalent:
		\begin{enumerate}
			\item $\varLambda$ and $\varGamma$ are stably equivalent.
			\item After removing their semisimple direct summands,
			$\varLambda$ and $\varGamma$ are stably equivalent of node type,
			up to algebra isomorphism.
			\item The algebras obtained from $\varLambda$ and $\varGamma$
			by deleting all nodes and removing all semisimple direct summands
			are isomorphic.
		\end{enumerate}
		Under these conditions, $\varGamma$ is also quadratic monomial. If $\varLambda$ is almost gentle \textup{(}resp.\ string almost gentle\textup{)}, then so is $\varGamma$. Moreover, there is a triangle equivalence
		\[
			\underline{\operatorname{Gproj}}\varLambda/
			\mathcal T_{\varLambda}
			\simeq
			\underline{\operatorname{Gproj}}\varGamma/
			\mathcal T_{\varGamma},
		\]
		where $\mathcal T_\varLambda$ is the thick subcategory of $\underline{\operatorname{Gproj}}\varLambda$ generated by the Gorenstein-projective node simples.
	\end{theorem}

	The key step is to show that two stably equivalent algebras with no nodes or semisimple summands are isomorphic if one is quadratic monomial; see Theorem~\ref{theorem of quadratic monomial rigidity}. 
	We prove this theorem by induction on the number of arrows. The results of Mart\'inez-Villa \cite[Section~3]{Mar90} show that these two algebras have isomorphic quivers. 
	If they have no projective-injective modules, then the diagrammatic morphisms show that maps from injective modules to projective modules of the quadratic monomial algebra have semisimple image.
	Using the induced equivalence of injectively stable categories, we then recover the relations and obtain an algebra isomorphism preserving the vertex correspondence; see Lemma~\ref{lemma of stable equivalences induces isomorphisms}.
	If one of the algebras has a projective-injective module, then the top of an indecomposable projective-injective module corresponds to a source vertex with exactly one outgoing arrow; see Proposition~\ref{prop of proj-injective modules}. 
	Removing this vertex and its outgoing arrow, then deleting any nodes that arise and removing semisimple summands, decreases the number of arrows by one.
	These operations preserve stable equivalence and the correspondence between the vertices of the two quivers, and allow us to apply the induction hypothesis and obtain the required isomorphism by reversing node deletion and taking the corresponding one-point extensions.
	This answers Mart\'inez-Villa's question on selfinjective-free algebras affirmatively in the quadratic monomial case; see Corollary~\ref{corollary of quadratic monomial Martinez-Villa question}.

	\medskip

	We also describe stable equivalences of Morita type for quadratic monomial algebras using simple-minded systems. These were introduced by Koenig and Liu \cite{KL12} as analogues in stable module categories of simple-minded collections in derived categories \cite{Ric02,KY14}. In contrast to the rich family of simple-minded collections for gentle algebras \cite{Cha26}, we show that the stable module category of any quadratic monomial algebra has a unique weakly simple-minded system. It consists precisely of the nonprojective simple modules. Combining this uniqueness with the lifting criteria \cite{Lin96,Liu03,Liu06}, we obtain the following result.

	\begin{theorem}[see Theorems~\ref{theorem of quadratic monomial sms}
		and~\ref{corollary of AG Morita rigidity}]
		\label{theorem of sms and Morita type}
		Let $A$ be a finite-dimensional basic quadratic monomial algebra over an algebraically closed field $k$.
		\begin{enumerate}
			\item Up to isomorphism, the nonprojective simple $A$-modules form the unique weakly simple-minded system in $\smod A$.  In particular, they form the unique simple-minded system.

			\item Let $B$ be a finite-dimensional basic $k$-algebra stably equivalent to $A$ of Morita type. If neither algebra has semisimple direct summands, every such equivalence lifts to a Morita equivalence. In particular, $A\cong B$.
		\end{enumerate}
	\end{theorem}

	The uniqueness statement also gives positive answers, for quadratic monomial algebras, to the questions on weak generation and cardinality in \cite{KL12} and to the simple-image problem of \cite{CKL15}; see Corollary~\ref{corollary of sms questions}.
	In particular, every stable equivalence involving a quadratic monomial algebra preserves the number of isomorphism classes of nonprojective simple modules, as predicted by the Auslander--Reiten conjecture.

	\medskip

	We further interpret these equivalences through combinatorial models. Although monomial algebras admit descriptions by fractional Brauer configurations \cite{LX26}, the resulting combinatorial data  can be difficult to represent diagrammatically. We therefore focus on almost gentle algebras, for which our construction recovers the marked configurations arising from the work of Green and Schroll \cite{GS19}.
	After omitting semisimple direct summands, stable equivalence classes correspond to equivalence classes under incidence comparisons, while stable equivalence classes of Morita type correspond to isomorphism classes of marked configurations.
	We also relate these configurations to relation quivers \cite{CSZ18,Che18}. Table~\ref{table of configuration relation comparison} summarizes the resulting descriptions.

	\begin{table}[htbp]
		\centering
		\small
		\setlength{\tabcolsep}{4pt}
		\renewcommand{\arraystretch}{1.16}
		\begin{tabular}{@{}
				>{\raggedright\arraybackslash}p{.26\textwidth}
				>{\raggedright\arraybackslash}p{.36\textwidth}
				>{\raggedright\arraybackslash}p{.32\textwidth}@{}}
			\toprule
			& Marked configuration $C_A$ & Relation quiver $\mathcal R_A$ \\
			\midrule
			Quiver vertices
			& Polygons & --- \\
			Quiver arrows
			& Uncut corners & Vertices \\
			Nonzero quadratic compositions
			& Consecutive uncut corners at the same configuration vertex
			& --- \\
			Quadratic relations
			& Transitions $h\rightsquigarrow g$ & Arrows \\
			\midrule
			Stable equivalence classes
			& Equivalence classes under incidence comparisons
			& --- \\
			Stable equivalence classes of Morita type
			& Isomorphism classes of marked configurations
			& --- \\
			Triangle equivalence classes of $\underline{\operatorname{Gproj}}A$
			& Multisets of lengths of perfect faces
			& Multisets of lengths of perfect components \\
			Singular equivalence classes (Hom-finite case)
			& Multisets of face lengths
			& Multisets of cycle lengths after deleting sources and sinks \\
			\bottomrule
		\end{tabular}
		\caption{Combinatorial descriptions for almost gentle algebras,
			with semisimple direct summands omitted.
			Dashes indicate data not determined by the relation quiver alone.}
		\label{table of configuration relation comparison}
	\end{table}

\begin{figure}[htbp]
		\centering
		\begin{tikzpicture}[scale=.95,>=stealth,font=\small,color=black,
			polygon/.style={fill=blue!4,draw=gray!65},
			occurrence/.style={draw=blue!65!black},
			quiver arrow/.style={->,thick},
			pics/digon/.style args={#1}{code={
					\filldraw[polygon] (0,0)
						.. controls (.55,.18) and ({#1-.55},.18) .. (#1,0)
						.. controls ({#1-.55},-.18) and (.55,-.18) .. cycle;
			}}]

			% The node polygon.
			\filldraw[polygon]
				(0,1)--(3.6,1)--(3.6,-1.2)--(0,-1.2)--cycle;
			\node at (1.8,-.1) {$P_0$};

			% The two polygons after node deletion.
			\begin{scope}[xshift=8cm]
				\filldraw[polygon] (0,1)--(.45,1)
					.. controls (1.2,1) and (1.2,-1.2) .. (.45,-1.2)
					--(0,-1.2)--cycle;
				\filldraw[polygon] (3.6,1)--(3.15,1)
					.. controls (2.4,1) and (2.4,-1.2) .. (3.15,-1.2)
					--(3.6,-1.2)--cycle;
				\node at (.47,-.1) {$P_{0^-}$};
				\node at (3.13,-.1) {$P_{0^+}$};
			\end{scope}

			% The remaining data are identical on both sides.
			\foreach \offset in {0,8}{
				\begin{scope}[xshift=\offset cm]
					\pic[rotate=45] at (0,1) {digon=2.545584};
					\pic[rotate=-45] at (1.8,2.8) {digon=2.545584};
					\pic[rotate=-135] at (0,-1.2) {digon=1.697056};
					\pic[rotate=-45] at (3.6,-1.2) {digon=1.697056};

					\node[above left] at (.9,1.9) {$P_1$};
					\node[above right] at (2.7,1.9) {$P_3$};
					\node[above left] at (-.6,-1.8) {$P_2$};
					\node[above right] at (4.2,-1.8) {$P_4$};

					% Occurrences.
					\foreach \xx/\yy/\start/\finish in {
						0/1/27/63,0/1/0/-90,
						0/-1.2/0/90,0/-1.2/207/243,
						3.6/1/117/153,3.6/1/180/270,
						3.6/-1.2/90/180,3.6/-1.2/-63/-27,
						1.8/2.8/207/243,1.8/2.8/297/333,
						-1.2/-2.4/27/63,4.8/-2.4/117/153}{
						\begin{scope}[shift={(\xx,\yy)}]
							\draw[occurrence] (\start:.23)
								arc[start angle=\start,end angle=\finish,radius=.23];
						\end{scope}
					}

					% Quiver arrows in the unmarked corners.
					\foreach \xx/\yy/\start/\finish in {
						0/1/23/4,3.6/1/176/157,
						0/-1.2/165/105,3.6/-1.2/75/15,
						1.8/2.8/145/35}{
						\begin{scope}[shift={(\xx,\yy)}]
							\draw[quiver arrow] (\start:.60)
								arc[start angle=\start,end angle=\finish,radius=.60];
						\end{scope}
					}
					\node at (.86,1.19) {$a$};
					\node at (-.57,-.62) {$b$};
					\node at (2.74,1.19) {$c$};
					\node at (4.17,-.62) {$d$};
					\node at (1.8,3.63) {$e$};

					% Marked corners.
					\foreach \xx/\yy in {
						-.43/.90,.38/-1.58,
						4.03/.90,3.22/-1.58,1.8/2.32}
					\node at (\xx,\yy) {$\times$};

					% Nontruncated configuration vertices.
					\foreach \xx/\yy/\z/\anchor in {
						0/1/a/south east,0/-1.2/b/east,
						3.6/1/c/south west,3.6/-1.2/d/west,
						1.8/2.8/e/south}{
						\fill (\xx,\yy) circle (1.2pt);
						\node[anchor=\anchor,inner sep=3pt]
							at (\xx,\yy) {$\vv_{\z}$};
					}

					% Truncated vertices.
					\fill (-1.2,-2.4) circle (1.2pt)
						(4.8,-2.4) circle (1.2pt);
					\node[below] at (-1.2,-2.4) {$\ww_2$};
					\node[below] at (4.8,-2.4) {$\ww_4$};
				\end{scope}
			}

			\draw[->,thick] (4.95,-.1)--(6.65,-.1)
				node[midway,above,font=\scriptsize] {node deletion};
			\node at (1.8,-3.05) {$C_A$};
			\node at (9.8,-3.05) {$C_G$};
		\end{tikzpicture}
		\caption{Node deletion for the algebras $A$ and $G$ in Example~\ref{example of closure counterexamples}.}
		\label{figure of introductory node deletion}
	\end{figure}
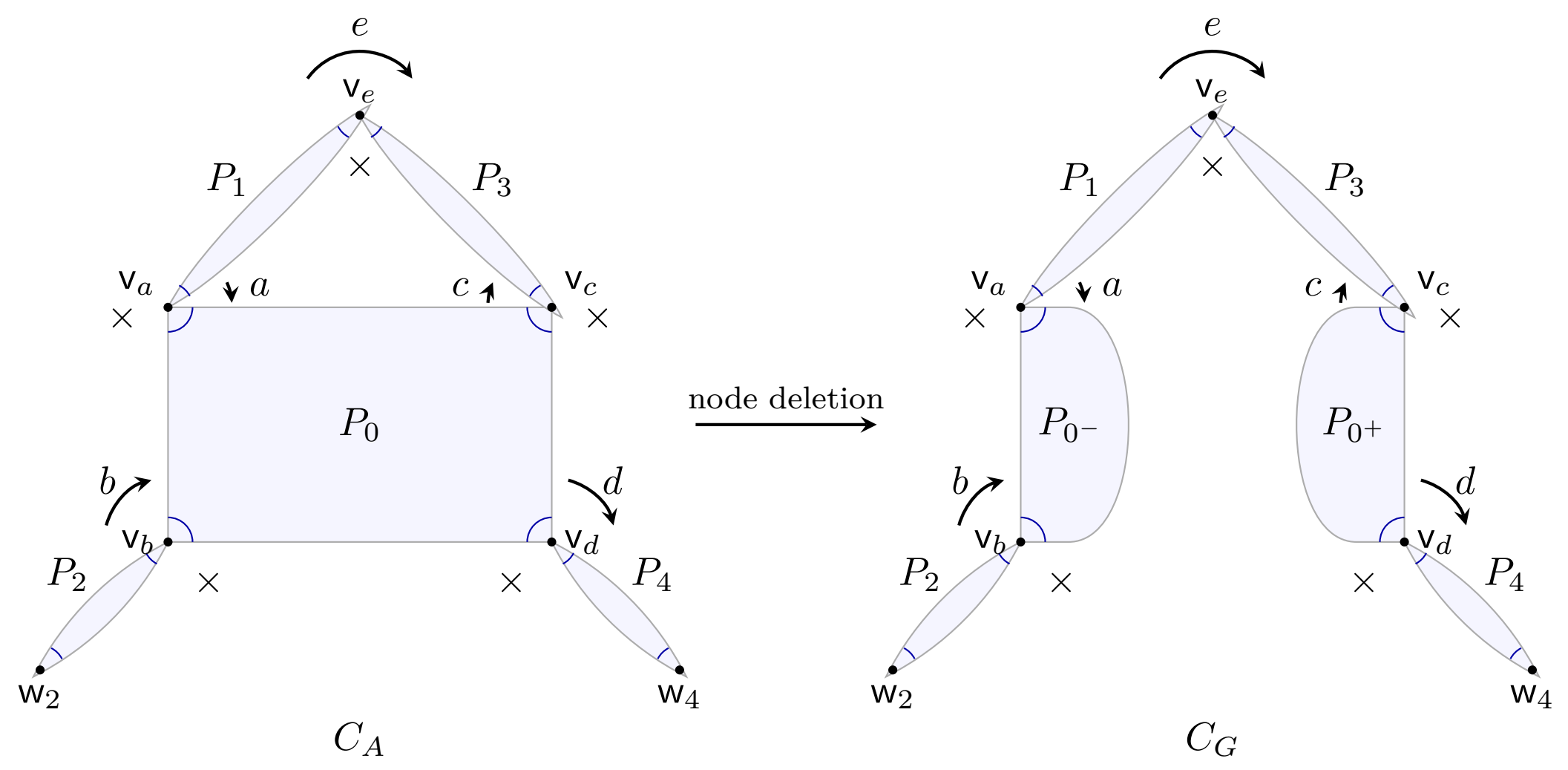

	For the algebras $A$ and $G$ in Example~\ref{example of closure counterexamples}, Figure~\ref{figure of introductory node deletion} realizes node deletion by splitting the polygon $P_0$ into $P_{0^-}$ and $P_{0^+}$.

	%\FloatBarrier

	\medskip
	\noindent\textbf{Outline.}
	In Section~\ref{section of preliminaries}, we recall stable module categories, node deletion, tree modules, and Mart\'inez-Villa's results on selfinjective-free algebras.
	In Section~\ref{section of Martinez-Villa question}, we prove the isomorphism theorem for algebras with no nodes or semisimple summands when one is quadratic monomial, and deduce that quadratic monomial algebras are closed under stable equivalence.
	In Section~\ref{section of sms and Morita type}, we prove the uniqueness of weakly simple-minded systems for quadratic monomial algebras and derive the lifting result for stable equivalences of Morita type.
	In Section~\ref{section of Gorenstein and singularity categories}, we recall the relation-quiver descriptions of stable Gorenstein-projective categories and singularity categories, and determine the effect of node deletion on these categories.
	In Section~\ref{section of marked configurations}, we prove stable closure for almost gentle and string almost gentle algebras and describe stable equivalences and stable equivalences of Morita type using marked configurations.
	We also compare their singularity classes and use faces to describe their stable Gorenstein-projective categories and Hom-finite singularity categories.

\section{Preliminaries}\label{section of preliminaries}

	Throughout, $k$ is an algebraically closed field, all algebras are finite-dimensional basic $k$-algebras, and all modules are finitely generated right modules.
	We denote the usual $k$-duality by $D=\Hom_k(-,k)$.
	For $\varLambda=kQ/I$ and a vertex $v_i$ of $Q$, let $e_i$ be the corresponding primitive idempotent, and put $P_i=e_i\varLambda$ and $S_i=\top P_i$.
	For a path $p$, we denote its starting and ending vertices by $s(p)$ and $t(p)$, respectively.
	We read products of paths in the order of traversal: $ab$ means first $a$ and then $b$.
	For homomorphisms, $gf$ denotes $f$ followed by $g$.

	{
		We denote by $\smod\varLambda$ the quotient of $\mod\varLambda$ by the ideal of homomorphisms factoring through projective modules, and by $\ismod\varLambda$ the quotient by those factoring through injective modules.
		Their Hom spaces are denoted by $\sHom_\varLambda$ and $\isHom_\varLambda$, respectively.
		Two algebras are called \textbf{stably equivalent} if their stable categories are equivalent as $k$-categories.
	}

	{
		Note that a quiver algebra $kQ/I$ is called \textbf{quadratic monomial} if $I$ is an admissible ideal generated by paths of length two.
	}

	\subsection{Node deletion}\label{subsection of node deletion}

		Recall that a simple module $S$ is a \textbf{node} if it is neither projective nor injective and the almost split sequence starting at $S$ has a projective middle term. For a bound quiver algebra, the corresponding vertex $v_i$ is a node precisely when it is neither a source vertex nor a sink vertex and every composition of two arrows through $v_i$ is zero (see, for example, \cite[Proposition~1.3]{Mar90} and \cite{MM05}).

		By \cite[Proposition~1.3(iii)]{Mar90}, a simple module $S$ which is neither projective nor injective is a node if and only if it is not a composition factor of $\rad P/\soc P$ for any nonsimple indecomposable projective module $P$.
        
		For a node $v_i$ of $\varLambda=kQ/I$, construct a quiver $Q'$ from $Q$ by replacing $v_i$ with a source vertex $v_{i^+}$ and a sink vertex $v_{i^-}$.
		The arrows of $Q'$ correspond to those of $Q$: the ending vertex of each arrow ending at $v_i$ is changed to $v_{i^-}$, and the starting vertex of each arrow starting at $v_i$ is changed to $v_{i^+}$.
		In particular, a loop at $v_i$ becomes an arrow $v_{i^+}\to v_{i^-}$.
		Retaining the relations whose paths have no internal vertex $v_i$ and removing those whose arrows are not composable in $Q'$ gives an ideal $I'$ and hence an algebra $\varLambda'=kQ'/I'$.
		This construction is called \textbf{node deletion}.
		The reverse operation glues the source vertex $v_{i^+}$ and the sink vertex $v_{i^-}$ together to form $v_i$ and sets every path of length two through $v_i$ equal to zero.
		By \cite[Theorem~2.10(c)]{Mar80}, node deletion induces an equivalence of $k$-categories $$\smod\varLambda\simeq\smod\varLambda'.$$ 		A finite composition of these equivalences and their inverses is said to be of \textbf{node type}.

		We now recall the node-deletion functor of Mart\'inez-Villa in \cite[Section~2]{Mar80}. Let $\varLambda=kQ/I$ be a bound quiver algebra,
		let $V$ be a set of its nodes, and let $\varLambda'$ be obtained by deleting the nodes in $V$. For a $\varLambda$-module $X$, write $X_i=Xe_i$, and let
		$\varphi_a:X_i\to X_j$ be the map $x\mapsto xa$ for an arrow $a:v_i\to v_j$. For each $v_j\in V$, put
		\[
			U_j(X)=\sum_{a:v_i\to v_j}\Im\varphi_a.
		\]
		The module $FX$ has vertex spaces
		\[
			(FX)_{j^+}=X_j/U_j(X),\qquad (FX)_{j^-}=U_j(X)
		\]
		at the vertices replacing $v_j$, and $(FX)_i=X_i$ if $v_i\notin V$. The arrow maps are induced by those of $X$: incoming arrows at $v_j$ take values in $U_j(X)$, and outgoing arrows factor through $X_j/U_j(X)$, since every composition of two arrows through $v_j$ is zero. On homomorphisms, $F$ takes the induced maps on these quotients and subspaces and leaves the other components unchanged. This defines a functor $F\colon\mod\varLambda\to\mod\varLambda'$. If $V=\varnothing$, we take $\varLambda'=\varLambda$ and $F$ to be the identity functor.

		We recall the following properties of $F$.

		\begin{proposition}[{\cite[Proposition~2.7, Lemma~2.8, and Theorem~2.10(c)]{Mar80}}]
			\label{proposition of node deletion}
			Let $\varLambda=kQ/I$ be a bound quiver algebra, let $V$ be a set of its nodes, let $\varLambda'$ be obtained by deleting the nodes in $V$, and let $F:\mod\varLambda\to\mod\varLambda'$ be the functor defined above. Then the following hold.
			\begin{enumerate}
				\item $F$ is full.
				\item If $f:X\to Y$ is a $\varLambda$-module homomorphism with $F(f)=0$, then $f$ factors through a projective module.
				\item $F$ induces an equivalence of $k$-categories $\underline F:\smod\varLambda\to\smod\varLambda'$.
			\end{enumerate}
		\end{proposition}

		\begin{remark}\label{remark of quadratic monomial node type}
			Deleting a node or gluing a source vertex and a sink vertex preserves the quadratic monomial property.
			Indeed, the relations after either operation are still generated by paths of length two.
		\end{remark}

	\subsection{Tree modules}\label{subsection of tree modules}
		Tree diagrams provide a description of projective and injective modules and their homomorphisms. We recall the construction from \cite{Kra91} in the notation used here.

		\begin{definition}\label{definition of diagrammatic basis}
			Let $\varLambda=kQ/I$ be a finite-dimensional monomial algebra.
			A \textbf{tree diagram} over $\varLambda$ is a finite oriented tree whose vertices and arrows are labelled by vertices and arrows of $Q$.
			The labelling preserves sources and targets for each arrow and satisfies the following conditions:
			\begin{enumerate}
				\item Two distinct arrows with a common source, or with a common target, have distinct labels.
				\item There is no directed path in the tree whose label lies in $I$.
			\end{enumerate}
			The diagram defines a $\varLambda$-module $X$ with basis $\{x_r\}$ indexed by its vertices.
			For a vertex $r$ with label $v_i$, put $x_re_i=x_r$ and $x_re_j=0$ for $j\ne i$.
			For an arrow $a$ of $Q$, we have $x_ra=x_s$ if the diagram contains an arrow $r\xrightarrow{a}s$, and $x_ra=0$ otherwise.
			We call $X$ a \textbf{tree module} and, by analogy with \cite[Section~1]{Ant07}, call $\{x_r\}$ a \textbf{diagrammatic basis} of $X$.
		\end{definition}

		\begin{remark}\label{lemma of projective injective tree bases}
			Let $\varLambda=kQ/I$ be a finite-dimensional monomial algebra.
			Then every indecomposable projective or injective $\varLambda$-module is a tree module.
		\end{remark}

		\begin{definition}\label{definition of diagrammatic morphism}
			Let $\varLambda=kQ/I$ be a monomial algebra, and let $X$ and $Y$ be tree $\varLambda$-modules with fixed diagrammatic bases $\{x_r\}$ and $\{x_s'\}$.
			Let $D_X$ and $D_Y$ be the corresponding diagrams.
			A $k$-linear map $f:X\to Y$ is called a \textbf{diagrammatic morphism} if there exist nonempty connected subtrees $U\subseteq D_X$ and $V\subseteq D_Y$ and an isomorphism $\varphi:U\to V$ preserving all vertex and arrow labels such that:
			\begin{enumerate}
				\item $U$ is predecessor-closed: every arrow of $D_X$ ending at a vertex of $U$ belongs to $U$.
				\item $V$ is successor-closed: every arrow of $D_Y$ starting at a vertex of $V$ belongs to $V$.
				\item For each vertex $r$ of $D_X$,
				\[
					f(x_r)=
					\begin{cases}
						x'_{\varphi(r)},& r\in U,\\
						0,& r\notin U.
					\end{cases}
				\]
			\end{enumerate}
		\end{definition}

		The first two conditions ensure that $f$ is a $\varLambda$-module homomorphism.

		Krause's theorem gives the following description of the homomorphisms between tree modules; see \cite[Theorem]{Kra91}.
		\begin{theorem}\label{theorem of diagrammatic Hom basis}
			Let $\varLambda=kQ/I$ be a finite-dimensional monomial algebra, and let $X$ and $Y$ be tree $\varLambda$-modules with fixed diagrammatic bases.
			The diagrammatic morphisms $X\to Y$ form a $k$-basis of $\Hom_{\varLambda}(X,Y)$.
		\end{theorem}

	\subsection{Selfinjective-free algebras}\label{section of selfinjective-free algebras}

		Following \cite[Section~3]{Mar90}, an algebra with no nodes and no semisimple summands is called \textbf{selfinjective-free} if every indecomposable projective module belongs to the Nakayama-functor orbit of some noninjective projective module.
		We denote the Nakayama functor by $\nu=D\Hom_{\varLambda}(-,\varLambda)$; in particular, $\nu(e_i\varLambda)\cong D(\varLambda e_i)$.

		For selfinjective-free algebras with no nodes and no semisimple summands, Mart\'inez-Villa constructs a correspondence between isomorphism classes of modules using $\Ext^1$ functors.
		We recall this construction for a stable equivalence $\alpha:\smod\varLambda\to\smod\varLambda'$, where both algebras have no nodes and no semisimple summands and $\varLambda$ is selfinjective-free.
		The following result defines the correspondence on noninjective projective modules.

		\begin{proposition}\textnormal{(\cite[Proposition~1.5]{Mar90})}\label{Mar90 proposition 1.5}
			Let $\varLambda=kQ/I$ and $\varLambda'=kQ'/I'$ be algebras with no nodes and no semisimple summands, and let $\alpha:\smod\varLambda\to\smod\varLambda'$ be a stable equivalence.
			For each indecomposable noninjective projective $\varLambda$-module $P$, there is an indecomposable noninjective projective $\varLambda'$-module $P'$ and a natural isomorphism
			\[
				\Ext^1_{\varLambda}(-,P)\cong\Ext^1_{\varLambda'}(\alpha(-),P').
			\]
			The correspondence $P\mapsto P'$ gives a bijection between the isomorphism classes of indecomposable noninjective projective modules of the two algebras and restricts to a bijection between their simple projective modules.
			It extends by finite direct sums to projective modules without projective-injective summands.
		\end{proposition}

		For each indecomposable nonprojective $\varLambda$-module $M$, take $\alpha'(M)$ to be an indecomposable nonprojective $\varLambda'$-module such that $\alpha'(M)\cong\alpha(M)$ in $\smod\varLambda'$.
		On indecomposable noninjective projective modules, $\alpha'$ is given by Proposition~\ref{Mar90 proposition 1.5}.
		Mart\'inez-Villa extends this correspondence to projective-injective modules along $\nu$-orbits and then to all modules by finite direct sums; see \cite[Proposition~2.4 and Section~3]{Mar90} for more details.
		The resulting bijection on isomorphism classes preserves indecomposable modules, simple modules, projective modules, and injective modules.

		For every $\varLambda$-module $M$ without projective summands, \cite[Proposition~1.4]{Mar90} gives a natural isomorphism
		\[
			\Ext^1_{\varLambda}(-,M)\cong\Ext^1_{\varLambda'}(\alpha(-),\alpha'(M)).
		\]
		Together with Proposition~\ref{Mar90 proposition 1.5} and the preservation of injective modules, this gives the same conclusion for every $\varLambda$-module $M$.

		We shall use the following properties. 

		\begin{lemma}\textnormal{(\cite[Lemma~3.3]{Mar90})}\label{Mar90 lemma 3.3}
			Let $\varLambda=kQ/I$ and $\varLambda'=kQ'/I'$ be algebras with no nodes and no semisimple summands, and suppose that $\varLambda$ is selfinjective-free.
			Let $\alpha:\smod\varLambda\to\smod\varLambda'$ be a stable equivalence, and let $\alpha'$ be its associated correspondence on isomorphism classes of modules.
			Then, for every $\varLambda$-module $X$, there is an isomorphism
			\[
				\alpha'(\rad X)\cong\rad\alpha'(X).
			\]
		\end{lemma}

		\begin{theorem}\textnormal{(\cite[Theorem~3.4]{Mar90})}\label{Mar90 thm 3.4}
			Let $\varLambda=kQ/I$ and $\varLambda'=kQ'/I'$ be algebras with no nodes and no semisimple summands, and suppose that $\varLambda$ is selfinjective-free.
			Let $\alpha:\smod\varLambda\to\smod\varLambda'$ be a stable equivalence, and let $\alpha'$ be its associated correspondence on isomorphism classes of modules.
			Then:
			\begin{enumerate}
				\item[(a)] the quivers $Q$ and $Q'$ are isomorphic;
				\item[(b)] for all $X,Y\in\mod\varLambda$,
				\[
					\dim_k\Hom_{\varLambda}(X,Y)
					=\dim_k\Hom_{\varLambda'}\bigl(\alpha'(X),\alpha'(Y)\bigr).
				\]
				In particular, $\varLambda$ and $\varLambda'$ have the same Cartan matrix.
			\end{enumerate}
		\end{theorem}

		We use the simple-module correspondence to label the vertices of $Q$ and $Q'$ in the same way.
		With this labelling, $S_i'=\alpha'(S_i)$ and $P_i'=\alpha'(P_i)$, since $\alpha'(P_i)$ is the projective cover of $\alpha'(S_i)$.

		\begin{definition}\label{def of uniserial modules}
			A nonzero finite-dimensional module $M$ is \textbf{uniserial} if, for each $i\geq0$, either $\rad^iM=0$ or $\rad^iM/\rad^{i+1}M$ is a simple module.
		\end{definition}

		\begin{lemma}\label{lemma of sending uniserial modules}
			Let $\varLambda=kQ/I$ and $\varLambda'=kQ'/I'$ be algebras with no nodes and no semisimple summands, and suppose that $\varLambda$ is selfinjective-free.
			Let $\alpha:\smod\varLambda\to\smod\varLambda'$ be a stable equivalence, and let $\alpha'$ be its associated correspondence on isomorphism classes of modules.
			Then $\alpha'$ sends uniserial $\varLambda$-modules to uniserial $\varLambda'$-modules.
		\end{lemma}

		\begin{proof}
			Let $M$ be a uniserial $\varLambda$-module, and let $\mathcal S$ be the direct sum of one simple $\varLambda$-module from each isomorphism class.
			Applying Lemma~\ref{Mar90 lemma 3.3} repeatedly gives
			\[
				\alpha'(\rad^iM)\cong\rad^i\alpha'(M)\qquad(i\geq0).
			\]
			For every nonzero $\rad^iM$, by Theorem~\ref{Mar90 thm 3.4}(b) we have
			\[
				\dim_k\Hom_{\varLambda'}
				\bigl(\rad^i\alpha'(M),\alpha'(\mathcal S)\bigr)
				=\dim_k\Hom_{\varLambda}(\rad^iM,\mathcal S)=1.
			\]
			Since $\alpha'(\mathcal S)$ is the direct sum of the simple $\varLambda'$-modules, each nonzero $\rad^i\alpha'(M)$ has simple top.
			Thus $\alpha'(M)$ is uniserial.
		\end{proof}

		We now describe the projective-injective modules of quadratic monomial algebras.

		\begin{proposition}\label{prop of proj-injective modules}
			Let $\varLambda=kQ/I$ be a quadratic monomial algebra with no nodes and no semisimple summands, and let $X$ be an indecomposable projective-injective $\varLambda$-module.
			Then the following hold.
			\begin{enumerate}
				\item $X$ is uniserial.
				\item Write
				\[
					v_0\xrightarrow{a_1}v_1\xrightarrow{a_2}\cdots\xrightarrow{a_n}v_n,\qquad n\geq1,
				\]
				for the diagram of $X$.
				Then in the quiver $Q$ of $\varLambda$, $v_0$ is a source vertex with exactly one outgoing arrow $a_1$, and $v_n$ is a sink vertex with exactly one incoming arrow $a_n$.
			\end{enumerate}
		\end{proposition}

		\begin{proof}
			Write $X=e_0\varLambda$ and $\soc X\cong S_n$.
			Since $\varLambda$ has no semisimple summands, $X$ is not simple.
			Since the relations are monomial, the nonzero paths starting at $v_0$ which cannot be extended on the right by an arrow form a basis of $\soc X$.
			The socle is simple, so there is exactly one such path $p=a_1\cdots a_n$.
			Since $\varLambda$ is finite dimensional, every nonzero path starting at $v_0$ can be extended to a path which cannot be extended further on the right, and this path must be $p$.
			Hence $a_1$ is the only arrow starting at $v_0$, and the paths $e_0,a_1,a_1a_2,\ldots,a_1\cdots a_n$ form a basis of $X$.
			Thus $X$ is uniserial.
			Since $D(X)\cong\varLambda e_n$ is also uniserial, exactly one arrow ends at $v_n$.

			Suppose that $b:v_r\to v_0$ is an arrow with $ba_1\notin I$.
			Since $a_1$ is the only arrow starting at $v_0$ and $I$ is generated by paths of length two, the left multiplication by $b$ defines a monomorphism $b\cdot - \colon e_0\varLambda\to e_r\varLambda$.
			The image of $b\cdot -$ is $b\varLambda$, which is contained in $\rad(e_r\varLambda)$.
			Since $X$ is injective, this monomorphism splits, which contradicts the indecomposability of $e_r\varLambda$.
			Thus $ba_1\in I$ for every arrow $b$ ending at $v_0$.
			Since $a_1$ is the only arrow starting at $v_0$ and $\varLambda$ has no nodes, $v_0$ must be a source vertex.
			The same argument over $\varLambda^{\mathrm{op}}$ shows that $v_n$ is a sink vertex.
		\end{proof}

		\begin{lemma}\label{lemma of AG selfinjective-free}
			Let $\varLambda=kQ/I$ be a quadratic monomial algebra with no nodes and no semisimple summands.
			Then $\varLambda$ is selfinjective-free.
		\end{lemma}

		\begin{proof}
			Let $P$ be an indecomposable projective $\varLambda$-module.
			If $P$ is noninjective, there is nothing to prove.
			Otherwise, Proposition~\ref{prop of proj-injective modules} shows that $\soc P\cong S_n$ is a simple projective module.
			It is noninjective because $\varLambda$ has no semisimple summands, and $\nu(S_n)=D(\varLambda e_n)\cong P$.
			Thus every indecomposable projective module belongs to the $\nu$-orbit of a noninjective projective module.
		\end{proof}

\section{Main result}\label{section of Martinez-Villa question}
	The aim of this section is to prove that quadratic monomial algebras are closed under stable equivalences.
	We first prove an isomorphism result for algebras with no nodes or semisimple summands, treating the case without projective-injective modules first and then using induction on the number of arrows.
	The closure result follows by deleting nodes, removing semisimple summands, and reversing these operations.

	\subsection{The case without projective-injective modules}\label{section of algebra reconstruction}
		We first describe morphisms from injective modules to projective modules of quadratic monomial algebras without projective-injective modules.

		\begin{lemma}\label{lemma of injective to projective}
			Let $\varLambda=kQ/I$ be an algebra with $I$ generated by paths of length two, and suppose that $\varLambda$ has no projective-injective modules.
			Then every morphism from an injective $\varLambda$-module to a projective $\varLambda$-module has semisimple image.
		\end{lemma}
		\begin{proof}
			Let $P_i=e_i\varLambda$ and $I_j$ be indecomposable projective and injective modules with fixed diagrammatic bases, respectively.
			By Theorem~\ref{theorem of diagrammatic Hom basis}, it suffices to show that every diagrammatic morphism $h:I_j\to P_i$ has image in $\soc P_i$.

			Let $U$ and $V$ be the subtrees corresponding to $h$ in the diagrams of $I_j$ and $P_i$, respectively.
			We use the same notation for corresponding vertices and arrows in $U$ and $V$, and denote each arrow or path by its label in $Q$.
			Suppose that $U$ contains an arrow $a$.
			Let $p$ be the path from the source vertex of the diagram of $P_i$ to $t(a)$, and let $q$ be the path from $t(a)$ to the sink vertex of the diagram of $I_j$.
			Since $p$ ends with $a$, $aq\notin I$, and $I$ is generated by paths of length two, we have $pq\notin I$.
			Thus the diagram of $P_i$ contains the path $q$ starting at $t(a)$.
			Since $V$ is successor-closed, this path belongs to $V$, and the corresponding path $q$ belongs to $U$.
			Hence $U$ contains the sink vertex.
			Since $U$ is predecessor-closed, it is the whole diagram of $I_j$.
			Thus $h$ is a monomorphism, which splits since $I_j$ is injective.
			Since $P_i$ is indecomposable, $h$ is an isomorphism, which contradicts the assumption that $P_i$ is not injective.

			Therefore $U$ contains no arrows.
			The subtree $V$ consists of a single vertex and is successor-closed, so the image of $h$ lies in $\soc P_i$.
		\end{proof}

		By Lemma~\ref{lemma of AG selfinjective-free} and Theorem~\ref{Mar90 thm 3.4}(a), the algebras in the following lemma have isomorphic quivers, with vertices corresponding under $\alpha'$.
		We show that the two algebras have the same relations after a suitable choice of arrows, and obtain an algebra isomorphism preserving the corresponding vertex idempotents.

		\begin{lemma}\label{lemma of stable equivalences induces isomorphisms}
			Let $\varLambda=kQ/I$ be a quadratic monomial algebra with no nodes, no semisimple summands, and no projective-injective modules.
			Let $\varLambda'=kQ'/I'$ be an algebra with no nodes and no semisimple summands, and let $\alpha:\smod\varLambda\to\smod\varLambda'$ be a stable equivalence between them.
			Let $\alpha'$ be the correspondence associated with $\alpha$, and denote by $e_i'$ the vertex idempotent of $\varLambda'$ corresponding to the simple module $S_i'=\alpha'(S_i)$.
			Then there is an algebra isomorphism $\varLambda\cong\varLambda'$ sending each $e_i$ to $e_i'$.
		\end{lemma}

		\begin{proof}
			By Lemma~\ref{lemma of AG selfinjective-free}, $\varLambda$ is selfinjective-free.
			Write $\varLambda=\bigoplus_iP_i$ and $\varLambda'=\bigoplus_iP_i'$ for the decompositions of the regular modules, where $P_i=e_i\varLambda$ and $P_i'=e_i'\varLambda'$.
			Let $\cal{F}\subseteq\End_{\varLambda}(\varLambda)$ and $\cal{F}'\subseteq\End_{\varLambda'}(\varLambda')$ be the ideals of endomorphisms factoring through injective modules.
			By definition, $\isEnd_{\varLambda}(\varLambda)=\End_{\varLambda}(\varLambda)/\cal{F}$ and $\isEnd_{\varLambda'}(\varLambda')=\End_{\varLambda'}(\varLambda')/\cal{F}'$.

			We first construct an isomorphism between $\isEnd_{\varLambda}(\varLambda)$ and $\isEnd_{\varLambda'}(\varLambda')$.
			Let \[\tau_\varLambda=D\operatorname{Tr}\colon\smod\varLambda\to\ismod\varLambda\] be the Auslander--Reiten equivalence with quasi-inverse $\tau_\varLambda^{-1}=\operatorname{Tr}D$, and define $\tau_{\varLambda'}$ similarly (see \cite[Chapter~IV, Proposition~1.9]{ARS}).
			The composite
			\[
				\beta=\tau_{\varLambda'}\alpha\tau_\varLambda^{-1}\colon\ismod\varLambda\longrightarrow\ismod\varLambda'
			\]
			is a $k$-linear equivalence.
			For every module $M$ without injective summands, by the Auslander--Reiten formula and the natural isomorphism $\Ext^1_\varLambda(-,M)\cong\Ext^1_{\varLambda'}(\alpha(-),\alpha'(M))$, we have
			\[
				\sHom_{\varLambda'}\bigl(\tau_{\varLambda'}^{-1}\alpha'(M),\alpha(X)\bigr)
				\cong\sHom_\varLambda\bigl(\tau_\varLambda^{-1}M,X\bigr)\cong \sHom_{\varLambda'}\bigl(\alpha\tau_\varLambda^{-1}M,\alpha(X)\bigr),
			\]
			where $X\in \mod\varLambda$.
			These are natural isomorphisms of functors on $\smod\varLambda$.
			By Yoneda's lemma, $\tau_{\varLambda'}^{-1}\alpha'(M)\cong\alpha\tau_\varLambda^{-1}M$ in $\smod\varLambda'$.
			Applying $\tau_{\varLambda'}$ gives $\beta(M)\cong\alpha'(M)$ in $\ismod\varLambda'$.

			By Proposition~\ref{Mar90 proposition 1.5} and Theorem~\ref{Mar90 thm 3.4}(a), $\varLambda'$ has the same numbers of noninjective projectives and of projectives as $\varLambda$.
			Hence $\varLambda'$ also has no projective-injective modules.
			Fix isomorphisms $\beta(P_i)\cong P_i'$ and $\beta(S_i)\cong S_i'$ in $\ismod\varLambda'$.
			The isomorphisms on the projective summands give an algebra isomorphism
			\[
				\varphi\colon\isEnd_{\varLambda}(\varLambda)
				\cong\isEnd_{\varLambda'}(\varLambda').
			\]

			In fact, $\cal{F}$ and $\cal{F}'$ are isomorphic as bimodules with respect to $\varphi$.
			By Lemma~\ref{lemma of injective to projective}, every morphism from an injective $\varLambda$-module to a projective $\varLambda$-module has semisimple image.
			We use this to describe $\cal{F}$.
			Fix projective covers $\pi_i\colon P_i\to S_i$ and $\pi_i'\colon P_i'\to S_i'$, and put
			\[
				T=\{i\mid\Hom_{\varLambda}(I,S_i)\ne0\text{ for some injective }I\}.
			\]
			We regard $\Hom_{\varLambda}(P_i,P_j)$ as the corresponding subspace of $\End_{\varLambda}(\varLambda)$.
			Then
			\begin{equation}\label{eq:reconstruction-F}
				\cal{F}_{ij}\coloneq\cal{F}\cap\Hom_{\varLambda}(P_i,P_j)=
				\begin{cases}
					\Hom_{\varLambda}(S_i,P_j)\pi_i,&i\in T,\\
					0,&i\notin T.
				\end{cases}
			\end{equation}
			For $i\in T$, choose an epimorphism $q:I\to S_i$ with $I$ injective.
			Since $P_i$ is projective, there is a map $f:P_i\to I$ with $qf=\pi_i$, so $t\pi_i=(tq)f$ factors through $I$ for every $t:S_i\to P_j$.
			Conversely, suppose that $gf\ne0$ for maps $f:P_i\to I$ and $g:I\to P_j$ with $I$ injective.
			Since $\Im g$ is semisimple, $gf=t\pi_i$ for some $t:S_i\to P_j$, and $\Im(gf)\cong S_i$ is a direct summand of $\Im g$.
			Hence $\Hom_{\varLambda}(I,S_i)\ne0$ and $i\in T$.

			In the following path calculations, we use the algebra isomorphism $\varLambda\cong\End_{\varLambda}(\varLambda)$ given by left multiplication, and similarly for $\varLambda'$.
			We write $f_x$ for left multiplication by an element $x$ of either algebra.
			Formula~\eqref{eq:reconstruction-F} shows that $\cal{F}$ has a basis consisting of the endomorphisms $f_p$, where $p$ is a nonzero path of positive length ending at a vertex $v_i$ with $i\in T$ and annihilated on the right by every arrow.

			We next show that $\cal{F}$ is generated as an ideal by endomorphisms corresponding to arrows.
			Let $p=a_1\cdots a_m$ be such a path.
			If $a_mb\notin I$ for some arrow $b$, then $pb\notin I$, since $I$ is generated by paths of length two, a contradiction.
			Thus $f_{a_m}\in\cal{F}$.
			Let $\cal{A}$ be the set of arrows $a$ for which $f_a\in\cal{F}$.
			Since $f_p=f_{a_1\cdots a_{m-1}}f_{a_m}$, the endomorphisms $f_a$ with $a\in\cal{A}$ generate $\cal{F}$.

			The endomorphisms $f_p$ corresponding to nonzero paths containing no arrows from $\cal{A}$, including the trivial paths, span a vector space complementary to $\cal{F}$.
			This space is a subalgebra, since $f_pf_q=f_{pq}$ and $pq$ is either zero or a path containing no arrows from $\cal{A}$.
			The quotient map restricts to an algebra isomorphism from this subalgebra to $\isEnd_{\varLambda}(\varLambda)$, whose relations are the original relations among the remaining arrows.
			Using this isomorphism, we write
			\[
				\End_{\varLambda}(\varLambda)=\isEnd_{\varLambda}(\varLambda)\oplus\cal{F}
			\]
			as vector spaces.
			Every map in $\cal{F}$ vanishes on $\rad\varLambda$, and $\cal{F}\subseteq\rad(\End_{\varLambda}(\varLambda))$ because $\varLambda$ has no projective-injective modules.
			Hence
			\[
				\cal{F}\rad(\End_{\varLambda}(\varLambda))=0\qquad \text{and}\qquad \cal{F}^2=0.
			\]

			For $\varLambda$-modules or $\varLambda'$-modules $M,N$, denote the canonical map by
			\[
				\theta_{M,N}\colon\Hom(M,N)\longrightarrow\isHom(M,N).
			\]
			Since $\varLambda$ and $\varLambda'$ have no  projective-injective modules, $\theta_{S,Q}$ is an isomorphism when $S$ is simple and $Q$ is projective.

			By Theorem~\ref{Mar90 thm 3.4}(b), $T$ coincides with the corresponding index set for $\varLambda'$.
			Since $P_i'$ is projective, for every $i\in T$ and every vertex $v_j'$ of $Q'$ we have
			\[
				\Hom_{\varLambda'}(S_i',P_j')\pi_i'\subseteq \cal{F}'\cap\Hom_{\varLambda'}(P_i',P_j').
			\]
			By Theorem~\ref{Mar90 thm 3.4}(b) and the equivalence $\beta$, we have
			\[
				\begin{aligned}
					\dim_k\Hom_{\varLambda}(P_i,P_j)
					&=\dim_k\Hom_{\varLambda'}(P_i',P_j'),\\
					\dim_k\isHom_{\varLambda}(P_i,P_j)
					&=\dim_k\isHom_{\varLambda'}(P_i',P_j').
				\end{aligned}
			\]
			Taking their difference shows that the corresponding blocks $\cal{F}_{ij}$ and $\cal{F}_{ij}'$ have equal dimensions, so \eqref{eq:reconstruction-F} also holds for the corresponding modules of $\varLambda'$.
			Since $\varLambda'$ has no  projective-injective modules, $\cal{F}'\subseteq \rad(\End_{\varLambda'}(\varLambda'))$, and consequently
			\[
				\cal{F}'\rad(\End_{\varLambda'}(\varLambda'))=0,\qquad (\cal{F}')^2=0.
			\]

			For $i\in T$ and $t:S_i\to P_j$, the isomorphism $\theta_{S_i',P_j'}$ gives a unique map $t':S_i'\to P_j'$ whose class is $\beta(\bar t)$ under the fixed isomorphisms.
			Define
			\[
				\psi\colon\cal{F}\longrightarrow\cal{F}',\qquad t\pi_i\longmapsto t'\pi_i'.
			\]
			The map $\psi$ is defined on each block $\cal F_{ij}$ and extended linearly. It is a vector space isomorphism by \eqref{eq:reconstruction-F}.
			By functoriality of $\beta$, we have $\psi(hu)=\varphi(h)\psi(u)$ for $h\in\isEnd_{\varLambda}(\varLambda)$ and $u\in\cal{F}$.
			Since \[\cal{F}\rad(\End_{\varLambda}(\varLambda))=\cal{F}'\rad(\End_{\varLambda'}(\varLambda'))=0,\] the chosen vertex labelling also gives $\psi(uh)=\psi(u)\varphi(h)$.
			Hence $\psi$ is a bimodule isomorphism with respect to $\varphi$.

			The algebra isomorphism $\varLambda\cong\varLambda'$ that we construct must send each $e_i$ to $e_i'$.
			We therefore choose a lift of $\varphi(\overline{f_a})$ in $\End_{\varLambda'}(\varLambda')$ for each arrow $a\notin\cal{A}$.
			For $a:v_i\to v_j$, take $f_a'\in\Hom_{\varLambda'}(P_j',P_i')$ with $\overline{f_a'}=\varphi(\overline{f_a})$, and regard it as an endomorphism of $\varLambda'$.
			Since $\cal{F}'\subseteq\rad(\End_{\varLambda'}(\varLambda'))$, each $f_a'$ belongs to $\rad(\End_{\varLambda'}(\varLambda'))$.

			Let $\iota:X\hookrightarrow Q$ be an inclusion into a projective module, and put $X'=\alpha'(X)$ and $Q'=\alpha'(Q)$.
			The module $X$ has no injective summand; otherwise, such a summand would split off $Q$ and give a projective-injective module.
			Composing with $\iota$ and using the isomorphism $\theta_{S_r,Q}$ shows that $\theta_{S_r,X}$ is an isomorphism for every simple $S_r$.
			Theorem~\ref{Mar90 thm 3.4}(b) and $\beta$ show that $\theta_{S_r',X'}$ is also an isomorphism.
			Choose $\iota':X'\to Q'$ whose class is $\beta(\bar\iota)$.
			By naturality on simple Hom spaces, $\iota'$ is injective on every simple submodule, and hence is a monomorphism.
			Similarly, if $g:P_j\to X$ induces a surjection $\Hom_{\varLambda}(S_r,P_j)\to\Hom_{\varLambda}(S_r,X)$, then every map $g':P_j'\to X'$ whose class is $\beta(\bar g)$ induces a surjection $\Hom_{\varLambda'}(S_r',P_j')\to\Hom_{\varLambda'}(S_r',X')$.

			Fix an arrow $b:v_j\to v_k$ not in $\cal{A}$.
			If $a\in\cal{A}$, then $\psi(f_a)f_b'=0$, since $\psi(f_a)\in\cal{F}'$ and $\cal{F}'\rad(\End_{\varLambda'}(\varLambda'))=0$.
			For every arrow $a\notin\cal{A}$ with $ab\in I$, we need the corresponding lifts to satisfy $f_a'f_b'=0$.
			Put
			\[
				U_b=\{a:v_i\to v_j\mid a\notin\cal{A},\ ab\in I\}.
			\]
			If $U_b$ is empty, take $x_b=0$.
			Otherwise, consider the restrictions of the maps $f_a$ to $P_j$ and write
			\[
				h=(f_a)_{a\in U_b}:P_j\longrightarrow Q_b:=\bigoplus_{\substack{a:v_i\to v_j\\a\in U_b}}P_i,
				\qquad h=\iota g,
			\]
			where $g:P_j\to H=\Im h$ is the quotient map and $\iota:H\hookrightarrow Q_b$ is the inclusion.

			We claim that, for every $r\in T$, the map
			\[
				\Hom_{\varLambda}(S_r,P_j)\longrightarrow\Hom_{\varLambda}(S_r,H),\qquad s\longmapsto gs,
			\]
			is surjective.
			To prove this claim, let $p$ run through the nonzero paths starting at $v_j$.
			The nonzero vectors $h(p)$ have disjoint supports in the bases of $Q_b$ given by paths, and hence form a basis of $H$.
			For every fixed arrow $c$, the nonzero vectors $h(p)c=h(pc)$ also have disjoint supports.
			Thus a linear combination of the $h(p)$ belongs to $\soc H$ only if each term with nonzero coefficient belongs to $\soc H$.
			If $p$ has positive length and $h(p)\ne0$, choose $a\in U_b$ with $ap\ne0$.
			For any arrow $c$ with $pc\ne0$, we have $apc\ne0$, since $I$ is generated by paths of length two.
			Hence $h(p)\in\soc H$ implies $p\in\soc P_j$.
			The trivial path $e_j$ can contribute only to $(\soc H)e_j$.
			If $j\in T$ and $h(e_j)\in\soc H$, every $a\in U_b$ is annihilated on the right by every arrow, so \eqref{eq:reconstruction-F} gives $f_a\in\cal{F}$, contrary to $a\notin\cal{A}$.
			Therefore every element of $(\soc H)e_r$ with $r\in T$ has a preimage in $(\soc P_j)e_r$.
			The claim then follows.

			Put $H'=\alpha'(H)$ and $Q_b'=\alpha'(Q_b)$.
			Choose maps $\iota':H'\to Q_b'$ and $g':P_j'\to H'$ whose classes are $\beta(\bar\iota)$ and $\beta(\bar g)$, respectively.
			By the preceding argument, $\iota'$ is a monomorphism and $g'$ induces a surjection on Hom spaces from $S_r'$ for every $r\in T$.
			Put
			\[
				h'=\iota'g',\qquad h_0=(f_a')_{a\in U_b}:P_j'\to Q_b'.
			\]
			The components of $h'-h_0$ belong to $\cal{F}'$.
			Since $f_b'$ belongs to the radical and $\cal{F}'\rad(\End_{\varLambda'}(\varLambda'))=0$, we have
			\[
				h'f_b'=h_0f_b'.
			\]
			The class of this map in the injectively stable category is zero, since $hf_b=0$.
			If $k\notin T$, formula~\eqref{eq:reconstruction-F} for the corresponding modules of $\varLambda'$ gives $h'f_b'=0$, and we take $x_b=0$.
			If $k\in T$, the same formula writes $h'f_b'$ as a map through $\pi_k'$.
			Its image lies in $\iota'(H')$, so there is a map $t:S_k'\to H'$ such that
			\[
				h'f_b'=\iota't\pi_k'.
			\]
			Choose $s:S_k'\to P_j'$ with $g's=t$, and put $x_b=s\pi_k'\in\cal{F}_{kj}'$.
			Then $h'x_b=h'f_b'$.
			Since $(\cal{F}')^2=0$, we also have $h_0x_b=h'x_b$.
			Consequently, for every $a\in U_b$,
			\[
				f_a'x_b=f_a'f_b'.
			\]
			Replace all lifts simultaneously by $f_b''=f_b'-x_b$.
			Since $\cal{F}'\rad(\End_{\varLambda'}(\varLambda'))=0$ and $(\cal{F}')^2=0$,
			\[
				(f_a'-x_a)(f_b'-x_b)=f_a'f_b'-f_a'x_b=0
			\]
			for every relation $\overline{f_a}\overline{f_b}=0$.
			Sending $\overline{f_{e_i}}$ to $f_{e_i'}$ and $\overline{f_b}$ to $f_b''$ for each arrow $b\notin\cal{A}$ respects all relations of $\isEnd_{\varLambda}(\varLambda)$.
			Hence it defines an algebra homomorphism
			\[
				\rho\colon\isEnd_{\varLambda}(\varLambda)
				\longrightarrow\End_{\varLambda'}(\varLambda')
			\]
			whose composition with the quotient map is $\varphi$.
			In particular, $\rho$ is injective.

			Finally, $\End_{\varLambda'}(\varLambda')=\rho\bigl(\isEnd_{\varLambda}(\varLambda)\bigr)\oplus\cal{F}'$ as vector spaces.
			Since $\psi$ is a bimodule isomorphism and $\cal{F}^2=(\cal{F}')^2=0$, the map
			\[
				\varPhi\colon\End_{\varLambda}(\varLambda)\longrightarrow\End_{\varLambda'}(\varLambda'),\qquad h+u\longmapsto\rho(h)+\psi(u)
			\]
			is an algebra isomorphism.
			Since $\varPhi(f_{e_i})=\rho(\overline{f_{e_i}})=f_{e_i'}$, the resulting isomorphism $\varLambda\cong\varLambda'$ sends each $e_i$ to $e_i'$.
		\end{proof}

	\subsection{Martínez-Villa's question for quadratic monomial algebras}\label{subsection of Martinez-Villa question}
		The next two lemmas will be used in the inductive proof of Theorem~\ref{theorem of quadratic monomial rigidity}.
		By Proposition~\ref{prop of proj-injective modules}, an indecomposable projective-injective module of a quadratic monomial algebra with no nodes or semisimple summands corresponds to a source vertex with only one outgoing arrow.
		For such a module $e_0\varLambda$, put $e=1-e_0$.
		The algebra $e\varLambda e$ is obtained by removing $v_0$ and its outgoing arrow, with the relations among the remaining arrows unchanged.
		The next lemma shows that the corresponding algebras on the two sides are stably equivalent.

		\begin{lemma}\label{lemma of stable equivalence for idempotent subalgebras}
			Let $\varLambda=kQ/I$ be a quadratic monomial algebra with no nodes and no semisimple summands.
			Let $\varLambda'=kQ'/I'$ be an algebra with no nodes and no semisimple summands, and let $\alpha\colon \smod \varLambda\ra\smod\varLambda'$ and $\beta\colon \smod \varLambda'\ra \smod\varLambda$ be a stable equivalence and its quasi-inverse, respectively.
			Let $v_0$ be a vertex of $Q$ such that $e_0\varLambda$ is an indecomposable projective-injective module.
			Put $e=1-e_0$.
			Then there is a primitive idempotent $e_0'$ of $\varLambda'$ such that, with $e'=1-e_0'$, the algebras $e\varLambda e$ and $e'\varLambda'e'$ are stably equivalent.
		\end{lemma}

		\begin{proof}
			Let $S_0=\top(e_0\varLambda)$.
			By Proposition~\ref{prop of proj-injective modules}, $v_0$ is a source vertex.
			Put $S_0'=\alpha'(S_0)$, so $\alpha(S_0)\cong S_0'$ in $\smod\varLambda'$.
			By Theorem~\ref{Mar90 thm 3.4}, the corresponding vertex $v_0'$ is also a source vertex.
			Put $e'=1-e_0'$, $\varGamma=e\varLambda e$, and $\varGamma'=e'\varLambda'e'$.
			Since $v_0$ is a source vertex, $e\varLambda e_0=0$ and $\varGamma\cong \varLambda/\varLambda e_0\varLambda$.

			Let $\tau\colon\mod \varGamma\to\mod\varLambda$ be the functor induced by the quotient homomorphism $\varLambda\to \varGamma$, and define $r\colon\mod\varLambda\to\mod \varGamma$ by
			\[
				r(M)=M\otimes_{\varLambda} \varGamma=M/(Me_0\varLambda).
			\]
			We have $r\tau\cong\id$.
			Moreover,
			\[
				r(e_0\varLambda)=0,\qquad r(e_i\varLambda)=e_i\varGamma\ (i\ne0),\qquad
				\tau(e_i\varGamma)=e_i\varLambda\ (i\ne0).
			\]
			Thus $r$ and $\tau$ preserve projectives and induce functors on stable categories, which we denote by the same symbols by abuse of notation.

			The functor $\tau$ is fully faithful both on module categories and on stable categories.

			To obtain a stable equivalence between $\varGamma$ and $\varGamma'$, we show that $\tau$ induces an equivalence between $\smod \varGamma$ and the full subcategory
			\[
				\mathcal C_0=\{M\in\smod \varLambda\mid\sHom_{\varLambda}(M,S_0)=0\}.
			\]
			Let $M$ be a $\varLambda$-module without projective summands.
			A nonzero map $M\to S_0$ cannot factor through a projective module: otherwise, a nonzero component $M\to e_0\varLambda\to S_0$ would make $M\to e_0\varLambda$ surjective by Nakayama's lemma, and this epimorphism would split since $e_0\varLambda$ is projective.
			Thus $\sHom_{\varLambda}(M,S_0)=\Hom_{\varLambda}(M,S_0)$.
			Since $v_0$ is a source vertex, $(\rad M)e_0=0$, so $S_0$ occurs in $\top M$ if and only if $Me_0\ne0$.
			It follows that $M$ belongs to $\mathcal C_0$ if and only if $Me_0=0$, which is precisely the condition that its $\varLambda$-action factors through $\varGamma$.

			Define $r'$, $\tau'$, and $\mathcal C_0'$ in the same way for $\varLambda'$.
			Since $\alpha(S_0)\cong S_0'$, the equivalence $\alpha$ restricts to an equivalence between $\mathcal C_0$ and $\mathcal C_0'$.
			Hence
			\[
				\widetilde\alpha=r'\alpha\tau\colon\smod \varGamma\to\smod \varGamma',
				\qquad
				\widetilde\beta=r\beta\tau'\colon\smod \varGamma'\to\smod \varGamma
			\]
			are quasi-inverse equivalences.
		\end{proof}

		The following lemma expresses $\varLambda$ as a one-point extension of $e\varLambda e$.
		An isomorphism between $e\varLambda e$ and $e'\varLambda'e'$ extends to the original algebras if it preserves the vertex idempotents associated with the socles of the corresponding projective-injective modules.

		\begin{lemma}\label{lemma of isomorphisms of one-point extensions}
			Let $\varLambda=kQ/I$, and let $v_0$ be a source vertex of $Q$ such that $P_0=e_0\varLambda$ is a non-simple indecomposable projective-injective module.
			Put $e=1-e_0$, $\varGamma=e\varLambda e$, and $M=e_0\varLambda e$.
			If $\soc P_0=S_t$, then $t\ne0$ and
			\[
				M\cong D(\varGamma e_t),\qquad
				\varLambda\cong\begin{pmatrix}k&D(\varGamma e_t)\\0&\varGamma\end{pmatrix}.
			\]

			Let $\varLambda'=kQ'/I'$ be another finite-dimensional basic algebra over $k$, and let $v_0'$ be a source vertex of $Q'$ such that $P_0'=e_0'\varLambda'$ is a non-simple indecomposable projective-injective module with socle $S_{t'}'$.
			Put $e'=1-e_0'$, $\varGamma'=e'\varLambda'e'$, and $M'=e_0'\varLambda'e'$.
			Every $k$-algebra isomorphism $\varphi:\varGamma\to \varGamma'$ satisfying $\varphi(e_t)=e_{t'}'$ extends to an isomorphism $\varPhi:\varLambda\to \varLambda'$ which sends $e_0$ to $e_0'$.
		\end{lemma}

		\begin{proof}
			Since $v_0$ is a source vertex, $e\varLambda e_0=0$ and $e_0\varLambda e_0=k$, so
			\[ 
				\varLambda\cong\begin{pmatrix}k&M\\0&\varGamma\end{pmatrix}. 
			\]
			The simple module $S_0$ is injective, so $t\ne0$ because $P_0$ is not simple.
			Since $P_0\cong D(\varLambda e_t)$, restriction of linear functionals gives
			\[
				M=P_0e\cong D(\varLambda e_t)e\cong D(e\varLambda e_t)=D(\varGamma e_t).
			\]

			Use these isomorphisms to write $M=D(\varGamma e_t)$ and $M'=D(\varGamma'e_{t'}')$.
			Since $\varphi(e_t)=e_{t'}'$, composition with $\varphi^{-1}$ defines a bijection
			\[
				\psi:M\to M',\qquad
				\psi(f)(x')=f(\varphi^{-1}(x')),
				\qquad \psi(mb)=\psi(m)\varphi(b).
			\]
			Consequently
			\[
				\varPhi\left(\begin{pmatrix}\lambda&m\\0&b\end{pmatrix}\right)
				=\begin{pmatrix}\lambda&\psi(m)\\0&\varphi(b)\end{pmatrix}
			\]
			is the required algebra isomorphism.
		\end{proof}

		At the end of \cite[Section~3]{Mar90}, Mart\'inez-Villa asked whether two nonisomorphic algebras satisfying the assumptions of Theorem~\ref{Mar90 thm 3.4} can be stably equivalent.
		We state the question in the following form.

		\begin{question}\label{question of Martinez-Villa}
			Let $\varLambda$ be a selfinjective-free algebra and let $\varLambda'$ be an algebra with no nodes and no semisimple summands.
			If $\varLambda$ and $\varLambda'$ are stably equivalent, are they isomorphic?
		\end{question}

		The preceding lemmas give the following result by induction.

		\begin{theorem}\label{theorem of quadratic monomial rigidity}\label{theorem of AG rigidity}\label{proposition of conditional labelled induction}
			Let $\varLambda=kQ/I$ be a quadratic monomial algebra with no nodes and no semisimple summands.
			Let $\varLambda'=kQ'/I'$ be an algebra with no nodes and no semisimple summands, and let $\alpha:\smod\varLambda\to\smod\varLambda'$ be a stable equivalence.
			Then there is an algebra isomorphism $\varLambda\cong\varLambda'$ preserving the vertex correspondence given by $\alpha'$.
		\end{theorem}
		\begin{proof}
			Label the vertices so that $\alpha'(S_i)\cong S_i'$.
			We prove by induction on the number of arrows that there is an algebra isomorphism sending each $e_i$ to $e_i'$.
			Assume that the assertion holds whenever the quiver of the quadratic monomial algebra has fewer arrows than $Q$.
			If $\varLambda$ has no  projective-injective modules, this follows from Lemma~\ref{lemma of stable equivalences induces isomorphisms}.

			Suppose that $P_0=e_0\varLambda$ is an indecomposable projective-injective module with socle $S_t$.
			By Proposition~\ref{prop of proj-injective modules} and Theorem~\ref{Mar90 thm 3.4}, $v_0$ and $v_0'$ are source vertices with exactly one outgoing arrow, and $P_0'=e_0'\varLambda'$ is projective-injective with socle $S_t'$.
			Put $e=1-e_0$, $e'=1-e_0'$, $\varGamma=e\varLambda e$, and $\varGamma'=e'\varLambda'e'$.
			Deleting $v_0$ and $v_0'$ from $Q$ and $Q'$ gives the isomorphic quivers of $\varGamma$ and $\varGamma'$, respectively.

			Since $v_0$ and $v_0'$ are source vertices, $e_j\varGamma=e_j\varLambda$ and $e_j'\varGamma'=e_j'\varLambda'$ for $j\ne0$.
			For each nonsimple $P_j$ and each surviving vertex $v_r$, the multiplicity of $S_r$ as a composition factor of $\rad P_j/\soc P_j$ is 
			\[\dim_k\Hom_{\varLambda}(P_r,\rad P_j)-\dim_k\Hom_{\varLambda}(S_r,P_j).\]
			By Lemma~\ref{Mar90 lemma 3.3} and Theorem~\ref{Mar90 thm 3.4}, these multiplicities are equal for the corresponding modules.
			Hence the nodes of $\varGamma$ and $\varGamma'$ occur at corresponding vertices.

			Let $\gamma=r'\alpha\tau:\smod\varGamma\to\smod\varGamma'$ be the stable equivalence constructed in Lemma~\ref{lemma of stable equivalence for idempotent subalgebras}.
			Then $\tau'\gamma\cong\alpha\tau$, and $\gamma(S_i)\cong S_i'$ for every nonprojective simple $\varGamma$-module $S_i$.
			Since $\tau$ and $\tau'$ are exact and preserve projectives, we obtain a natural isomorphism
			\begin{equation}\label{eq simple Ext correspondence for idempotent subalgebras}
				\Ext^1_\varGamma(-,S_i)
				\cong\Ext^1_{\varGamma'}(\gamma(-),S_i')
			\end{equation}
			for each surviving vertex $v_i$.

			Let $F:\mod\varGamma\to\mod \widehat{\varGamma}$ and $F':\mod\varGamma'\to\mod \widehat{\varGamma}'$ be the node-deletion functors.
			Put $i^*=i$ if $v_i$ is not a node and $i^*=i^-$ otherwise.
			Write $I_i=D(\varGamma e_i)$, $\widehat{I}_i=D(\widehat{\varGamma}e_{i^*})$, and let $p_i:I_i\to I_i/S_i$ be the quotient map.
			From the construction of $F$, we can check that
			\[
				F(I_i)\cong \widehat{I}_i,\qquad
				F(I_i/S_i)\cong \widehat{I}_i/S_{i^*},
			\]
			and that $F(p_i)$ is the corresponding quotient map.
			By Proposition~\ref{proposition of node deletion}, $F$ is full, and every map to $I_i/S_i$ sent to zero by $F$ factors through a projective module and hence lifts through $p_i$.
			Thus $F$ induces an isomorphism between the cokernels computing the two Ext functors, giving a natural isomorphism
			\begin{equation}\label{eq node splitting Ext}
				\Ext^1_\varGamma(-,S_i)
				\cong\Ext^1_{\widehat{\varGamma}}(F(-),S_{i^*}).
			\end{equation}
			The same argument applies to $F'$.

			The node-deletion functors induce a stable equivalence $\delta:\smod \widehat{\varGamma}\to\smod \widehat{\varGamma}'$ such that $\delta\underline F\cong\underline F'\gamma$.
			Combining \eqref{eq simple Ext correspondence for idempotent subalgebras} and \eqref{eq node splitting Ext} gives
			\begin{equation}\label{eq split simple Ext correspondence}
				\Ext^1_{\widehat{\varGamma}}(-,S_{i^*})
				\cong\Ext^1_{\widehat{\varGamma}'}(\delta(-),S_{i^*}').
			\end{equation}
			For a nonprojective simple $S_i$, we have $F(S_i)=S_i$ if $v_i$ is not a node and $F(S_i)=S_{i^+}$ otherwise.
			Together with $\gamma(S_i)\cong S_i'$, this shows that $\delta$ preserves the labels of the nonprojective simple modules.

			Remove the corresponding isolated vertices of $\widehat{\varGamma}$ and $\widehat{\varGamma}'$.
			The remaining algebra on the quadratic monomial side is quadratic monomial by Remark~\ref{remark of quadratic monomial node type}, and is selfinjective-free by Lemma~\ref{lemma of AG selfinjective-free}.
			Every remaining projective simple module is noninjective and has the form $S_{i^*}$.
			By \eqref{eq split simple Ext correspondence} and Proposition~\ref{Mar90 proposition 1.5}, its corresponding module is $S_{i^*}'$.
			Hence the induced correspondence preserves all remaining vertex labels.

			The quiver of $e\varLambda e$ has one fewer arrow than $Q$, whereas node deletion and removal of isolated vertices do not change the number of arrows.
			If no vertices remain, the required isomorphism is immediate; otherwise, the induction hypothesis gives an isomorphism preserving the corresponding vertex idempotents.
			Adjoining the corresponding isolated vertices gives such an isomorphism $\widehat{\varGamma}\cong \widehat{\varGamma}'$.
			We replace each pair $e_{j^+},e_{j^-}$ by their sum and retain the radical.
			Since the isomorphism sends each vertex idempotent to the corresponding one, it restricts to an isomorphism $\varphi:\varGamma\to\varGamma'$ with $\varphi(e_i)=e_i'$.
			In particular, $\varphi(e_t)=e_t'$, so Lemma~\ref{lemma of isomorphisms of one-point extensions} extends $\varphi$ to the required isomorphism $\varLambda\cong\varLambda'$.
		\end{proof}

		\begin{corollary}\label{corollary of quadratic monomial Martinez-Villa question}
			Question~\ref{question of Martinez-Villa} has a positive answer for quadratic monomial algebras with no nodes and no semisimple summands.
		\end{corollary}
		\begin{proof}
			By Lemma~\ref{lemma of AG selfinjective-free}, these algebras are selfinjective-free, and Theorem~\ref{theorem of quadratic monomial rigidity} gives the required isomorphism.
		\end{proof}

	\subsection{Closure under stable equivalences}\label{sec:induction}
		Our main result now follows from Theorem~\ref{theorem of quadratic monomial rigidity} by deleting nodes and removing semisimple summands.

		\begin{theorem}\label{theorem of quadratic monomial stable closure}
			Let $\varLambda$ be a quadratic monomial algebra and let $\varLambda'$ be a basic algebra stably equivalent to $\varLambda$.
			Then $\varLambda'$ is quadratic monomial.
			Moreover, the algebras obtained from $\varLambda$ and $\varLambda'$ by deleting all nodes and removing all semisimple direct summands are isomorphic.
		\end{theorem}

		\begin{proof}
			Let $\varLambda=kQ/I$ be quadratic monomial, and let $\varLambda'=kQ'/I'$ be a basic algebra stably equivalent to $\varLambda$.
			Delete all nodes on both sides and remove the semisimple direct summands, obtaining stably equivalent algebras $\widehat{\varLambda},\widehat{\varLambda}'$ with no nodes or semisimple summands.
			By Remark~\ref{remark of quadratic monomial node type}, $\widehat{\varLambda}$ is quadratic monomial if it is nonzero.
			If $\widehat{\varLambda}=0$, then $\smod\widehat{\varLambda}'=0$, so $\widehat{\varLambda}'$ is semisimple and hence zero.
			Otherwise, applying Theorem~\ref{theorem of quadratic monomial rigidity} to $\widehat{\varLambda}$ and $\widehat{\varLambda}'$ gives $\widehat{\varLambda}\cong\widehat{\varLambda}'$.
			Restore the semisimple factors and reverse the node deletions on the $\varLambda'$ side.
			These operations preserve the quadratic monomial property, so $\varLambda'$ is quadratic monomial.
		\end{proof}

		Thus, up to algebra isomorphism and semisimple summands, the algebras in Theorem~\ref{theorem of quadratic monomial stable closure} are related by stable equivalences of node type.

\section{Simple-minded systems and stable equivalences of Morita type}\label{section of sms and Morita type}

	% We prove that a quadratic monomial algebra has a unique weakly simple-minded system, consisting of its nonprojective simple modules.
	% This uniqueness will be used to prove that algebras stably equivalent of Morita type have isomorphic nonsemisimple factors whenever one of them is quadratic monomial.

	\subsection{Definitions and questions}

		We use the definitions of Koenig--Liu for right modules; equivalently, we apply their definitions to the opposite algebra.
		For an algebra $A$, let $\mathcal S_A$ be a set of representatives of the isomorphism classes of nonprojective simple right $A$-modules.
		Systems of modules are understood up to isomorphism.

		\begin{definition}[{\cite[Definitions~2.1 and~5.3]{KL12}}]\label{definition of weak sms}
			An \textbf{orthogonal system} in $\smod A$ is a set $\mathcal S$ of pairwise nonisomorphic indecomposable nonprojective modules such that
			\[
				\sHom_A(S,T)=0\quad(S\ne T),\qquad
				\underline{\End}_A(S)\text{ is a division ring}\quad(S\in\mathcal S).
			\]
			It is a \textbf{weakly simple-minded system} if, for every indecomposable nonprojective module $X$, there is an $S\in\mathcal S$ with $\sHom_A(X,S)\ne0$.
		\end{definition}

		To recall the generating condition for a simple-minded system over an arbitrary algebra, let $\mathsf e(A)$ be the additive subcategory generated by the node simples and the indecomposable projective noninjective modules.

		\begin{definition}[{\cite[Definition~2.1]{KL12}}]\label{definition of sms}
			For a set $\mathcal S$ of modules, put $\langle\mathcal S\rangle_1=\add\mathcal S$.
			For $n\geq2$, let $\langle\mathcal S\rangle_n$ be the additive subcategory generated by the indecomposable modules $Y$ for which there is an exact sequence
			\[
				0\longrightarrow X\longrightarrow Y\oplus P\longrightarrow Z\longrightarrow0
			\]
			with $P$ projective, $X\in\add(\langle\mathcal S\rangle_{n-1}\cup\mathsf e(A))$, and $Z\in\add\mathcal S$. An orthogonal system $\mathcal S$ is a \textbf{simple-minded system} if every indecomposable nonprojective module belongs to $\langle\mathcal S\rangle_n$ for some $n\geq1$.
		\end{definition}

		The set $\mathcal S_A$ is a simple-minded system, and every simple-minded system is weakly simple-minded; see \cite[Lemma~2.3 and Remark~5.4]{KL12}. Koenig--Liu posed the following questions.

		\begin{question}[{\cite[Remark~5.4]{KL12}}]\label{question of weak sms generation}
			Is every weakly simple-minded system in $\smod A$ a simple-minded system?
		\end{question}

		The answer is negative in general: Zhang constructed infinite weakly simple-minded systems which are not simple-minded systems over a $2$-domestic Brauer graph algebra \cite{Zha24BGA}. On the other hand, Koenig--Liu proved that every weakly simple-minded system over a representation-finite selfinjective algebra is a simple-minded system \cite[Theorem~5.6]{KL12}. For domestic Brauer graph algebras, Zhang proved that every finite weakly simple-minded system is a simple-minded system \cite[Corollary~3.5]{Zha26BGAII}.
		Zhang further conjectured that every finite weakly simple-minded system (over a selfinjective algebra) is a simple-minded system \cite[Conjecture~5.18]{Zha25}.

		\begin{question}[{\cite[Question~6.1]{KL12}}]\label{question of sms cardinality}
			Does every simple-minded system in $\smod A$ have the same cardinality as $\mathcal S_A$?
		\end{question}

		This question has a  positive answer for triangular algebras \cite{KL12}, representation-finite selfinjective algebras \cite{CKL15} and for domestic Brauer graph algebras \cite{Zha25,Zha26BGAII}. It is closely related to the \textbf{Auslander--Reiten conjecture} that stable equivalences preserve the number of isomorphism classes of nonprojective simple modules.
		Indeed, by \cite[Theorem~3.2]{KL12}, a stable equivalence $\alpha:\smod B\longrightarrow\smod A$ sends $\mathcal S_B$ to a simple-minded system in $\smod A$. A positive answer to Question~\ref{question of sms cardinality} for $A$ would therefore give
		\[
			|\mathcal S_B|
			=|\alpha(\mathcal S_B)|
			=|\mathcal S_A|.
		\]
		Thus a positive answer for all algebras would imply the Auslander--Reiten conjecture.

		Meanwhile, Chan--Koenig--Liu formulate the simple-image problem for selfinjective algebras in \cite{CKL15}. For arbitrary algebras, its corresponding formulation is as follows.

		\begin{question}[{cf.~\cite[Problem~2.2]{CKL15}}]\label{question of sms simple images}
			Given a simple-minded system $\mathcal S$ in $\smod A$, do there exist an algebra $B$ and a stable equivalence $\alpha:\smod B\longrightarrow\smod A$ such that $\alpha(\mathcal S_B)=\mathcal S$? Can $\alpha$ be chosen to be of Morita type?
		\end{question}

		We also note that both questions in Question~\ref{question of sms simple images} have positive answers for triangular algebras \cite{KL12} and for representation-finite selfinjective algebras \cite{CKL15}.

	\subsection{Simple-minded systems over quadratic monomial algebras}

		\begin{definition}[{see \cite[Introduction]{Rin12}}]
			A module $M\in\mod A$ is called \textbf{torsionless} if it is isomorphic to a submodule of a finitely generated projective module.
		\end{definition}

		The following lemma holds for any finite-dimensional algebra.

		\begin{lemma}\label{lemma of weak sms torsionless}
			Every nonsimple member of a weakly simple-minded system is torsionless.
		\end{lemma}

		\begin{proof}
			Let $\mathcal S$ be a weakly simple-minded system in $\smod A$, let $M\in\mathcal S$ be nonsimple, and put
			\[
				R=\bigcap_{f\in\Hom_A(M,A)}\ker f.
			\]
			We prove that $R=0$ by contradiction. Suppose that $R\ne0$. Choose a simple submodule $U\subseteq R$ and let $\pi:M\twoheadrightarrow N=M/U$ be the quotient map. Since $M$ is nonsimple, $U$ is proper. Moreover, $N$ is not projective.
			By the definition of $R$, every map from $M$ to a projective module vanishes on $U$. It follows that, for every module $T$, $\pi$ induces an injection
			\[
				\sHom_A(N,T)\longrightarrow\sHom_A(M,T),
				\qquad \underline g\longmapsto\underline{g\pi}.
			\]
			Indeed, if $g\pi=vu$ with $u:M\to P$, $v:P\to T$ and $P$ projective, then $u=\bar u\pi$ for some $\bar u:N\to P$. Since $\pi$ is surjective, $g=v\bar u$, so $g$ also factors through a projective module.

			For $T\in\mathcal S\setminus\{M\}$, the target is zero by orthogonality. For $T=M$, every composite $g\pi$ is noninvertible because it vanishes on the nonzero submodule $U$. Since $\End_A(M)$ is local and $\underline{\End}_A(M)$ is a division ring, every noninvertible endomorphism of $M$ has zero stable class.
			Hence $\sHom_A(N,T)=0$ for every $T\in\mathcal S$. This contradicts weak generation applied to a nonprojective indecomposable summand of $N$. Therefore $R=0$.

			Choose a $k$-basis $f_1,\ldots,f_d$ of $\Hom_A(M,A)$. The map
			\[
				M\longrightarrow A^d,\qquad
				m\longmapsto(f_1(m),\ldots,f_d(m))
			\]
			has kernel $R=0$, so $M$ is torsionless.
		\end{proof}

		We now turn to quadratic monomial algebras.

		\begin{lemma}\label{lemma of torsionless simple maps}
			Let $A=kQ/I$ be a quadratic monomial algebra and let $M$ be an indecomposable nonsimple torsionless module. Then $\sHom_A(S_i,M)=0$ for every simple module $S_i$.
		\end{lemma}

		\begin{proof}
			Fix an embedding $M\subseteq P$ with $P$ projective. Since $M$ is indecomposable and nonsimple, $\soc M\subseteq\rad M=M\rad A$. Note that the arrows of $Q$ generate $\rad A$ as a left ideal, so we may write any $u\in(\soc M)e_i$ as
\[
    u=\sum_{\substack{a\in Q_1\\a:v_r\to v_i}}m_a a,
    \qquad m_a\in Me_r.
\]

			Choose the path basis in each indecomposable projective summand of $P$. The nonzero paths with different final arrows span independent subspaces, and $\soc P$ is spanned by paths admitting no nonzero continuation on the right.
			Since the relations are quadratic and monomial, a nonzero path ending in an arrow $a$ lies in $\soc P$ if and only if $a\rad A=0$. Thus $u\in\soc P$ implies that $m_a a=0$ whenever $a\rad A\ne0$.

			For each arrow $a:v_r\to v_i$ with $a\rad A=0$, the map $S_i\to M$ sending $\overline e_i$ to $m_a a$ factors through the projective module $e_rA$:
			\[
				S_i\longrightarrow e_rA\longrightarrow M,
				\qquad
				\overline e_i\longmapsto a,\qquad p\longmapsto m_a p.
			\]
			Summing these factorizations shows that the map $S_i\to M$ sending $\overline e_i$ to $u$ factors through a projective module. Since every map $S_i\to M$ has this form, $\sHom_A(S_i,M)=0$.
		\end{proof}

		\begin{theorem}\label{theorem of quadratic monomial sms}
			Let $A$ be a quadratic monomial algebra. Then $\mathcal S_A$ is the unique weakly simple-minded system in $\smod A$. In particular, it is the unique simple-minded system.
		\end{theorem}

		\begin{proof}
			Let $\mathcal S$ be a weakly simple-minded system and let $S_i\in\mathcal S_A$. Weak generation gives a $T\in\mathcal S$ with $\sHom_A(S_i,T)\ne0$. If $T$ were nonsimple, Lemma~\ref{lemma of weak sms torsionless} would make it torsionless, contradicting Lemma~\ref{lemma of torsionless simple maps}. Thus $T$ is simple and $T\cong S_i$. Hence $\mathcal S_A\subseteq\mathcal S$.

			If $T\in\mathcal S\setminus\mathcal S_A$, weak generation for the standard system $\mathcal S_A$ gives an $S_i\in\mathcal S_A$ with $\sHom_A(T,S_i)\ne0$, contradicting orthogonality in $\mathcal S$. Therefore $\mathcal S=\mathcal S_A$.
		\end{proof}

		The preceding theorem gives the following corollary.

		\begin{corollary}\label{corollary of sms questions}
			Questions~\ref{question of weak sms generation}, \ref{question of sms cardinality}, and~\ref{question of sms simple images} have positive answers for quadratic monomial algebras.
		\end{corollary}

		The uniqueness result also shows that stable equivalences of quadratic monomial algebras preserve nonprojective simple modules.

		\begin{lemma}\label{lemma of AG simple preservation}
			Let $A$ be a quadratic monomial algebra and let $\alpha:\smod A\to\smod B$ be a stable equivalence. Then $\alpha$ induces a bijection between the isomorphism classes of nonprojective simple modules.
		\end{lemma}

		\begin{proof}
			Note that the defining conditions of a weakly simple-minded system are also preserved by equivalences of the stable categories. Thus $\alpha^{-1}(\mathcal S_B)$ is a weakly simple-minded system in $\smod A$. Theorem~\ref{theorem of quadratic monomial sms} gives $\alpha^{-1}(\mathcal S_B)=\mathcal S_A$.
		\end{proof}

	\subsection{Stable equivalences of Morita type}

		\begin{definition}[{see \cite[Section~5.3]{Zim14}}]
			\label{definition of stable equivalence of Morita type}
			Bimodules ${}_AM_B$ and ${}_BN_A$ define a \textbf{stable equivalence of Morita type} if they are finitely generated and projective on both sides and
			\[
				M\otimes_BN\cong A\oplus P,\qquad
				N\otimes_AM\cong B\oplus Q,
			\]
			where $P$ and $Q$ are projective as $A$-bimodules and $B$-bimodules, respectively.
			The induced equivalence $\smod A\to\smod B$ is given by $-\otimes_AM$, with inverse $-\otimes_BN$.
		\end{definition}

		% We note that Liu proved that a stable equivalence of Morita type is a Morita equivalence if its defining tensor functor sends every simple module to a simple module \cite{Liu03}. We use the following version of Hu--Xi, which only requires
		% preservation of nonprojective simple modules in the stable category.

		% \begin{lemma}[{cf.~\cite[Proposition~3.3]{HX18}}]
		% \label{lemma of simple preserving Morita type}
		% Let $A$ and $B$ be algebras without semisimple direct summands and let $\alpha:\smod A\longrightarrow\smod B$ be a stable equivalence of Morita type.
		% If $\alpha(S)$ is isomorphic in $\smod B$ to a nonprojective simple $B$-module for every nonprojective simple $A$-module $S$, then $\alpha$ lifts to a Morita equivalence.
		% \end{lemma}

		Liu proved that a stable equivalence of Morita type is a Morita equivalence if its defining tensor functor sends every simple module to a simple module \cite[Theorem~1.1]{Liu03}. Together with \cite[Lemma~3.1]{Liu06} and \cite[Lemma~4.4]{KL12}, this yields the following criterion; see also the proof of \cite[Proposition~4.3]{KL12}.

		\begin{lemma}
			\label{lemma of simple preserving Morita type}
			Let $A$ and $B$ be algebras without semisimple direct summands, and let $\alpha:\smod A\longrightarrow\smod B$ be a stable equivalence of Morita type.
			If $\alpha(S)$ is isomorphic in $\smod B$ to a nonprojective simple $B$-module for every nonprojective simple $A$-module $S$, then $\alpha$ lifts to a Morita equivalence.
		\end{lemma}

		We shall need the following result, which is well known to experts. For completeness, we recall a bimodule construction from unpublished notes of Liu.

		\begin{lemma}\label{lemma of adjoining semisimple factors}
			The algebras $A$ and $B:=A\times k$ are stably equivalent of Morita type.
		\end{lemma}

		\begin{proof}
			The required stable equivalence is defined by the bimodules
			\[
				{}_AM_B={}_A(A\oplus A)_B,
				\qquad
				{}_BN_A={}_B(A\oplus A)_A,
			\]
			with actions $a(x,y)(b,\lambda)=(axb,ay\lambda)$, and $(b,\lambda)(x,y)a=(bxa,\lambda ya)$.
		\end{proof}

		We now prove the main result of this subsection.

		\begin{theorem}\label{corollary of AG Morita rigidity}
			Let $A$ be a quadratic monomial algebra and let $B$ be an algebra stably equivalent to $A$ of Morita type. Then the following hold.
			\begin{enumerate}
				\item After removing their semisimple direct summands, $A$ and $B$ are isomorphic.
				If $A$ and $B$ have no semisimple direct summands, every stable equivalence of Morita type between them lifts to a Morita equivalence. 

				\item If $A$ and $B$ are also stably equivalent of node type, then $A\cong B$.
			\end{enumerate}
		\end{theorem}

		\begin{proof}
			\underline{(1)}
			If $A$ is semisimple, then $\smod B=0$, so $B$ is also semisimple. Otherwise, Lemma~\ref{lemma of adjoining semisimple factors} allows us to remove the semisimple direct summands of $A$ and $B$ while preserving stable equivalence of Morita type. We may therefore assume that neither algebra has a semisimple direct summand.

			Let $\alpha:\smod A\longrightarrow\smod B$ be a stable equivalence of Morita type.
			By Lemma~\ref{lemma of AG simple preservation}, $\alpha$ sends nonprojective simple modules to nonprojective simple modules in the stable category. Thus by Lemma~\ref{lemma of simple preserving Morita type}, $\alpha$ lifts to a Morita equivalence. Since $A$ and $B$ are basic, $A\cong B$.

			\medskip
			\noindent\underline{(2)}
			Note that stable equivalences of node type preserve the number of semisimple direct factors. Together with (1), this gives $A\cong B$.
		\end{proof}

		\begin{remark}\label{remark of node and Morita type}
			Let $A$ and $B$ be algebras without semisimple direct summands. If a stable equivalence $\alpha:\smod A\longrightarrow\smod B$ is both of node type and of Morita type, then $\alpha$ lifts to a Morita equivalence. Indeed, node-deletion equivalences and their inverses preserve nonprojective simple modules in the stable category, so Lemma~\ref{lemma of simple preserving Morita type} applies.
			Here both conditions concern the same stable functor. For quadratic monomial algebras, Lemma~\ref{lemma of AG simple preservation} applies to every stable equivalence, so the two equivalences in Theorem~\ref{corollary of AG Morita rigidity}(2) may be chosen independently.
		\end{remark}

\section{Gorenstein-projective modules and singularity categories}\label{section of Gorenstein and singularity categories}

	We now turn to the triangulated categories associated to quadratic monomial algebras. We first recall the general definitions and then use relation quivers to describe these categories and the effect of node deletion.

	\subsection{Comparison functors}\label{subsection of stable and singularity categories}

		We use cochain complexes, with differentials of degree $1$, and write $\operatorname{proj}\varLambda$ for the category of finitely generated projective right $\varLambda$-modules.

		A $\varLambda$-module $G$ is called \textbf{Gorenstein-projective} if $G\cong\ker d^0$ for an acyclic complex
		\[
			\cdots\longrightarrow P^{-1}\xrightarrow{d^{-1}}P^0\xrightarrow{d^0}P^1\longrightarrow\cdots
		\]
		of finitely generated projective modules such that $\Hom_{\varLambda}(P^\bullet,\varLambda)$ is also acyclic. The category $\operatorname{Gproj}\varLambda$, with the short exact sequences inherited from $\mod\varLambda$, is Frobenius, and its projective-injective objects are precisely the projective $\varLambda$-modules. Its stable category $\underline{\operatorname{Gproj}}\varLambda$ is therefore triangulated, with suspension $[1]=\Omega_{\varLambda}^{-1}$; see, for example, \cite[Proposition~2.1.11 and Remark~2.1.12]{Che17}. Here $\Omega_{\varLambda}X$ denotes the kernel of a projective cover of $X$.

		The \textbf{singularity category} of $\varLambda$ is $D_{\mathrm{sg}}(\varLambda)=D^b(\mod\varLambda)/K^b(\operatorname{proj}\varLambda).$
		Here $D^b(\mod\varLambda)$ is the bounded derived category and $K^b(\operatorname{proj}\varLambda)$ is identified with its full subcategory of perfect complexes. Denote the quotient functor by $q_{\varLambda}$ and regard modules as complexes concentrated in degree zero.
		For a collection of objects $\mathcal X$ in a triangulated category $\mathcal D$, write $\operatorname{thick}_{\mathcal D}(\mathcal X)$ for the smallest full triangulated subcategory containing $\mathcal X$ and closed under direct summands.

        Unless explicitly stated otherwise, equivalences between singularity categories or between stable categories of Gorenstein-projective modules are understood to be $k$-linear triangle equivalences. Stable equivalences of algebras refer only to $k$-linear equivalences of their stable module categories, even when the algebras are selfinjective.
        
        	An algebra is \textbf{Gorenstein} (or \textbf{Iwanaga--Gorenstein}) if its left and right regular modules have finite injective dimension.  A $k$-linear category is called \textbf{Hom-finite} if all its morphism spaces are finite dimensional over $k$.

		\begin{proposition}\label{proposition of stable and singularity comparison}
			Let $\varLambda$ be an algebra. Inclusion and passage to stalk complexes give functors
			\[
				\iota_{\varLambda}:\underline{\operatorname{Gproj}}\varLambda\longrightarrow\smod\varLambda,
				\qquad
				\Phi_{\varLambda}:\underline{\operatorname{Gproj}}\varLambda\longrightarrow D_{\mathrm{sg}}(\varLambda).
			\]
			Then the following hold.
			\begin{enumerate}
				\item \textnormal{(cf. \cite[Proposition~2.3.1]{Che17})}
				The functor $\iota_{\varLambda}$ is fully faithful. It is an equivalence of $k$-categories if and only if $\varLambda$ is selfinjective.

				\item \textnormal{(cf. \cite[Theorem~4.4.1]{Buc87},
					\cite[Theorems~3.1 and~3.6]{BJO15})}
				The functor $\Phi_{\varLambda}$ is fully faithful. It is a triangle equivalence if and only if $\varLambda$ is Gorenstein.

				\item \textnormal{(cf. \cite[Section~1 and Lemma~3.5]{Che18}, \cite[Corollary~2.1.3(2)]{Che17})}
				If $\varLambda$ is selfinjective, all three categories agree up to their canonical triangle equivalences, with suspension $\Omega_{\varLambda}^{-1}$ on $\smod\varLambda$.
				If $\varLambda$ has finite global dimension, then $\underline{\operatorname{Gproj}}\varLambda=D_{\mathrm{sg}}(\varLambda)=0$, whereas $\smod\varLambda=0$ holds precisely when $\varLambda$ is semisimple.
			\end{enumerate}
		\end{proposition}

		\begin{remark}\label{remark of stabilization and singularity categories}
			We note that a stable equivalence $\smod\varLambda\simeq\smod\varGamma$ compatible with the syzygy functors and the standard left triangulated structures induces a triangle equivalence $D_{\mathrm{sg}}(\varLambda)\simeq D_{\mathrm{sg}}(\varGamma)$; see \cite[Example~2.3]{KV87} and \cite[Proposition~3.2]{Che18}. An equivalence of the underlying $k$-categories alone need not induce such an equivalence, as Example~\ref{example of node gluing different periods} shows.
		\end{remark}

		Two algebras are called \textbf{singularly equivalent} if their singularity categories are triangle equivalent.

	\subsection{Relation quivers of quadratic monomial algebras}\label{subsection of relation quiver criterion}

		Let $\varLambda=kQ/I$ be a quadratic monomial algebra.
		Following \cite[Definition~5.2]{CSZ18}, the \textbf{relation quiver} $\mathcal R_{\varLambda}$ has vertex set $Q_1$ and one arrow $a\rightarrow b$ for each pair of composable arrows $a,b$ satisfying $ab\in I$. A connected component of $\mathcal R_{\varLambda}$ is called \textbf{perfect} if it consists of a single directed cycle with no additional arrows; a single vertex with one loop is included. Connectedness refers to the underlying undirected graph, and a component is \textbf{acyclic} if it contains no directed cycle. 

		An arrow of $Q$ is called \textbf{perfect} if the corresponding vertex of $\mathcal R_{\varLambda}$ belongs to a perfect component. Let $\mathcal P(\varLambda)$ be the set of perfect components in $\mathcal R_{\varLambda}$. For each $C\in\mathcal P(\varLambda)$, let $\ell(C)$ denote the number of vertices of $C$.
		For $\ell\geq1$, let $C_\ell$ be the oriented cycle with $\ell$ vertices, let $J_\ell$ be the arrow ideal of $kC_\ell$, and put
		\[
			N_\ell=kC_\ell/J_\ell^2.
		\]
		The algebra $N_\ell$ is selfinjective, and its indecomposable nonprojective modules are its simple modules $T_i$, indexed by $i\in\mathbb Z/\ell\mathbb Z$ so that the arrows are $v_i\to v_{i+1}$. Then
		\[
			\sHom_{N_\ell}(T_i,T_j)=
			\begin{cases}k,&i=j,\\0,&i\ne j,\end{cases}
			\qquad \Omega_{N_\ell}T_i\cong T_{i+1},\qquad T_i[1]\cong T_{i-1}.
		\]
		Thus $\smod N_\ell$ is semisimple abelian, and its suspension has one orbit of length $\ell$ on the simple objects. Here a triangulated category is called \textbf{semisimple} if its underlying category is abelian and every short exact sequence in that abelian category splits.

		\begin{theorem}[{\cite[Proposition~5.1, Lemma~5.3, and Theorem~5.7]{CSZ18}}]
			\label{theorem of Gorenstein projectives and periods}
			Let $\varLambda=kQ/I$ be a quadratic monomial algebra.
			Then the following hold.
			\begin{enumerate}
				\item
				The indecomposable nonprojective Gorenstein-projective $\varLambda$-modules are precisely the modules $a\varLambda$ for perfect arrows $a$, up to isomorphism.
				If $a:v_i\to v_j$ and $a\to b$ is an arrow of a perfect component, the projective cover
				\[
					e_j\varLambda\longrightarrow a\varLambda,
					\qquad p\longmapsto ap,
				\]
				has kernel $b\varLambda$. In particular, $\Omega_{\varLambda}(a\varLambda)\cong b\varLambda$.

				\item
				There is a triangle equivalence
				\[
					\underline{\operatorname{Gproj}}\varLambda
					\simeq
					\prod_{C\in\mathcal P(\varLambda)}\smod N_{\ell(C)},
				\]
				where an empty product denotes the zero category.

				\item
				Put $\mathcal L(\varLambda)
				=\{\!\{\ell(C)\mid C\in\mathcal P(\varLambda)\}\!\}$. For any quadratic monomial algebra $\varGamma$, there is a triangle equivalence
				$\underline{\operatorname{Gproj}}\varLambda
				\simeq
				\underline{\operatorname{Gproj}}\varGamma$ if and only if $\mathcal L(\varLambda)=\mathcal L(\varGamma)$.
			\end{enumerate}
		\end{theorem}

		The relation quiver also determines the singularity category, by the following results of Chen.

		\begin{proposition}[{\cite[Theorem~4.4]{Che18},
				\cite[Theorem~B]{Che11}}]
			\label{proposition of singularity category via relation quiver}
			Let $\varLambda=kQ/I$ be a quadratic monomial algebra, and let $J$ be the arrow ideal of $k\mathcal R_{\varLambda}$. There is a triangle equivalence
			\[
				D_{\mathrm{sg}}(\varLambda)
				\simeq
				D_{\mathrm{sg}}(k\mathcal R_{\varLambda}/J^2).
			\]
			In particular, $D_{\mathrm{sg}}(\varLambda)$ is semisimple.
		\end{proposition}

		We next recall the criteria for Gorensteinness and finite global dimension in terms of the relation quiver.

		\begin{proposition}[{\cite[Proposition~5.5]{CSZ18}}]
			\label{proposition of relation quiver Gorenstein criterion}
			Let $\varLambda=kQ/I$ be a quadratic monomial algebra. Then the following hold.
			\begin{enumerate}
				\item
				$\varLambda$ is Gorenstein if and only if every connected component of $\mathcal R_{\varLambda}$ is perfect or acyclic.

				\item
				$\varLambda$ has finite global dimension if and only if every connected component of $\mathcal R_{\varLambda}$ is acyclic.
			\end{enumerate}
		\end{proposition}

		The Gorenstein criterion has the following equivalent formulation in terms of the arrows and relations of $(Q,I)$.

		\begin{proposition}\label{proposition of forbidden cycle criterion}
			Let $\varLambda=kQ/I$ be a quadratic monomial algebra. Then $\varLambda$ is not Gorenstein if and only if there is a cycle
			\[
				v_0\xrightarrow{a_0}v_1\xrightarrow{a_1}\cdots
				\xrightarrow{a_{\ell-2}}v_{\ell-1}
				\xrightarrow{a_{\ell-1}}v_0
			\]
			with $\ell\geq1$, pairwise distinct arrows
			$a_0,\ldots,a_{\ell-1}$, and $a_i a_{i+1}\in I$ for every $i$,
			such that at least one of the following holds for some $i$:
			\begin{enumerate}
				\item
				There is an arrow $b$ ending at $v_i$ with
				$b\ne a_{i-1}$ and $ba_i\in I$.

				\item
				There is an arrow $c$ starting at $v_i$ with
				$c\ne a_i$ and $a_{i-1}c\in I$.
			\end{enumerate}
			Indices are read modulo $\ell$, and the vertices $v_i$
			need not be distinct.
			For $\ell=1$, the displayed cycle is a loop.
		\end{proposition}

	\subsection{Node deletion}\label{subsection of triangulated node deletion}

		We now examine how node deletion, recalled in Subsection~\ref{subsection of node deletion}, affects stable categories of Gorenstein-projective modules and singularity categories.

		\subsubsection{Node deletion on Gorenstein-projective modules}

			Let $\varLambda=kQ/I$ be quadratic monomial, let $V$ be a set of its nodes, and let $\varLambda'$ be obtained by deleting the nodes in $V$, as in Subsection~\ref{subsection of node deletion}.

			Let $\mathcal P_V(\varLambda)$ consist of the perfect components of $\mathcal R_{\varLambda}$ containing an arrow of $Q$ which starts or ends at a vertex in $V$, and put
			\[
				\mathcal T_{\varLambda,V}	=
				\operatorname{thick}_{\underline{\operatorname{Gproj}}\varLambda}
				\{S_i\mid v_i\in V\text{ and }
				S_i\in\operatorname{Gproj}\varLambda\}.
			\]
			The thick subcategory generated by the empty set is zero.

			\begin{proposition}
				\label{proposition of node deletion Gorenstein projectives}
				Let $\varLambda=kQ/I$ be an algebra with $I$ generated by paths of length two. Let $V$ be a set of its nodes, and let $\varLambda'$ be obtained by deleting the nodes in $V$. 	Then there is a triangle equivalence
				\[
					\underline{\operatorname{Gproj}}\varLambda
					\simeq
					\underline{\operatorname{Gproj}}\varLambda'
					\times
					\prod_{C\in\mathcal P_V(\varLambda)}
					\smod N_{\ell(C)},
				\]
				under which $\mathcal T_{\varLambda,V}$ corresponds to the second factor.
				Consequently, $\underline{\operatorname{Gproj}}\varLambda/
				\mathcal T_{\varLambda,V}
				\simeq \underline{\operatorname{Gproj}}\varLambda'$.
			\end{proposition}

			\begin{proof}
				The relation quiver $\mathcal R_{\varLambda'}$ is obtained from $\mathcal R_{\varLambda}$ by deleting precisely the arrows $a\to b$ for which $t(a)=s(b)\in V$. In $\mathcal R_{\varLambda'}$, an arrow of $Q$ ending at a vertex in $V$ has no successor, and one starting at a vertex in $V$ has no predecessor. Thus, for every deleted arrow $a\to b$, neither $a$ nor $b$ lies on a directed cycle in $\mathcal R_{\varLambda'}$. Such a cycle is therefore a perfect component precisely when it was perfect before deletion, and
				\[
					\mathcal P(\varLambda')	=
					\mathcal P(\varLambda)\setminus
					\mathcal P_V(\varLambda).
				\]
				The product decomposition follows from Theorem~\ref{theorem of Gorenstein projectives and periods}(2).

				For $C\in\mathcal P_V(\varLambda)$, choose consecutive arrows $a,b$ on $C$ with $t(a)=s(b)=v_i\in V$. Every composition through $v_i$ is zero, so
				$a\varLambda=ka\cong S_i$. By Theorem~\ref{theorem of Gorenstein projectives and periods}(1), successive syzygies of $a\varLambda$ run through the modules
				corresponding to all vertices of $C$. Thus the image of $S_i$ generates the factor $\smod N_{\ell(C)}$ as a thick subcategory. Conversely, every Gorenstein-projective node simple $S_i$ with $v_i\in V$ is nonprojective, so the same theorem gives $S_i\cong a\varLambda$ for a perfect arrow $a$. The top of $a\varLambda$ is the simple module at $t(a)$, hence $t(a)=v_i$. The perfect component containing $a$ therefore belongs to $\mathcal P_V(\varLambda)$. It follows that $\mathcal T_{\varLambda,V}$ corresponds exactly to $\prod_{C\in\mathcal P_V(\varLambda)}\smod N_{\ell(C)}$. Projection onto the remaining factors gives the asserted Verdier quotient.
			\end{proof}

			\begin{remark}
				\label{remark of node functor and Gorenstein projectives}
				Note that the node-deletion functor $F$ induces an equivalence of
				stable module categories, but need not preserve Gorenstein-projective
				modules.
				The equivalences in
				Proposition~\ref{proposition of node deletion Gorenstein projectives}
				are constructed from the perfect-component decomposition.

				For example, take $\varLambda=k[x]/(x^2)$.
				Its unique simple module $S$ is a node and is Gorenstein-projective,
				since $\varLambda$ is selfinjective.
				Deleting the node gives $\varLambda'=k(v_+\rightarrow v_-)$.
				Since $Sx=0$, the definition of $F$ gives
				\[
					(FS)_{v_+}=k,\qquad (FS)_{v_-}=0,
				\]
				so $FS=S_+$ is the simple module at the source.
				This module is nonprojective, whereas
				$\operatorname{Gproj}\varLambda'=\operatorname{proj}\varLambda'$
				because $\varLambda'$ is hereditary.
				Thus $FS$ is not Gorenstein-projective.
			\end{remark}

		\subsubsection{Node deletion on singularity categories}

			We now describe the effect of node deletion on singularity categories.

			\begin{proposition}\label{proposition of node singularity quotients}
				Let $\varLambda=kQ/I$ be a bound quiver algebra, let $V$ be a set of nodes, and let $\varLambda'$ be obtained by deleting the nodes in $V$. Then there is a triangle equivalence
				\begin{equation}\label{eq node singularity quotients}
					\frac{D_{\mathrm{sg}}(\varLambda)}
					{\operatorname{thick}(q_{\varLambda}S_i\mid v_i\in V)}
					\simeq
					\frac{D_{\mathrm{sg}}(\varLambda')}
					{\operatorname{thick}(q_{\varLambda'}S_{i^+}\mid v_i\in V)}.
				\end{equation}
				If $\varLambda$ is quadratic monomial and Gorenstein, then $\varLambda'$ is Gorenstein and
				\[
					\frac{D_{\mathrm{sg}}(\varLambda)}
					{\operatorname{thick}(q_{\varLambda}S_i\mid v_i\in V)}
					\simeq D_{\mathrm{sg}}(\varLambda').
				\]
			\end{proposition}

			\begin{proof}
				For an algebra $\varGamma$ and an idempotent $h\in\varGamma$, restriction $X\mapsto Xh$ gives a triangle equivalence
				\[
					\frac{D_{\mathrm{sg}}(\varGamma)}
					{\operatorname{thick}(q_{\varGamma}S\mid
						S\text{ simple},\ Sh=0)}
					\simeq
					\frac{D^b(\mod h\varGamma h)}
					{\operatorname{thick}(\varGamma h)}
				\]
				by \cite[Proposition~3.3]{KY16}, with $\varGamma h$ viewed as a right $h\varGamma h$-module.
				Let $e\in\varLambda$ and $f\in\varLambda'$ be the sums of the vertex idempotents outside $V$, excluding both copies of each split node on the $\varLambda'$-side. If $e=0$, both quotients in \eqref{eq node singularity quotients} are zero. For $e\ne0$, since a path with a node as an internal vertex is zero, the path correspondence gives
				\[
					R:=e\varLambda e\cong f\varLambda'f,
					\qquad
					\varLambda e\cong\varLambda'f\quad\text{in }\mod R.
				\]
				In the module isomorphism, a path starting at $v_i\in V$ corresponds to the path starting at $v_{i^+}$; no path starting at $v_{i^-}$ ends outside the split nodes.
				Thus the two idempotent quotients are both $D^b(\mod R)/\operatorname{thick}(\varLambda e)$. The simples annihilated by $e$ are precisely the $S_i$ with
				$v_i\in V$, and those annihilated by $f$ are the $S_{i^+}$ and $S_{i^-}$ with $v_i\in V$. Since each $S_{i^-}$ is projective, its image in $D_{\mathrm{sg}}(\varLambda')$ is zero. This proves \eqref{eq node singularity quotients}.

				Suppose now that $\varLambda$ is quadratic monomial and Gorenstein. Node deletion cuts the components in $\mathcal P_V(\varLambda)$ into paths, leaves the other perfect components unchanged, and preserves acyclicity of the acyclic components.
				Hence $\varLambda'$ is Gorenstein by Proposition~\ref{proposition of relation quiver Gorenstein criterion}. Each $S_{i^+}$ is injective and therefore has finite projective dimension over $\varLambda'$. Thus the right-hand denominator in \eqref{eq node singularity quotients} is zero, giving the second equivalence.
			\end{proof}

			%In general, both denominators in \eqref{eq node singularity quotients} must be retained. For a quadratic monomial algebra, the left-hand denominator uses all simples $S_i$ with $v_i\in V$, whereas $\mathcal T_{\varLambda,V}$ is generated only by those which are Gorenstein-projective.

			The following example shows that an equivalence of stable module
			categories as $k$-categories need not give a triangle equivalence
			of singularity categories.

			\begin{example}\label{example of node gluing different periods}
				Let $\varLambda_0$ be the product of the path algebras of $v_{s_1}\xrightarrow{a}v_{t_1}$ and $v_{s_2}\xrightarrow{b}v_{t_2}$. Put
				\[
					\varLambda_1=k[x]/(x^2)\times k[y]/(y^2),
					\qquad
					\varLambda_2=N_2.
				\]
				Deleting all nodes of either algebra gives $\varLambda_0$. Hence $\varLambda_1$ and $\varLambda_2$ are stably equivalent by Proposition~\ref{proposition of node deletion}. Both algebras are selfinjective, and their relation quivers consist of two loops and one directed cycle of length two, respectively. Thus
				\[
					D_{\mathrm{sg}}(\varLambda_1)
					\simeq\smod N_1\times\smod N_1,
					\qquad
					D_{\mathrm{sg}}(\varLambda_2)
					\simeq\smod N_2.
				\]
				The suspension fixes both simple objects in the first category
				and exchanges them in the second, so the two singularity categories
				are not triangle equivalent.
			\end{example}

		\subsubsection{Consequences of stable equivalence}

			The preceding results give the following consequences for stably equivalent quadratic monomial algebras. For a quadratic monomial algebra $A=kQ/I$, let $V_A$ be the set of all node vertices of $Q$ and put $\mathcal T_A:=\mathcal T_{A,V_A}$. Write $\widehat A$ for the algebra obtained by deleting all nodes of $A$ and removing all semisimple direct summands. We take $\underline{\operatorname{Gproj}}0=0$.

			\begin{corollary}
				\label{corollary of stable Gorenstein quotient invariant}
				Let $\varLambda$ and $\varGamma$ be stably equivalent algebras, and suppose that $\varLambda$ is quadratic monomial. Then $\varGamma$ is quadratic monomial, and there is a triangle equivalence
				\[
					\underline{\operatorname{Gproj}}\varLambda/\mathcal T_{\varLambda}
					\simeq
					\underline{\operatorname{Gproj}}\varGamma/\mathcal T_{\varGamma}.
				\]
			\end{corollary}

			\begin{proof}
				By Theorem~\ref{theorem of quadratic monomial stable closure}, $\varGamma$ is quadratic monomial, and that theorem gives $\widehat{\varLambda}\cong\widehat{\varGamma}$.
				Since semisimple direct summands have zero stable Gorenstein-projective category, Proposition~\ref{proposition of node deletion Gorenstein projectives} gives
				$\underline{\operatorname{Gproj}}A/\mathcal T_A\simeq\underline{\operatorname{Gproj}}\widehat A$ where $A\in\{\varLambda,\varGamma\}$. The isomorphism of the reduced algebras therefore gives the required triangle equivalence.
			\end{proof}

			\begin{remark}
				\label{remark of stable Gorenstein periods}
				For a quadratic monomial algebra $A$, put $\mathcal L_{\mathrm{nd}}(A):=
				\{\!\{\ell(C)\mid C\in\mathcal P_{V_A}(A)\}\!\}$. Each perfect component meeting a node is counted once, even if it meets several nodes. For stably equivalent quadratic monomial algebras $\varLambda$ and $\varGamma$, we have
				\[
					\underline{\operatorname{Gproj}}\varLambda
					\simeq
					\underline{\operatorname{Gproj}}\varGamma
					\quad\Longleftrightarrow\quad
					\mathcal L_{\mathrm{nd}}(\varLambda)
					=
					\mathcal L_{\mathrm{nd}}(\varGamma).
				\]
				If both algebras are Gorenstein, these conditions are also equivalent to $D_{\mathrm{sg}}(\varLambda)\simeq D_{\mathrm{sg}}(\varGamma)$.
			\end{remark}

\section{Marked configurations associated with almost gentle algebras}\label{section of marked configurations}

	We now apply the results for quadratic monomial algebras to almost gentle and string almost gentle algebras. After recalling their definitions and proving their stable-closure properties, we describe equivalences in terms of marked configurations.

	\subsection{Almost gentle and string almost gentle algebras}\label{subsection of almost gentle algebras}

		Our proof also applies to almost gentle and string almost gentle algebras; see Theorems~{\ref{theorem of AG stable closure}} and~{\ref{theorem of SAG stable closure}}.
		We recall their definitions and the properties needed below.

		\begin{definition}[{cf.~\cite{GR05,GS19,FGR22}}]
			An algebra $\varLambda$ is called a \textbf{gentle algebra} if it admits a presentation $\varLambda=kQ/I$, where $I$ is an admissible ideal, satisfying the following conditions:
			\begin{enumerate}
				\item At most two arrows start and at most two arrows end at each vertex of $Q$.

				\item For each arrow $a$, there is at most one composable successor $b$ with $ab\notin I$ and at most one composable	predecessor $c$ with $ca\notin I$.

				\item For each arrow $a$, there is at most one composable successor $b$ with $ab\in I$ and at most one composable predecessor $c$ with $ca\in I$.

				\item The ideal $I$ is generated by paths of length two.
			\end{enumerate}
			An algebra is called a \textbf{string almost gentle algebra} if it admits a bound quiver presentation satisfying conditions (1), (2), and (4), and an \textbf{almost gentle algebra} if it admits a bound quiver presentation satisfying conditions (2) and (4).
		\end{definition}

		\begin{lemma}\label{lemma of AG node type}
			Let $\varLambda=kQ/I$ be an almost gentle algebra, and let $\varLambda'=kQ'/I'$ be an algebra stably equivalent to $\varLambda$ by a stable equivalence of node type.
			Then $\varLambda'$ is an almost gentle algebra.
		\end{lemma}

		\begin{proof}
			It suffices to consider a single node deletion or its reverse.
			Suppose that $\varLambda'$ is obtained from $\varLambda$ by deleting a node $v_i$.
			The ideal $I'$ is generated by the paths of length two in $I$ whose arrows are composable in $Q'$.
			For arrows $a,b$ composable in $Q'$, we have $ab\in I'$ if and only if $ab\in I$.
			Thus, for each arrow $a$ in $Q'$, there is at most one arrow $b$ with $ab\notin I'$ and at most one arrow $c$ with $ca\notin I'$.
			Hence $\varLambda'$ is an almost gentle algebra.

			Conversely, suppose that $\varLambda'$ is obtained by gluing a source vertex $v_{i^+}$ and a sink vertex $v_{i^-}$ to form $v_i$.
			The incoming arrows at $v_i$ in $Q'$ are precisely those ending at $v_{i^-}$ in $Q$, and the outgoing arrows at $v_i$ in $Q'$ are precisely those starting at $v_{i^+}$ in $Q$.
			In particular, arrows $v_{i^+}\to v_{i^-}$ become loops.
			The ideal $I'$ is generated by the paths of length two in $I$ and all paths of length two through $v_i$.
			If $ab\notin I'$ for arrows $a,b$ in $Q'$, then $a,b$ are composable in $Q$ and $ab\notin I$.
			Consequently, for each arrow $a$ in $Q'$, there is at most one arrow $b$ with $ab\notin I'$ and at most one arrow $c$ with $ca\notin I'$.
			Hence $\varLambda'$ is an almost gentle algebra.
		\end{proof}

		\begin{remark}\label{remark of SAG node type}
			The same conclusion holds for string almost gentle algebras.
			Indeed, node deletion separates the incoming and outgoing arrows at $v_i$ between $v_{i^-}$ and $v_{i^+}$, while gluing combines the incoming arrows at the sink vertex with the outgoing arrows at the source vertex.
			There are still at most two arrows entering and at most two arrows leaving each vertex, so the assertion follows from Lemma~\ref{lemma of AG node type}.
		\end{remark}

	\subsection{Stable equivalence of almost gentle algebras}

		We now give the proofs using the quadratic monomial case and the preservation of these classes under node deletion and its reverse.

		\begin{theorem}\label{theorem of AG stable closure}
			Let $\varLambda$ be an almost gentle algebra and let $\varLambda'$ be a basic algebra stably equivalent to $\varLambda$.
			Then $\varLambda'$ is an almost gentle algebra.
		\end{theorem}

		\begin{proof}
			Let $\varLambda=kQ/I$ be an almost gentle algebra, and let $\varLambda'=kQ'/I'$ be a basic algebra stably equivalent to $\varLambda$.
			By Theorem~\ref{theorem of quadratic monomial stable closure}, deleting all nodes and removing semisimple summands gives isomorphic algebras $\widehat{\varLambda}$ and $\widehat{\varLambda}'$.
			By Lemma~\ref{lemma of AG node type}, $\widehat{\varLambda}$ is almost gentle if it is nonzero, and the same holds for $\widehat{\varLambda}'$.
		
			Restore the semisimple factors and reverse the node deletions on the $\varLambda'$ side.
			These operations preserve the almost gentle property, so $\varLambda'$ is almost gentle.
		\end{proof}

		\begin{theorem}\label{theorem of SAG stable closure}
			Let $\varLambda$ be a string almost gentle algebra and let $\varLambda'$ be a basic algebra stably equivalent to $\varLambda$.
			Then $\varLambda'$ is a string almost gentle algebra.
		\end{theorem}

		\begin{proof}
			Let $\varLambda=kQ/I$ be a string almost gentle algebra, and let $\varLambda'=kQ'/I'$ be a basic algebra stably equivalent to $\varLambda$.
			Repeat the proof of Theorem~\ref{theorem of AG stable closure}, using Remark~\ref{remark of SAG node type} in place of Lemma~\ref{lemma of AG node type}.
		\end{proof}

		\begin{remark}\label{remark of Morita type closure}
			By Theorem~\ref{corollary of AG Morita rigidity}(1), the classes of quadratic monomial, gentle, string almost gentle, and almost gentle algebras are closed under stable equivalence of Morita type.
		\end{remark}

		The Gorenstein-projective quotient equivalence in Corollary~\ref{corollary of stable Gorenstein quotient invariant} also applies to these algebras.

	\subsection{Marked configurations}

		Whenever we use marked configurations below, the algebras are assumed to have no semisimple direct summands.

		We use the Brauer configurations of Green--Schroll \cite[Definition~1.5]{GS17BC} with all vertex multiplicities equal to one, and omit the multiplicity function.

		\begin{definition}\label{definition of marked configuration}
			A \textbf{configuration} is a triple $\Gamma=(\Gamma_0,\Gamma_1,\mathfrak o)$ consisting of a finite vertex set $\Gamma_0$, a finite indexed collection $\Gamma_1$ of finite multisets of vertices, called \textbf{polygons}, and a cyclic ordering $\mathfrak o$ of the occurrences at each vertex. Write $\operatorname{val}(\vv)$ for the total number of occurrences of a vertex $\vv$. 
			We assume that each vertex occurs in a polygon, each polygon has at least two occurrences, and each polygon contains a vertex $\vv$ with $\operatorname{val}(\vv)>1$.
			A vertex $\vv$ is \textbf{truncated} if $\operatorname{val}(\vv)=1$. A \textbf{marking} chooses a corner between consecutive occurrences at each nontruncated vertex; cutting at that corner gives a linear ordering. A configuration equipped with a marking is a \textbf{marked configuration}.
		\end{definition}

		Let $A=kQ/I$ be almost gentle. We recall its marked configuration $C_A$ from \cite{GS19}. For each $v\in Q_0$, take a polygon $P_v$. Each nontrivial maximal nonzero path
		\[
			p:\quad v_0\xrightarrow{a_1}v_1\longrightarrow\cdots
			\xrightarrow{a_r}v_r\qquad(r\geq1)
		\]
		gives a configuration vertex $\vv_p$ with occurrences $h_{p,0},\ldots,h_{p,r}$, where $h_{p,i}$ belongs to $P_{v_i}$. Close this list cyclically and mark the corner $h_{p,r}\to h_{p,0}$. If a polygon has only one occurrence, add one truncated vertex to that polygon. Occurrences remain distinct when a path visits a quiver vertex more than once. Since $Q$ has no isolated vertices, every polygon receives an occurrence.

		Note that the underlying configuration is that of the trivial extension of $A$, and the marking records the admissible cut recovering $A$ \cite{GS19}. When $A$ is gentle, every polygon has two occurrences; viewing polygons as edges gives the marked ribbon graph $C_A$  \cite{Sch15,OPS}.

		For such a configuration $C$, let $H:=H(C)$ be the set of occurrences at nontruncated vertices and let $\Gamma_1:=\Gamma_1(C)$ be its polygon set. Denote cyclic successor and polygon incidence by
		\[
			\sigma:H\longrightarrow H,\qquad
			\pi:H\longrightarrow\Gamma_1.
		\]
		Let $D\subseteq H$ consist of the occurrences whose outgoing corners are marked, and put $U=H\setminus D$. We identify $h\in U$ with the uncut corner $h\to\sigma h$. For another configuration $C'$, use primed notation. An \textbf{isomorphism} $C\cong C'$ is given by bijections $f:H\to H'$ and $g:\Gamma_1\to\Gamma_1'$ satisfying
		\[
			f\sigma=\sigma'f,\qquad \pi'f=g\pi,\qquad f(D)=D'.
		\]
		For configurations constructed in this way, the truncated vertices are uniquely determined and require no additional data.

		For $h,g\in U$, write $h\rightsquigarrow g$ if $\pi(\sigma h)=\pi(g)$ and $g\ne\sigma h$.

		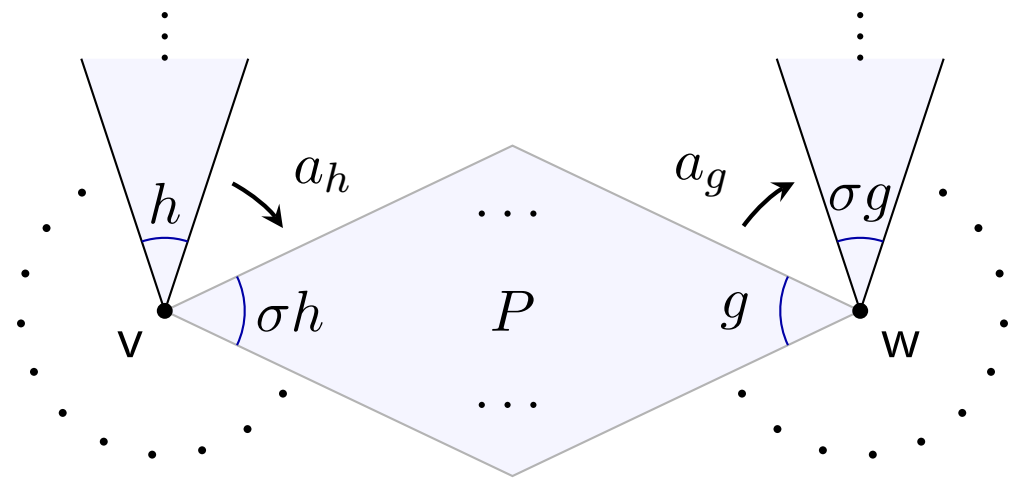
\begin{figure}[htbp]
			\centering
			\begin{tikzpicture}[scale=1.15,>=stealth,
				occurrence/.style={draw=blue!65!black},
				quiver arrow/.style={->,thick}]
				% Polygons; the upper polygons continue beyond the drawing.
				\filldraw[fill=blue!4,draw=gray!60]
					(0,0)--(2,.95)--(4,0)--(2,-.95)--cycle;
				\fill[blue!4] (0,0)--(-.48,1.45)--(.48,1.45)--cycle;
				\draw (0,0)--(-.48,1.45) (0,0)--(.48,1.45);
				\fill[blue!4] (4,0)--(3.52,1.45)--(4.48,1.45)--cycle;
				\draw (4,0)--(3.52,1.45) (4,0)--(4.48,1.45);
				\node at (0,1.67) {$\vdots$};
				\node at (4,1.67) {$\vdots$};
				% Occurrences are represented by angles.
				\draw[occurrence] (108.3:.42)
					arc[start angle=108.3,end angle=71.7,radius=.42];
				\node at (0,.62) {$h$};
				\draw[occurrence] (25.4:.46)
					arc[start angle=25.4,end angle=-25.4,radius=.46];
				\node at (.72,0) {$\sigma h$};
				\begin{scope}[shift={(4,0)}]
					\draw[occurrence] (154.6:.46)
						arc[start angle=154.6,end angle=205.4,radius=.46];
					\draw[occurrence] (108.3:.42)
						arc[start angle=108.3,end angle=71.7,radius=.42];
				\end{scope}
				\node at (3.28,0) {$g$};
				\node at (4,.62) {$\sigma g$};
				\fill (0,0) circle (1.3pt) (4,0) circle (1.3pt);
				\node[below left] at (0,0) {$\vv$};
				\node[below right] at (4,0) {$\ww$};
				\node at (2,0) {$P$};
				\node at (2,.55) {$\cdots$};
				\node at (2,-.55) {$\cdots$};
				% Each arrow lies in the uncut corner at its own configuration vertex.
				\draw[quiver arrow] (62:.83)
					arc[start angle=62,end angle=35,radius=.83];
				\node at (.91,.78) {$a_h$};
				\begin{scope}[shift={(4,0)}]
					\draw[quiver arrow] (144:.83)
						arc[start angle=144,end angle=117,radius=.83];
				\end{scope}
				\node at (3.09,.78) {$a_g$};
				\foreach \ang in {125,145,...,325}
				\fill (\ang:.83) circle (.65pt);
				\begin{scope}[shift={(4,0)}]
					\foreach \ang in {-145,-125,...,55}
					\fill (\ang:.83) circle (.65pt);
				\end{scope}
			\end{tikzpicture}
			\caption{A transition $h\rightsquigarrow g$. The quiver arrows $a_h$ and $a_g$ lie at different configuration vertices, and $a_ha_g=0$.}
			\label{figure of occurrence transition}
		\end{figure}

		\begin{lemma}\label{lemma of configuration relation quiver}
			For $C=C_A$, the arrows of $Q$ are in bijection with $U$: the corner $h\to\sigma h$ corresponds to an arrow $a_h:\pi(h)\to\pi(\sigma h)$, where $P_v$ is identified with $v$. The relation quiver $\mathcal R_A$ is the directed graph with vertex set $U$ and arrows given by the transitions $h\rightsquigarrow g$.
		\end{lemma}

		\begin{proof}
			The assertions follow directly from the definitions.
		\end{proof}

	\subsection{Stable equivalence}

		We first note that the relation quiver alone does not determine stable equivalence. For example, the path algebras of $1\to2\to3$ and $1\to2\leftarrow3$ both have relation quivers consisting of two isolated vertices. However, these algebras have no nodes and are not isomorphic, so they are not stably equivalent by Theorem~\ref{theorem of AG rigidity}.
		Thus, to characterize stable equivalence, we also need to keep track of which occurrences belong to the same polygon in the marked configuration.

		We begin by describing algebra isomorphisms in terms of quiver presentations.

		\begin{lemma}\label{lemma of monomial presentation isomorphism}
			Let $A=kQ/I$ and $B=kQ'/I'$, where $I$ and $I'$ are admissible ideals, and let $\varphi:A\to B$ be an algebra isomorphism.
			\begin{enumerate}
				\item
				There is a bijection $\theta_0:Q_0\to Q'_0$ such that, for every $u,v\in Q_0$, $\varphi$ induces a linear isomorphism between the subspaces of $\rad A/\rad^2 A$ and $\rad B/\rad^2 B$ spanned by the residue classes of arrows from $u$ to $v$ and from $\theta_0(u)$ to $\theta_0(v)$, respectively.

				\item
				If both $I$ and $I'$ are monomial, then $\theta_0$ extends to a quiver isomorphism $\theta:Q\to Q'$ whose extension to path algebras sends $I$ onto $I'$.
			\end{enumerate}
		\end{lemma}

		\begin{proof}
			\underline{(1)}
			The isomorphism $\varphi$ induces an isomorphism $A/\rad A\cong B/\rad B$, which permutes the primitive idempotents and determines a bijection $\theta_0:Q_0\to Q'_0$. It also induces a linear isomorphism
			\[
				\rad A/\rad^2 A\cong\rad B/\rad^2 B.
			\]
			Since $I$ and $I'$ are admissible, the residue classes of the arrows form bases of these spaces. Compatibility with multiplication by the primitive idempotents gives the asserted restrictions.

			\underline{(2)}
			Since $I$ and $I'$ are homogeneous with respect to path length, $A$ and $B$ identify naturally with their associated graded algebras for the radical filtrations. Replacing $\varphi$ by its induced graded isomorphism, we may therefore assume that $\varphi$ is graded and $\varphi(e_u)=e'_{\theta_0(u)}$.

			For each pair of vertices, $\varphi$ sends the arrows between them to linear combinations of arrows between the corresponding vertices, with an invertible coefficient matrix. Fix $u,v\in Q_0$, and let $a_1,\ldots,a_m$ be the arrows from $u$ to $v$ and $b_1,\ldots,b_m$ the arrows from $\theta_0(u)$ to $\theta_0(v)$. Write
			\[
				\varphi(a_i)=\sum_{j=1}^m c_{ij}b_j.
			\]
			The matrix $(c_{ij})$ is invertible, so its determinant expansion contains a nonzero term. Thus there is a permutation $\tau\in S_m$ such that
			$$\prod_{i=1}^m c_{i,\tau(i)}\ne0.$$
			Define $\theta(a_i)=b_{\tau(i)}$. Since $\tau$ is a permutation, this gives a bijection between the two sets of arrows, and the coefficient of $\theta(a_i)$ in $\varphi(a_i)$ is nonzero for every $i$.

			Repeating this construction for each pair of vertices, together with $\theta(u)=\theta_0(u)$ on vertices, defines a quiver isomorphism $\theta:Q\to Q'$. We denote its extension to path algebras by the same symbol.
			Let $p=a_1\cdots a_r$ be a path in $I$. If $\theta(p)\notin I'$, its coefficient in $\varphi(a_1)\cdots\varphi(a_r)$ is the product of the selected nonzero coefficients, since no other choice of arrows in the expansion produces the same path. The residue classes of paths outside $I'$ form a basis of $B$, so this contradicts $\varphi(p)=0$. Since $I$ is generated by paths, we obtain $\theta(I)\subseteq I'$. Finally,
			$\dim_k kQ'/\theta(I)=\dim_k kQ/I=\dim_k kQ'/I'$,
			and hence $\theta(I)=I'$.
		\end{proof}

		\begin{corollary}\label{corollary of marked configuration isomorphism}
			Let $A$ and $B$ be almost gentle algebras. The following are equivalent.
			\begin{enumerate}
				\item $A\cong B$.
				\item $A$ and $B$ are stably equivalent of Morita type.
				\item $C_A\cong C_B$ as marked configurations.
			\end{enumerate}
		\end{corollary}

		\begin{proof}
			The equivalence $(1)\Leftrightarrow(2)$ follows from Theorem~\ref{corollary of AG Morita rigidity}(1). It remains to prove $(1)\Rightarrow(3)$ and $(3)\Rightarrow(1)$. Write $A=kQ/I$ and $B=kQ'/I'$.

			\noindent\underline{$(1)\Rightarrow(3)$}
			By Lemma~\ref{lemma of monomial presentation isomorphism}(2), an algebra isomorphism $A\cong B$ gives a quiver isomorphism $\theta:Q\to Q'$ sending $I$ onto $I'$. By the construction of $C_A$ and $C_B$, this induces an isomorphism of marked configurations.

			\noindent\underline{$(3)\Rightarrow(1)$}
			Let $(f,g):C_A\to C_B$ be an isomorphism of marked configurations, where $f$ is the occurrence bijection and $g$ is the polygon bijection. Identifying polygons with quiver vertices, the assignments
			\[
				P\longmapsto g(P),\qquad a_h\longmapsto a'_{f(h)}
			\]
			define a quiver isomorphism $Q\cong Q'$. By Lemma~\ref{lemma of configuration relation quiver}, this isomorphism preserves the quadratic relations. Since $I$ and $I'$ are generated by quadratic paths, it induces an algebra isomorphism $A\cong B$.
		\end{proof}

		To describe stable equivalences, we now consider how node deletion changes the polygon incidences.
		For a marked configuration $C$, cutting at the marking at a nontruncated vertex $\vv$ gives a list
		\[
			h_{\vv,0}<\cdots<h_{\vv,r_\vv},\qquad r_\vv\geq1.
		\]
		Let $H^+$, $H^-$, and $H^\circ$ be the sets of first, last, and internal occurrences in these lists. For $P\in\Gamma_1$, put $H_P=\pi^{-1}(P)$ and define
		\[
			\begin{aligned}
				\mathscr I(C)&=\{H_P\mid H_P\cap H^\circ\ne\varnothing\},\\
				\mathscr S(C)&=\{H_P\cap H^+\mid H_P\cap H^\circ=\varnothing,
				\ H_P\cap H^+\ne\varnothing\},\\
				\mathscr T(C)&=\{H_P\cap H^-\mid H_P\cap H^\circ=\varnothing,
				\ H_P\cap H^-\ne\varnothing\}.
			\end{aligned}
		\]

		\begin{definition}\label{definition of incidence comparison}
			An \textbf{incidence comparison} $C\to C'$ is a bijection $f:H(C)\to H(C')$ which sends each cut-linear list onto a cut-linear list in its given order and satisfies
			\[
				f\bigl(\mathscr I(C)\bigr)=\mathscr I(C'),\qquad
				f\bigl(\mathscr S(C)\bigr)=\mathscr S(C'),\qquad
				f\bigl(\mathscr T(C)\bigr)=\mathscr T(C').
			\]
		\end{definition}

		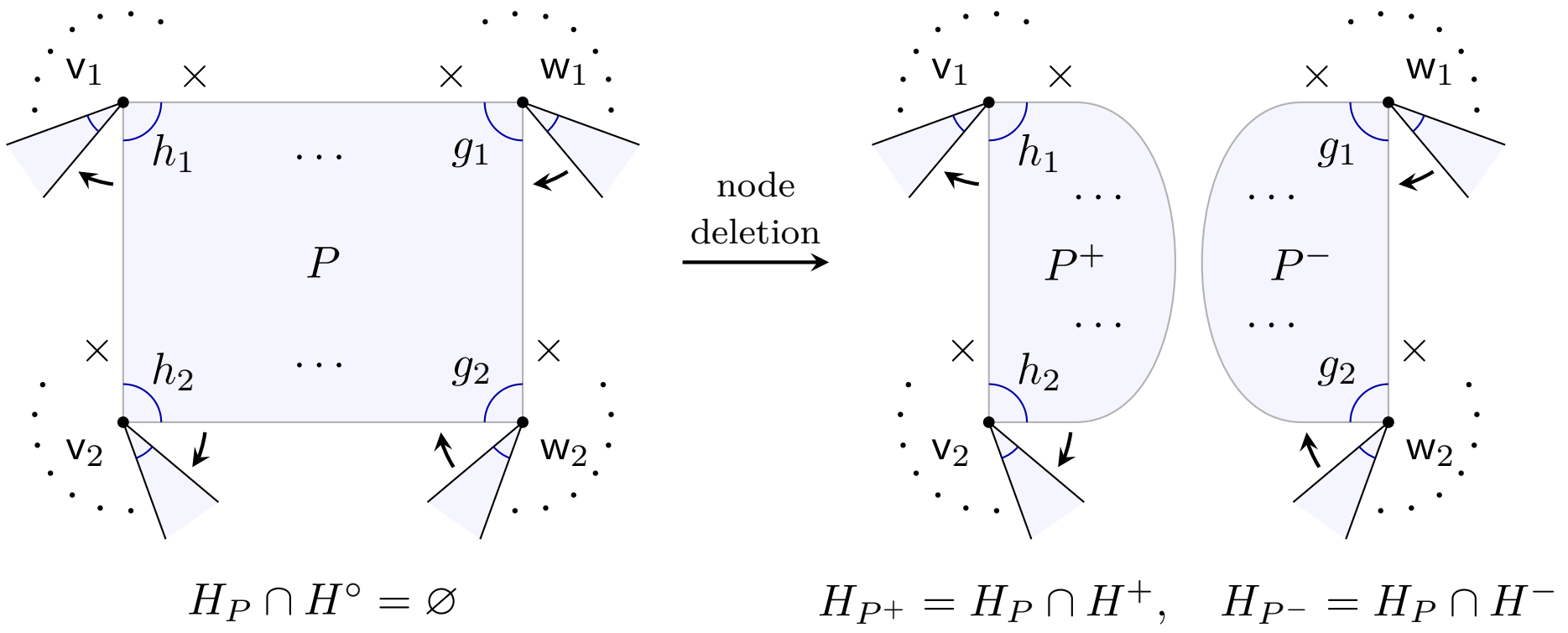
\begin{figure}[htbp]
			\centering
			\begin{tikzpicture}[scale=1.05,>=stealth,font=\small,
				polygon/.style={fill=blue!4,draw=gray!60},
				occurrence/.style={draw=blue!65!black},
				quiver arrow/.style={->,thick},
				pics/first/.style={code={
						% The occurrence in the central polygon is first in its cut-linear list.
						\fill[blue!4] (0,0)--(-115:.98)--(-85:.98)--cycle;
						\draw (0,0)--(-115:.98) (0,0)--(-85:.98);
						\draw[occurrence] (-115:.30)
							arc[start angle=-115,end angle=-85,radius=.30];
						\draw[occurrence] (45:.30)
							arc[start angle=45,end angle=-45,radius=.30];
						\draw[quiver arrow] (-52:.65)
							arc[start angle=-52,end angle=-78,radius=.65];
						\node at (65:.60) {$\times$};
						\foreach \ang in {110,130,...,230}
						\fill (\ang:.70) circle (.6pt);
						\fill (0,0) circle (1.3pt);
				}},
				pics/last/.style={code={
						% The occurrence in the central polygon is last in its cut-linear list.
						\fill[blue!4] (0,0)--(85:.98)--(115:.98)--cycle;
						\draw (0,0)--(85:.98) (0,0)--(115:.98);
						\draw[occurrence] (115:.30)
							arc[start angle=115,end angle=85,radius=.30];
						\draw[occurrence] (45:.30)
							arc[start angle=45,end angle=-45,radius=.30];
						\draw[quiver arrow] (78:.65)
							arc[start angle=78,end angle=52,radius=.65];
						\node at (-65:.60) {$\times$};
						\foreach \ang in {130,150,...,250}
						\fill (\ang:.70) circle (.6pt);
						\fill (0,0) circle (1.3pt);
				}}]
				% Before node deletion: one polygon contains first and last occurrences.
				\filldraw[polygon] (0,1.2)--(3,1.2)--(3,-1.2)--(0,-1.2)--cycle;
				\node at (1.5,0) {$P$};
				\node at (1.5,.78) {$\cdots$};
				\node at (1.5,-.78) {$\cdots$};
				% After node deletion: separate the two types of occurrences.
				\begin{scope}[shift={(6.5,0)}]
					\filldraw[polygon] (0,1.2)--(.65,1.2)
						.. controls (1.65,1.2) and (1.65,-1.2) .. (.65,-1.2)
						--(0,-1.2)--cycle;
					\filldraw[polygon] (3,1.2)--(2.35,1.2)
						.. controls (1.35,1.2) and (1.35,-1.2) .. (2.35,-1.2)
						--(3,-1.2)--cycle;
					\node at (.65,0) {$P^+$};
					\node at (2.35,0) {$P^-$};
					\node at (.85,.48) {$\cdots$};
					\node at (.85,-.48) {$\cdots$};
					\node at (2.15,.48) {$\cdots$};
					\node at (2.15,-.48) {$\cdots$};
				\end{scope}
				% The same configuration vertices, angles, markings and arrows on both sides.
				\foreach \offset in {0,6.5}{
					\begin{scope}[shift={(\offset,0)}]
						\pic[rotate=-45] at (0,1.2) {first};
						\pic[rotate=45] at (0,-1.2) {first};
						\pic[rotate=225] at (3,1.2) {last};
						\pic[rotate=135] at (3,-1.2) {last};
						\node at (.38,.82) {$h_1$};
						\node at (.38,-.82) {$h_2$};
						\node at (2.62,.82) {$g_1$};
						\node at (2.62,-.82) {$g_2$};
						\node[above left] at (0,1.2) {$\vv_1$};
						\node[below left] at (0,-1.2) {$\vv_2$};
						\node[above right] at (3,1.2) {$\ww_1$};
						\node[below right] at (3,-1.2) {$\ww_2$};
					\end{scope}
				}
				\draw[->,thick] (4.2,0)--(5.3,0)
					node[midway,above,align=center,font=\scriptsize] {node\\deletion};
				\node at (1.5,-2.55) {$H_P\cap H^\circ=\varnothing$};
				\node at (8,-2.55)
					{$H_{P^+}=H_P\cap H^+,\quad H_{P^-}=H_P\cap H^-$};
			\end{tikzpicture}
			\caption{Node deletion in a marked configuration.}
			\label{figure of configuration node deletion}
		\end{figure}

		\begin{proposition}\label{proposition of configuration stable criterion}
			Let $A$ and $B$ be almost gentle algebras. Then $A$ and $B$ are stably equivalent if and only if there is an incidence comparison $C_A\to C_B$.
		\end{proposition}

		\begin{proof}
			Let $\widehat A$ and $\widehat B$ be obtained by deleting all nodes of $A$ and $B$.
			In $C_A$, a polygon $P$ corresponds to the middle vertex of a path of length two that is nonzero in $A$ if and only if $H_P\cap H^\circ\ne\varnothing$.
			If $H_P\cap H^\circ=\varnothing$, that vertex is a node exactly when $H_P$ meets both $H^+$ and $H^-$. Node deletion splits $H_P$ into the two subsets $H_P\cap H^+$ and $H_P\cap H^-$, as illustrated in Figure~\ref{figure of configuration node deletion}.

			Thus $C_{\widehat A}$ has the same nontruncated occurrences and cut-linear lists as $C_A$, and the sets of nontruncated occurrences belonging to its polygons form the family $\mathscr I(C_A)\cup\mathscr S(C_A)\cup\mathscr T(C_A)$. The truncated vertices are determined by these sets.
			The same description applies to $C_{\widehat B}$. Consequently, an incidence comparison $C_A\to C_B$ is precisely an isomorphism $C_{\widehat A}\cong C_{\widehat B}$. By Corollary~\ref{corollary of marked configuration isomorphism}, such an isomorphism exists if and only if $\widehat A\cong\widehat B$.

			\noindent\underline{$(\Rightarrow)$}
			If $A$ and $B$ are stably equivalent, so are $\widehat A$ and $\widehat B$. Theorem~\ref{theorem of AG rigidity} gives $\widehat A\cong\widehat B$, and hence an incidence comparison $C_A\to C_B$.

			\noindent\underline{$(\Leftarrow)$}
			An incidence comparison $C_A\to C_B$ gives $\widehat A\cong\widehat B$. Composing this isomorphism with the node-deletion equivalences and their inverses gives a stable equivalence between $A$ and $B$.
		\end{proof}

	\subsection{{Singularity classes}}

		We now turn from stable equivalence to singular equivalence.
		We first recall the Gorenstein properties of gentle, string almost gentle, and almost gentle algebras, and then compare the singularity classes of gentle and almost gentle algebras.

		\begin{remark}\label{remark of three classes Gorenstein}
			Note that, for all three classes above, finite global dimension is equivalent to acyclicity of the relation quiver by Proposition~\ref{proposition of relation quiver Gorenstein criterion}(2). Now we discuss the Gorenstein property of these three classes.
			\begin{enumerate}
				\item
				Every gentle algebra is Gorenstein by Gei\ss{}--Reiten 	\cite[Theorem~3.4]{GR05}.

				\item
				By Proposition~\ref{proposition of relation quiver Gorenstein criterion}(1), 	an almost gentle algebra, including the string almost gentle case, is Gorenstein if and only if every connected component of its	relation quiver is perfect or acyclic. These algebras need not be Gorenstein.
				For example, take the string almost gentle algebra $kQ/I$ with 
				\[
					Q:\quad
					\begin{tikzpicture}[baseline=(v1.base),>=stealth]
						\node (v1) at (0,0) {$1$};
						\node (v2) at (1.8,0) {$2$};
						\draw[->] (v1.north west)
							.. controls (-.45,.60) and (-.95,.45) .. (-.95,0)
							.. controls (-.95,-.45) and (-.45,-.60) .. (v1.south west);
						\node[above] at (-.60,.48) {$a$};
						\draw[->] (v1) -- node[above] {$b$} (v2);
					\end{tikzpicture}
					\qquad I=(a^2,ab).
				\]

				\item Note that for almost gentle algebras, the cycles in Proposition~\ref{proposition of forbidden cycle criterion} are called \emph{forbidden cycles} in \cite{WHL25}. Proposition~\ref{proposition of forbidden cycle criterion} recovers the equivalence of conditions~(1) and~(3)
				in \cite[Theorem~1.3]{WHL25}.
			\end{enumerate}
		\end{remark}

		For a class $\mathcal X$ of algebras, put
		\[
			\mathcal C_{\mathcal X}=\{\varLambda\mid D_{\mathrm{sg}}(\varLambda)\simeq D_{\mathrm{sg}}(\varGamma)\text{ for some }\varGamma\in\mathcal X\}.
		\]
		Write $\mathcal C_{\mathrm{gent}}$, $\mathcal C_{\mathrm{algent}}$, and $\mathcal C_{\mathrm{Gor}}$ when $\mathcal X$ is the class of gentle, almost gentle, or Gorenstein algebras, respectively.

		By Proposition~\ref{proposition of singularity category via relation quiver}, $\mathcal C_{\mathrm{algent}}$ is precisely the class of algebras singularly equivalent to algebras with radical square zero, since every algebra with radical square zero is almost gentle. The same proposition implies that $D_{\mathrm{sg}}(\varLambda)$ is semisimple for every $\varLambda\in\mathcal C_{\mathrm{algent}}$.

		\begin{proposition}\label{proposition of semisimple singularity categories}
			Let $\varLambda$ be an algebra.	The following conditions are equivalent.
			\begin{enumerate}
				\item $\varLambda\in\mathcal C_{\mathrm{gent}}$.
				\item $D_{\mathrm{sg}}(\varLambda)$ is Hom-finite and semisimple.
				\item There exist positive integers $\ell_1,\ldots,\ell_r$, with $r\geq0$, such that
				\[
					D_{\mathrm{sg}}(\varLambda)\simeq\prod_{i=1}^r\smod N_{\ell_i}.
				\]
			\end{enumerate}
			The multiset $\{\!\{\ell_1,\ldots,\ell_r\}\!\}$ is uniquely determined by the singularity category.
		\end{proposition}

		\begin{proof}
			\underline{$(1)\Rightarrow(3)$} follows from \cite[Theorem~2.5 and Remark~2.6]{Kal15}.

			\underline{$(3)\Rightarrow(1)$.}
			For $r>0$, the selfinjective gentle algebra $\prod_{i=1}^r N_{\ell_i}$ has singularity category $\prod_{i=1}^r\smod N_{\ell_i}$ by Proposition~\ref{proposition of stable and singularity comparison}(3). For $r=0$, use $k$.

			\underline{$(2)\Rightarrow(3)$.}
			Put $E=q_{\varLambda}(\varLambda/\rad\varLambda)$. In a semisimple triangulated category, the middle object of a triangle $X\to Y\to Z\to X[1]$ is a direct summand of $X\oplus Z$. Applying this to composition series gives $q_{\varLambda}M\in\add E$ for every module $M$. In particular,
			\[
				E[-1]\cong
				q_{\varLambda}\Omega_{\varLambda}(\varLambda/\rad\varLambda)
				\in\add E.
			\]
			Since $D_{\mathrm{sg}}(\varLambda)$ is Hom-finite and semisimple, $E$ has finitely
			many pairwise nonisomorphic simple summands $E_1,\ldots,E_n$. The autoequivalence $[-1]$ induces an injection of their isomorphism classes into the same finite set, hence a permutation. Thus $\add E$ is closed under shifts and, by semisimplicity,
			under cones. Since every bounded complex is generated by its terms, $D_{\mathrm{sg}}(\varLambda)=\add E$. Let $\ell_1,\ldots,\ell_r$ be the lengths of the suspension
			orbits of the $E_i$. Since $k$ is algebraically closed, $\End_{D_{\mathrm{sg}}(\varLambda)}(E_i)=k$. The description of $\smod N_\ell$ therefore gives a
			$k$-linear equivalence $D_{\mathrm{sg}}(\varLambda)\simeq\prod_{i=1}^r\smod N_{\ell_i}$ commuting with suspension. This is a triangle equivalence by \cite[Lemma~3.4]{Che11}.

			\underline{$(3)\Rightarrow(2)$} is immediate.

			Finally, a triangle equivalence preserves the suspension orbits of simple objects, so the multiset of their lengths is uniquely determined.
		\end{proof}

		Then we have the following corollary.

		\begin{corollary}\label{corollary of gentle almost gentle singularity classes}
			There are equalities
			\[
				\mathcal C_{\mathrm{gent}}
				=\{\varLambda\in\mathcal C_{\mathrm{algent}}\mid D_{\mathrm{sg}}(\varLambda)\text{ is Hom-finite}\}
				=\mathcal C_{\mathrm{algent}}\cap\mathcal C_{\mathrm{Gor}}.
			\]
		\end{corollary}

		\begin{proof}
			Every algebra in $\mathcal C_{\mathrm{algent}}$ has a semisimple singularity category, so the first equality follows from Proposition~\ref{proposition of semisimple singularity categories}. By Proposition~\ref{proposition of stable and singularity comparison}(2), every algebra in $\mathcal C_{\mathrm{Gor}}$ has a Hom-finite singularity category; hence $\mathcal C_{\mathrm{algent}}\cap\mathcal C_{\mathrm{Gor}}\subseteq\mathcal C_{\mathrm{gent}}$. The reverse inclusion holds because gentle algebras are almost gentle and are Gorenstein.
		\end{proof}

		The preceding results describe the Hom-finite part of $\mathcal C_{\mathrm{algent}}$ and lead to the following question.

		\begin{question}\label{question of semisimple singularity classes}
			Let $\varLambda$ be a finite-dimensional basic $k$-algebra.
			If $D_{\mathrm{sg}}(\varLambda)$ is semisimple, must $\varLambda$ belong to $\mathcal C_{\mathrm{algent}}$?
		\end{question}

		Equivalently, the question asks whether every finite-dimensional algebra with a semisimple singularity category is singularly equivalent to a finite-dimensional algebra with radical square zero.
		These two questions are equivalent by Proposition~\ref{proposition of singularity category via relation quiver} and the fact that every algebra with radical square zero is almost gentle.

	\subsection{Faces and triangle equivalences}

		We describe the stable Gorenstein-projective categories and the Hom-finite singularity categories of almost gentle algebras in terms of their marked configurations. 

		Let $C$ be a marked configuration. Recall that $H=H(C)$ is the set of occurrences at nontruncated vertices, and $D\subseteq H$ consists of the occurrences whose outgoing corners are marked. Thus $U=H\setminus D$ indexes the uncut corners $h\to\sigma h$.

		\begin{definition}\label{definition of relation face}
			A \textbf{face} of $C$ is a primitive closed walk (considered up to cyclic rotation)
			\[
				F:\quad h_0\rightsquigarrow h_1\rightsquigarrow\cdots
				\rightsquigarrow h_{\ell-1}\rightsquigarrow h_0,
				\qquad h_i\in U,\quad \ell\geq1.
			\]
			Write $\mathcal F(C)$ for the set of faces and $\ell(F)=\ell$ for the length of $F$. 
			A face $F$ is \textbf{perfect} if, for every $0\leq i<\ell$,
			\[
				\{g\in U\mid g\rightsquigarrow h_i\}=\{h_{i-1}\},
				\qquad
				\{g\in U\mid h_i\rightsquigarrow g\}=\{h_{i+1}\},
			\]
			where indices are read modulo $\ell$. Write $\mathcal F_{\mathrm{perf}}(C)$ for the set of perfect faces.
		\end{definition}

		By Lemma~\ref{lemma of configuration relation quiver}, faces of $C_A$ correspond to primitive closed walks in $\mathcal R_A$, and perfect faces correspond to its perfect components. Both correspondences preserve lengths. Moreover, in the ribbon graph case, these are precisely the boundary walks that avoid the marked corners.

		As in Subsection~\ref{subsection of relation quiver criterion}, for $\ell\geq1$, let $C_\ell$ be the oriented cycle with $\ell$ vertices, let $J_\ell$ be the arrow ideal of $kC_\ell$, and put $N_\ell=kC_\ell/J_\ell^2$.

		\begin{corollary}\label{proposition of configuration Gorenstein criterion}
			Let $A$ be an almost gentle algebra. There is a triangle equivalence
			\[
				\underline{\operatorname{Gproj}}A
				\simeq
				\prod_{F\in\mathcal F_{\mathrm{perf}}(C_A)}\smod N_{\ell(F)}.
			\]
			Consequently, for almost gentle algebras $A$ and $B$, their stable Gorenstein-projective categories are triangle equivalent if and only if there is a length-preserving bijection $\mathcal F_{\mathrm{perf}}(C_A) \rightarrow\mathcal F_{\mathrm{perf}}(C_B)$.
		\end{corollary}

		\begin{proof}
			The assertions follow directly from Lemma~\ref{lemma of configuration relation quiver} and Theorem~\ref{theorem of Gorenstein projectives and periods}(2)--(3).
		\end{proof}

		To describe the Hom-finite singularity categories, we use the following subset of $U$. An occurrence $h\in U$ is \textbf{essential} if there exists an infinite sequence $(h_i)_{i\in\mathbb Z}$ of occurrences in $U$ such that
		\[
			h_0=h,\qquad h_i\rightsquigarrow h_{i+1}
			\quad\text{for every }i\in\mathbb Z.
		\]
		The occurrences $h_i$ are allowed to repeat. Write $E(C)$ for the set of essential occurrences.

		Equivalently, consider the directed graph with vertex set $U$ and arrows given by $\rightsquigarrow$. At each step, delete all sources and sinks, including isolated vertices, together with their incident arrows. Since $U$ is finite, this process terminates. The remaining vertices form $E(C)$, and the remaining graph is the full subquiver on $E(C)$.

		Every occurrence in a face belongs to $E(C)$, since repeating the closed walk in both directions gives the required sequence.

		\begin{corollary}\label{proposition of configuration singular criterion}
			Let $A$ be an almost gentle algebra. Then $D_{\mathrm{sg}}(A)$ is Hom-finite if and only if, for every occurrence $h\in E(C_A)$, there is exactly one $g\in E(C_A)$ with $g\rightsquigarrow h$ and exactly one $g'\in E(C_A)$ with $h\rightsquigarrow g'$. Under this condition, there is a triangle equivalence
			\[
				D_{\mathrm{sg}}(A)\simeq
				\prod_{F\in\mathcal F(C_A)}\smod N_{\ell(F)}.
			\]
			Consequently, almost gentle algebras $A$ and $B$ with Hom-finite singularity categories are singularly equivalent if and only if there is a length-preserving bijection $\mathcal F(C_A)\longrightarrow\mathcal F(C_B)$.
		\end{corollary}

		\begin{proof}
			Using the identification in Lemma~\ref{lemma of configuration relation quiver}, let $Q_E$ be the full subquiver of $\mathcal R_A$ on $E(C_A)$, and let $J_E$ be its arrow ideal. By Proposition~\ref{proposition of singularity category via relation quiver} and \cite[Lemma~4.3]{Che11}, $D_{\mathrm{sg}}(A)\simeq D_{\mathrm{sg}}(kQ_E/J_E^2)$, where the right-hand side is zero if $Q_E$ is empty.
			By \cite[Theorem~C]{Che11}, this category is Hom-finite if and only if $Q_E$ is a disjoint union of directed cycles, or is empty. This is precisely the stated uniqueness of predecessors and successors within $E(C_A)$.

			Under this condition, the directed cycles of $Q_E$ correspond bijectively to the faces of $C_A$, with the same lengths. Hence
			\[
				kQ_E/J_E^2\cong
				\prod_{F\in\mathcal F(C_A)}N_{\ell(F)}.
			\]
			Since each $N_\ell$ is selfinjective, Proposition~\ref{proposition of stable and singularity comparison}(3) gives the asserted triangle equivalence. The last assertion follows from the uniqueness of the length multiset in Proposition~\ref{proposition of semisimple singularity categories}.
		\end{proof}

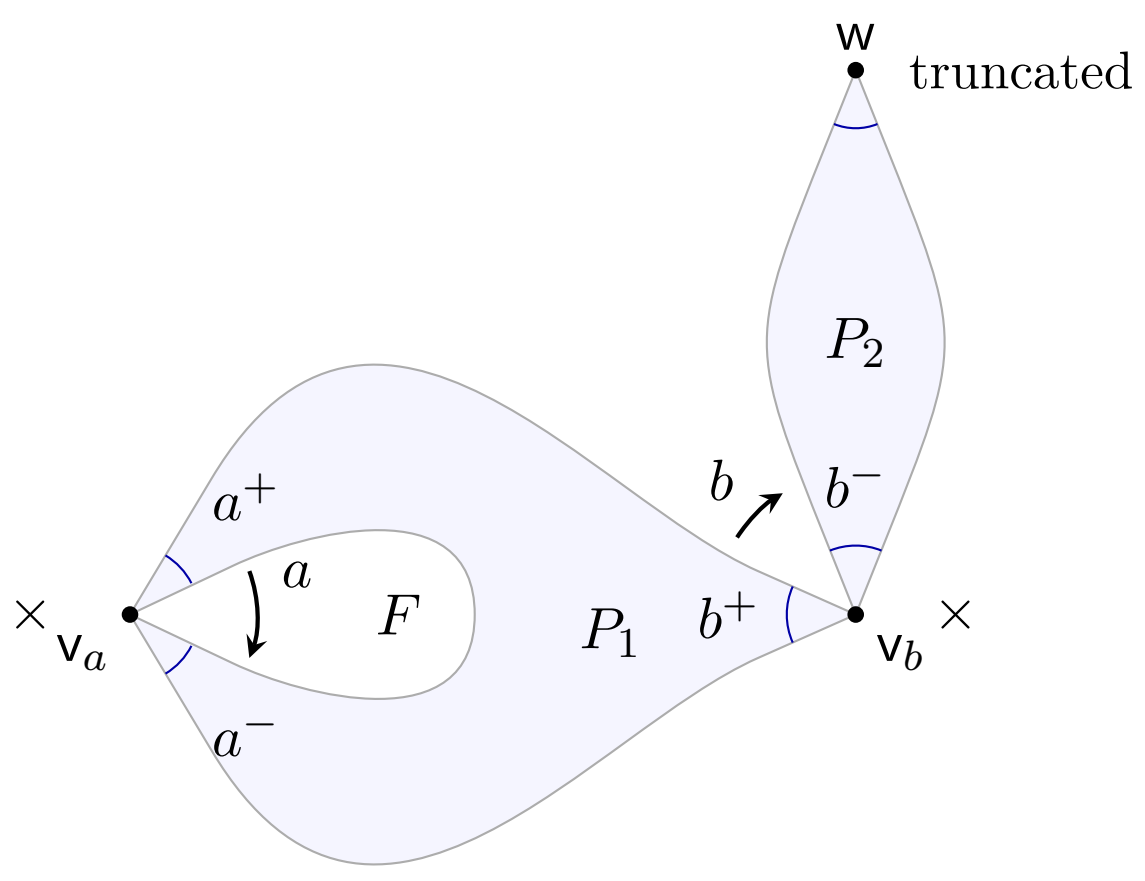
\begin{figure}[htbp]
			\centering
			\begin{tikzpicture}[scale=1.2,>=stealth,
				polygon/.style={fill=blue!4,draw=gray!65},
				occurrence/.style={draw=blue!65!black},
				quiver arrow/.style={->,thick}]
				% Straight segments preserve the angles at configuration vertices;
				% the connecting Bezier curves have matching tangents.
				\filldraw[polygon,even odd rule]
					(0,0)--(.45,.75)
					.. controls (1.35,2.25) and (2.56,.64) .. (3.46,.24)
					--(4,0)--(3.46,-.24)
					.. controls (2.56,-.64) and (1.35,-2.25) .. (.45,-.75)
					--cycle
					(0,0)--(.57,.27)
					.. controls (.95,.45) and (1.9,.70) .. (1.9,0)
					.. controls (1.9,-.70) and (.95,-.45) .. (.57,-.27)
					--cycle;
				\filldraw[polygon]
					(4,0)--(3.84,.4)
					.. controls (3.4,1.5) and (3.4,1.5) .. (3.84,2.6)
					--(4,3)--(4.16,2.6)
					.. controls (4.6,1.5) and (4.6,1.5) .. (4.16,.4)
					--cycle;
				\node at (2.65,-.10) {$P_1$};
				\node at (4,1.50) {$P_2$};
				\node at (1.48,0) {$F$};
				% Occurrences at v_a.
				\draw[occurrence] (27:.38)
					arc[start angle=27,end angle=59,radius=.38];
				\draw[occurrence] (-27:.38)
					arc[start angle=-27,end angle=-59,radius=.38];
				\node at (.64,.65) {$a^+$};
				\node at (.64,-.65) {$a^-$};
				\draw[quiver arrow] (20:.70)
					arc[start angle=20,end angle=-20,radius=.70];
				\node at (.92,.22) {$a$};
				\node at (-.55,0) {$\times$};
				\fill (0,0) circle (1.3pt);
				\node[below left] at (0,0) {$\vv_a$};
				% Occurrences at v_b; b goes clockwise from b^+ to b^-.
				\begin{scope}[shift={(4,0)}]
					\draw[occurrence] (156:.38)
						arc[start angle=156,end angle=204,radius=.38];
					\draw[occurrence] (112:.38)
						arc[start angle=112,end angle=68,radius=.38];
					\draw[quiver arrow] (147:.78)
						arc[start angle=147,end angle=121,radius=.78];
					\node at (135:1.05) {$b$};
				\end{scope}
				\node at (3.30,0) {$b^+$};
				\node at (4,.72) {$b^-$};
				\node at (4.55,0) {$\times$};
				\fill (4,0) circle (1.3pt);
				\node[below right] at (4,0) {$\vv_b$};
				% The occurrence at w is not in H(C_R).
				\begin{scope}[shift={(4,3)}]
					\draw[occurrence] (248:.32)
						arc[start angle=248,end angle=292,radius=.32];
				\end{scope}
				\fill (4,3) circle (1.3pt);
				\node[above] at (4,3) {$\ww$};
				\node[anchor=west,font=\small] at (4.18,3) {truncated};
			\end{tikzpicture}
			\caption{The marked configuration $C_R$ and its unique face $F$.}
			\label{figure of Hom finite non Gorenstein configuration}
		\end{figure}

		\begin{remark}\label{remark of configuration further singular results}
			For an almost gentle algebra $A$, the Hom-finiteness of $D_{\mathrm{sg}}(A)$ is equivalent to $A\in\mathcal C_{\mathrm{gent}}$ by Corollary~\ref{corollary of gentle almost gentle singularity classes}, but does not imply that $A$ is Gorenstein. For example, let $R=kQ/(a^2,ab)$, where
			\[
				Q:\quad
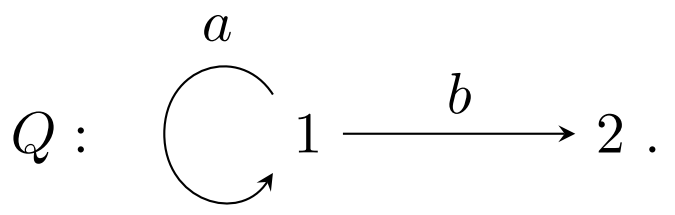
				\begin{tikzpicture}[baseline=(one.base),>=stealth]
					\node (one) at (0,0) {$1$};
					\node (two) at (2,0) {$2$};
					\draw[->] (one.north west)
						.. controls (-.45,.60) and (-.95,.45) .. (-.95,0)
						.. controls (-.95,-.45) and (-.45,-.60) .. (one.south west);
					\node[above] at (-.60,.48) {$a$};
					\draw[->] (one)--node[above] {$b$}(two);
				\end{tikzpicture}.
			\]
			The maximal nonzero paths are $a$ and $b$. The corresponding configuration vertices $\vv_a$ and $\vv_b$ have cut-linear lists $a^+<a^-$ and $b^+<b^-$, respectively, and
			\[
				H_{P_1}=\{a^+,a^-,b^+\},\qquad H_{P_2}=\{b^-\}.
			\]
			The polygon $P_2$ also contains a truncated vertex $\ww$; see Figure~\ref{figure of Hom finite non Gorenstein configuration}.

			Here $U=\{a^+,b^+\}$, and the transitions are $a^+\rightsquigarrow a^+$ and $a^+\rightsquigarrow b^+$. Thus the relation quiver has a loop at $a$ and an arrow $a\to b$. The unique face is $F:a^+\rightsquigarrow a^+$, of length one. It is not perfect, since $a^+$ has two successors in $U$. Hence $R$ is not Gorenstein by Proposition~\ref{proposition of relation quiver Gorenstein criterion}(1). Deleting the sink $b^+$ from the transition graph leaves only the loop at $a^+$, so $E(C_R)=\{a^+\}$. Corollaries~\ref{proposition of configuration Gorenstein criterion} and~\ref{proposition of configuration singular criterion} give
			\[
				\underline{\operatorname{Gproj}}R=0,
				\qquad D_{\mathrm{sg}}(R)\simeq\smod N_1.
			\]
			In particular, $D_{\mathrm{sg}}(R)$ is Hom-finite.

			We also note that if $A$ is Gorenstein, every face of $C_A$ is perfect by Proposition~\ref{proposition of relation quiver Gorenstein criterion}(1). The decompositions in Corollaries~\ref{proposition of configuration Gorenstein criterion} and~\ref{proposition of configuration singular criterion} therefore agree under the canonical equivalence $\underline{\operatorname{Gproj}}A\simeq D_{\mathrm{sg}}(A)$ from Proposition~\ref{proposition of stable and singularity comparison}(2).
		\end{remark}

		\begin{remark}\label{remark of face lengths without Hom finiteness}
			The Hom-finiteness assumption is essential for classifying singularity categories by face lengths. Consider
			\[
				A=k\langle x,y\rangle/(x^2,y^2,yx),
				\qquad
				B=k[x]/(x^2)\times k[y]/(y^2).
			\]
			Both $C_A$ and $C_B$ have exactly two faces, each of length one. The relation quiver of $A$ has loops at $x$ and $y$ and an additional arrow $y\to x$.
			Both vertices remain after deleting sources and sinks, but $x$ has two predecessors and $y$ has two successors. Thus $D_{\mathrm{sg}}(A)$ is not Hom-finite by Corollary~\ref{proposition of configuration singular criterion}. On the other hand, $D_{\mathrm{sg}}(B)\simeq\smod N_1\times\smod N_1$ is Hom-finite.
			Hence $A$ and $B$ are not singularly equivalent, although their face-length multisets coincide.
		\end{remark}

		We can also describe node deletion in terms of perfect faces. We say that a face $F:h_0\rightsquigarrow\cdots\rightsquigarrow h_{\ell-1} \rightsquigarrow h_0$ \textbf{passes through a polygon} $P$ if $P=\pi(\sigma h_i)=\pi(h_{i+1})$ for some $i$, with indices read modulo $\ell$. 

		By Proposition~\ref{proposition of node deletion Gorenstein projectives}, deleting a node removes precisely the perfect faces passing through its polygon and creates no new perfect faces.
		An incidence comparison $C_A\to C_B$ induces an isomorphism $C_{\widehat A}\cong C_{\widehat B}$, so it preserves the lengths of the perfect faces that remain after all nodes are deleted. This recovers Corollary~\ref{corollary of stable Gorenstein quotient invariant}.

		The following example shows that the lengths of the removed faces need not be preserved.

		\begin{example}\label{example of marked ribbon comparison}
			Let $A=k[x]/(x^2)\times k[y]/(y^2)$ and $B=N_2=kQ/(ab,ba)$, where $Q$ has arrows $a:1\to2$ and $b:2\to1$, as in Example~\ref{example of node gluing different periods}. For $z\in\{x,y,a,b\}$, the maximal path $z$ gives a configuration vertex $\vv_z$ with cut-linear list $z^+<z^-$. The polygon incidences are
			\[
				\begin{aligned}
					C_A:&\qquad H_{P_x}=\{x^+,x^-\},\quad
					H_{P_y}=\{y^+,y^-\},\\
					C_B:&\qquad H_{P_1}=\{a^+,b^-\},\quad
					H_{P_2}=\{a^-,b^+\}.
				\end{aligned}
			\]
			Every quiver vertex of $A$ and $B$ is a node. Deleting all nodes separates each displayed pair into two singleton sets of nontruncated occurrences. Each resulting polygon is completed by adding a truncated vertex. After identifying $a^\pm$ with $x^\pm$ and $b^\pm$ with $y^\pm$, we obtain the common marked configuration in Figure~\ref{figure of marked configuration comparison}.
			Thus the assignments $x^\pm\mapsto a^\pm$ and $y^\pm\mapsto b^\pm$ define an incidence comparison $C_A\to C_B$. By Proposition~\ref{proposition of configuration stable criterion}, $A$ and $B$ are stably equivalent.

			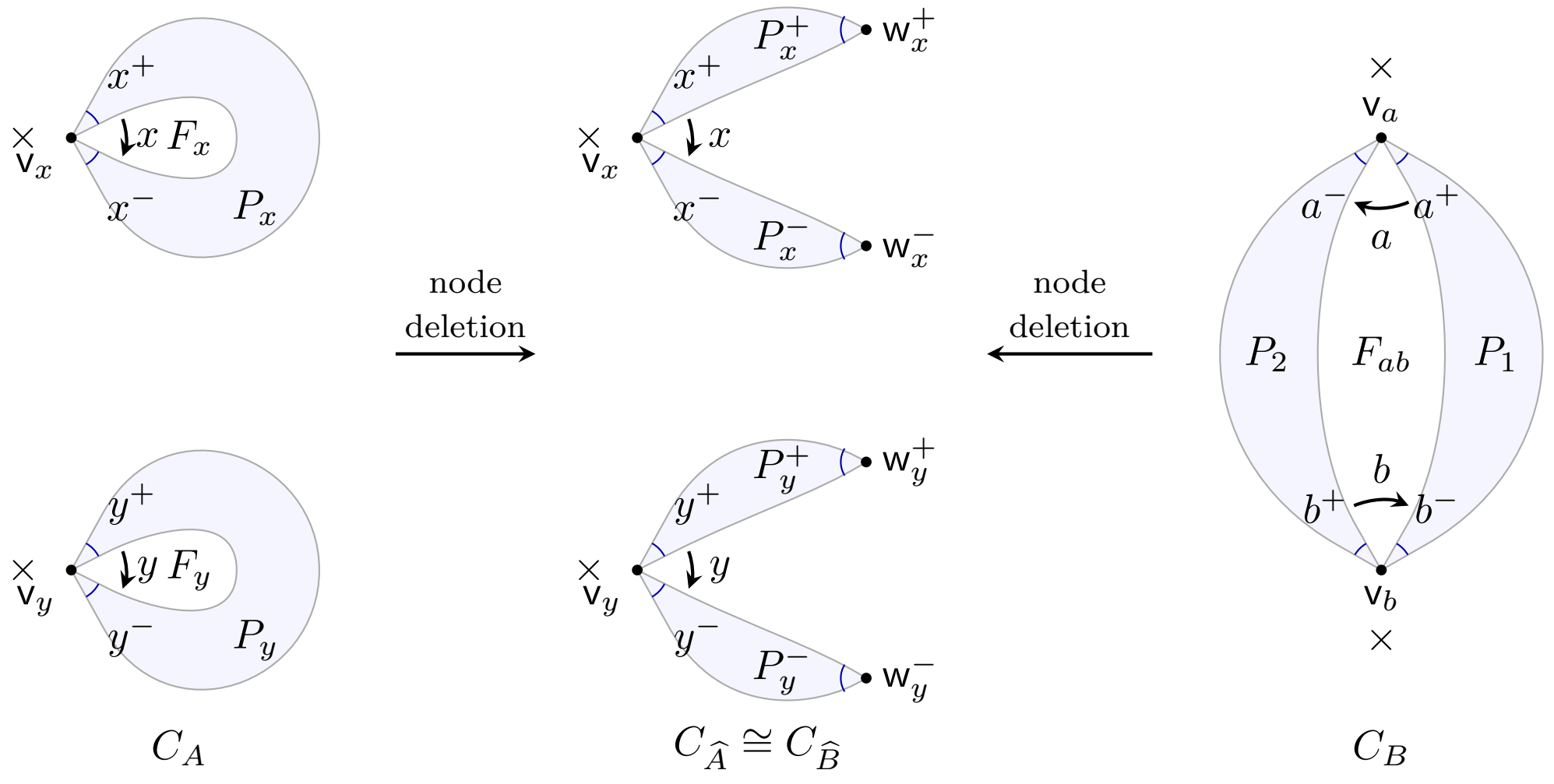
\begin{figure}[htbp]
				\centering
				\begin{tikzpicture}[scale=1.08,>=stealth,font=\small,
					polygon/.style={fill=blue!4,draw=gray!65},
					occurrence/.style={draw=blue!65!black},
					quiver arrow/.style={->,thick}]
					% C_A: a polygon with two occurrences at each nontruncated vertex.
					\foreach \yy/\z in {1.7/x,-1.7/y}{
						\begin{scope}[shift={(0,\yy)}]
							\filldraw[polygon,even odd rule]
								(0,0)--(.25,.45)
								.. controls (.75,1.35) and (1.95,.90) .. (1.95,0)
								.. controls (1.95,-.90) and (.75,-1.35) .. (.25,-.45)
								--cycle
								(0,0)--(.30,.15)
								.. controls (.60,.30) and (1.30,.50) .. (1.30,0)
								.. controls (1.30,-.50) and (.60,-.30) .. (.30,-.15)
								--cycle;
							\node at (1.45,-.55) {$P_{\z}$};
							\node at (.93,0) {$F_{\z}$};
							\draw[occurrence] (27:.24)
								arc[start angle=27,end angle=61,radius=.24];
							\draw[occurrence] (-27:.24)
								arc[start angle=-27,end angle=-61,radius=.24];
							\node at (.48,.52) {$\z^+$};
							\node at (.48,-.52) {$\z^-$};
							\draw[quiver arrow] (20:.43)
								arc[start angle=20,end angle=-20,radius=.43];
							\node at (.60,0) {$\z$};
							\node at (-.38,0) {$\times$};
							\fill (0,0) circle (1.2pt);
							\node[below left] at (0,0) {$\vv_{\z}$};
						\end{scope}
					}
					% The common configuration after deleting all nodes.
					% Each singleton set of nontruncated occurrences is completed by a
					% new truncated vertex w_z^+ or w_z^-.
					\foreach \yy/\z in {1.7/x,-1.7/y}{
						\begin{scope}[shift={(4.45,\yy)}]
							\filldraw[polygon]
								(0,0)--(.25,.45)
								.. controls (.65,1.17) and (1.40,1.10) .. (1.80,.85)
								.. controls (1.40,.60) and (.90,.45) .. (.30,.15)
								--cycle;
							\filldraw[polygon]
								(0,0)--(.30,-.15)
								.. controls (.90,-.45) and (1.40,-.60) .. (1.80,-.85)
								.. controls (1.40,-1.10) and (.65,-1.17) .. (.25,-.45)
								--cycle;
							\node at (1.14,.79) {$P_{\z}^{+}$};
							\node at (1.14,-.79) {$P_{\z}^{-}$};
							\draw[occurrence] (27:.24)
								arc[start angle=27,end angle=61,radius=.24];
							\draw[occurrence] (-27:.24)
								arc[start angle=-27,end angle=-61,radius=.24];
							\node at (.48,.52) {$\z^+$};
							\node at (.48,-.52) {$\z^-$};
							\draw[quiver arrow] (20:.43)
								arc[start angle=20,end angle=-20,radius=.43];
							\node at (.65,0) {$\z$};
							\node at (-.38,0) {$\times$};
							\fill (0,0) circle (1.2pt);
							\node[below left] at (0,0) {$\vv_{\z}$};
							\foreach \sign/\height in {+/.85,-/-.85}{
								\begin{scope}[shift={(1.8,\height)}]
									\draw[occurrence] (148:.20)
										arc[start angle=148,end angle=212,radius=.20];
									\fill (0,0) circle (1.2pt);
									\node[right] at (0,0) {$\ww_{\z}^{\sign}$};
								\end{scope}
							}
						\end{scope}
					}
					% C_B: two polygons joining v_a and v_b, with a face between them.
					\begin{scope}[shift={(10.3,0)}]
						\filldraw[polygon]
							(0,1.7)--(.35,1.5)
							.. controls (1.575,.80) and (1.575,-.80) .. (.35,-1.5)
							--(0,-1.7)--(.20,-1.35)
							.. controls (.60,-.65) and (.60,.65) .. (.20,1.35)
							--cycle;
						\filldraw[polygon]
							(0,1.7)--(-.20,1.35)
							.. controls (-.60,.65) and (-.60,-.65) .. (-.20,-1.35)
							--(0,-1.7)--(-.35,-1.5)
							.. controls (-1.575,-.80) and (-1.575,.80) .. (-.35,1.5)
							--cycle;
						\node at (.90,0) {$P_1$};
						\node at (-.90,0) {$P_2$};
						\node at (0,0) {$F_{ab}$};
						\begin{scope}[shift={(0,1.7)}]
							\draw[occurrence] (-60:.24)
								arc[start angle=-60,end angle=-30,radius=.24];
							\draw[occurrence] (-150:.24)
								arc[start angle=-150,end angle=-120,radius=.24];
							\draw[quiver arrow] (-67:.55)
								arc[start angle=-67,end angle=-113,radius=.55];
							\node at (0,-.80) {$a$};
							\node at (.44,-.50) {$a^+$};
							\node at (-.44,-.50) {$a^-$};
							\fill (0,0) circle (1.2pt);
							\node[above] at (0,0) {$\vv_a$};
							\node at (0,.55) {$\times$};
						\end{scope}
						\begin{scope}[shift={(0,-1.7)}]
							\draw[occurrence] (30:.24)
								arc[start angle=30,end angle=60,radius=.24];
							\draw[occurrence] (120:.24)
								arc[start angle=120,end angle=150,radius=.24];
							\draw[quiver arrow] (113:.55)
								arc[start angle=113,end angle=67,radius=.55];
							\node at (0,.80) {$b$};
							\node at (-.44,.50) {$b^+$};
							\node at (.44,.50) {$b^-$};
							\fill (0,0) circle (1.2pt);
							\node[below] at (0,0) {$\vv_b$};
							\node at (0,-.55) {$\times$};
						\end{scope}
					\end{scope}
					\draw[->,thick] (2.55,0)--(3.65,0)
						node[midway,above,align=center,font=\scriptsize] {node\\deletion};
					\draw[->,thick] (8.50,0)--(7.20,0)
						node[midway,above,align=center,font=\scriptsize] {node\\deletion};
					\node at (.85,-3.10) {$C_A$};
					\node at (5.40,-3.10) {$C_{\widehat A}\cong C_{\widehat B}$};
					\node at (10.3,-3.10) {$C_B$};
				\end{tikzpicture}
				\caption{The marked configurations $C_A$ and $C_B$ and their common
					configuration after node deletion.
					The vertices $\ww_z^\pm$ are truncated.}
				\label{figure of marked configuration comparison}
			\end{figure}

			The faces of $C_A$ are
			\[
				F_x:\ x^+\rightsquigarrow x^+,\qquad
				F_y:\ y^+\rightsquigarrow y^+,
			\]
			whereas $C_B$ has the single face
			\[
				F_{ab}:\ a^+\rightsquigarrow b^+\rightsquigarrow a^+.
			\]
			All three faces are perfect, but their length multisets are $\{\!\{1,1\}\!\}$ and $\{\!\{2\}\!\}$, respectively. Corollary~\ref{proposition of configuration Gorenstein criterion} therefore shows that the stable Gorenstein-projective categories of $A$ and $B$ are not triangle equivalent. Since both algebras are selfinjective, they are not singularly equivalent either, by Proposition~\ref{proposition of stable and singularity comparison}(3).
		\end{example}

\section*{Acknowledgments}

	We express our gratitude to our supervisor, Professor Yuming Liu, for helpful discussions, for contributing key ideas to the proof of Lemma~\ref{lemma of stable equivalences induces isomorphisms}, which shows that the stably equivalent algebras in the reduced case are isomorphic, and for sharing his knowledge of lifting problems for stable equivalences of Morita type.

	We thank Zhen Zhang for suggesting that we investigate whether our results extend to arbitrary quadratic monomial algebras.
	The initial draft focused on string almost gentle algebras, and at that stage we expected the argument to extend only to almost gentle algebras.

	We acknowledge the assistance of GPT-6 Astra in extending the argument to arbitrary quadratic monomial algebras. The authors have independently checked all proofs and take full responsibility for the mathematical results and the final manuscript.

\begingroup\linespread{1}\selectfont
\bibliographystyle{alpha}
\bibliography{ref}

@article{AH81,
  author       = {Assem, Ibrahim and Happel, Dieter},
  title        = {Generalized tilted algebras of type {$A_n$}},
  journal      = {Comm. Algebra},
  volume       = {9},
  number       = {20},
  pages        = {2101--2125},
  year         = {1981},
  doi          = {10.1080/00927878108822697}
}

@article{Ant07,
  author       = {Antipov, Mikhail},
  title        = {Invariants of the stable equivalence of symmetric special biserial algebras},
  journal      = {J. Math. Sci. (N.Y.)},
  volume       = {140},
  number       = {5},
  pages        = {611--621},
  year         = {2007}
}

@book{ARS,
  author       = {Auslander, Maurice and Reiten, Idun and Smal{\o}, Sverre Olaf},
  title        = {Representation theory of {Artin} algebras},
  series       = {Cambridge Studies in Advanced Mathematics},
  volume       = {36},
  publisher    = {Cambridge University Press},
  year         = {1997}
}

@article{AS87,
  author       = {Assem, Ibrahim and Skowro{\'n}ski, Andrzej},
  title        = {Iterated tilted algebras of type {$\widetilde{A}_n$}},
  journal      = {Math. Z.},
  volume       = {195},
  number       = {2},
  pages        = {269--290},
  year         = {1987},
  doi          = {10.1007/BF01166463}
}

@article{BJO15,
  author       = {Bergh, Petter Andreas and Jorgensen, David Alan and Oppermann, Steffen},
  title        = {The {Gorenstein} defect category},
  journal      = {Q. J. Math.},
  volume       = {66},
  number       = {2},
  pages        = {459--471},
  year         = {2015}
}

@incollection{Bro94,
  author       = {Brou{\'e}, Michel},
  title        = {Equivalences of blocks of group algebras},
  editor       = {Dlab, Vlastimil and Scott, Leonard L.},
  booktitle    = {Finite Dimensional Algebras and Related Topics},
  series       = {NATO Advanced Science Institutes Series C: Mathematical and Physical Sciences},
  volume       = {424},
  pages        = {1--26},
  publisher    = {Kluwer Academic Publishers},
  address      = {Dordrecht},
  year         = {1994},
  doi          = {10.1007/978-94-017-1556-0_1}
}

@unpublished{Buc87,
  author       = {Buchweitz, Ragnar-Olaf},
  title        = {Maximal {Cohen--Macaulay} modules and {Tate}-cohomology over {Gorenstein} rings},
  year         = {1987},
  note         = {Unpublished manuscript}
}

@article{Cha26,
  author       = {Chang, Wen},
  title        = {Geometric models for the algebraic hearts in the derived category of a gentle algebra},
  journal      = {Israel J. Math.},
  volume       = {271},
  number       = {2},
  pages        = {941--1000},
  year         = {2026},
  doi          = {10.1007/s11856-025-2801-7}
}

@article{Che11,
  author       = {Chen, Xiao-Wu},
  title        = {The singularity category of an algebra with radical square zero},
  journal      = {Doc. Math.},
  volume       = {16},
  pages        = {921--936},
  year         = {2011}
}

@misc{Che17,
  author       = {Chen, Xiao-Wu},
  title        = {{Gorenstein} Homological Algebra of {Artin} Algebras},
  year         = {2017},
  eprint       = {1712.04587},
  archivePrefix= {arXiv},
  primaryClass = {math.RT},
  url          = {https://arxiv.org/abs/1712.04587},
  note         = {arXiv:\allowbreak 1712.04587}
}

@article{Che18,
  author       = {Chen, Xiao-Wu},
  title        = {The singularity category of a quadratic monomial algebra},
  journal      = {Q. J. Math.},
  volume       = {69},
  number       = {3},
  pages        = {1015--1033},
  year         = {2018}
}

@article{CKL15,
  author       = {Chan, Aaron and Koenig, Steffen and Liu, Yuming},
  title        = {Simple-minded systems, configurations and mutations for representation-finite self-injective algebras},
  journal      = {J. Pure Appl. Algebra},
  volume       = {219},
  number       = {6},
  pages        = {1940--1961},
  year         = {2015},
  doi          = {10.1016/j.jpaa.2014.07.018}
}

@article{CSZ18,
  author       = {Chen, Xiao-Wu and Shen, Dawei and Zhou, Guodong},
  title        = {The {Gorenstein}-projective modules over a monomial algebra},
  journal      = {Proc. Roy. Soc. Edinburgh Sect. A},
  volume       = {148},
  number       = {6},
  pages        = {1115--1134},
  year         = {2018}
}

@article{FGR22,
  author       = {Franco, Andr{\'e}s and Giraldo, Hern{\'a}n and Rizzo, Pedro},
  title        = {String and Band Complexes Over String Almost Gentle Algebras},
  journal      = {Appl. Categ. Structures},
  volume       = {30},
  pages        = {417--452},
  year         = {2022},
  doi          = {10.1007/s10485-021-09661-x}
}

@incollection{GR05,
  author       = {Gei{\ss}, Christof and Reiten, Idun},
  title        = {Gentle algebras are {Gorenstein}},
  booktitle    = {Representations of Algebras and Related Topics},
  series       = {Fields Institute Communications},
  volume       = {45},
  pages        = {129--133},
  publisher    = {American Mathematical Society},
  address      = {Providence, RI},
  year         = {2005}
}

@article{GS16MS,
  author       = {Green, Edward L. and Schroll, Sibylle},
  title        = {Multiserial and special multiserial algebras and their representations},
  journal      = {Adv. Math.},
  volume       = {302},
  pages        = {1111--1136},
  year         = {2016},
  doi          = {10.1016/j.aim.2016.07.006}
}

@article{GS17BC,
  author       = {Green, Edward L. and Schroll, Sibylle},
  title        = {{Brauer} configuration algebras: A generalization of {Brauer} graph algebras},
  journal      = {Bull. Sci. Math.},
  volume       = {141},
  number       = {6},
  pages        = {539--572},
  year         = {2017},
  doi          = {10.1016/j.bulsci.2017.06.001}
}

@article{GS19,
  author       = {Green, Edward L. and Schroll, Sibylle},
  title        = {Almost gentle algebras and their trivial extensions},
  journal      = {Proc. Edinb. Math. Soc. (2)},
  volume       = {62},
  number       = {2},
  pages        = {489--504},
  year         = {2019}
}

@incollection{Hap91,
  author       = {Happel, Dieter},
  title        = {On {Gorenstein} algebras},
  editor       = {Michler, Gerhard O. and Ringel, Claus Michael},
  booktitle    = {Representation Theory of Finite Groups and Finite-Dimensional Algebras},
  series       = {Progress in Mathematics},
  volume       = {95},
  pages        = {389--404},
  publisher    = {Birkh{\"a}user},
  address      = {Basel},
  year         = {1991},
  doi          = {10.1007/978-3-0348-8658-1_16}
}

@article{HKK17,
  author       = {Haiden, Fabian and Katzarkov, Ludmil and Kontsevich, Maxim},
  title        = {Flat surfaces and stability structures},
  journal      = {Publ. Math. Inst. Hautes {\'E}tudes Sci.},
  volume       = {126},
  pages        = {247--318},
  year         = {2017},
  doi          = {10.1007/s10240-017-0095-y}
}

@article{HU25,
  author       = {Honma, Takahiro and Usui, Satoshi},
  title        = {The stable category of {Gorenstein}-projective modules over a monomial algebra},
  journal      = {J. Lond. Math. Soc. (2)},
  volume       = {111},
  number       = {6},
  pages        = {e70204},
  year         = {2025},
  doi          = {10.1112/jlms.70204}
}

@article{Kal15,
  author       = {Kalck, Martin},
  title        = {Singularity categories of gentle algebras},
  journal      = {Bull. Lond. Math. Soc.},
  volume       = {47},
  number       = {1},
  pages        = {65--74},
  year         = {2015}
}

@article{Kel94,
  author       = {Keller, Bernhard},
  title        = {Deriving {DG} categories},
  journal      = {Ann. Sci. {\'E}c. Norm. Sup{\'e}r. (4)},
  volume       = {27},
  number       = {1},
  pages        = {63--102},
  year         = {1994},
  doi          = {10.24033/asens.1689}
}

@incollection{Kel07,
  author       = {Keller, Bernhard},
  title        = {Derived categories and tilting},
  editor       = {Angeleri H{\"u}gel, Lidia and Happel, Dieter and Krause, Henning},
  booktitle    = {Handbook of Tilting Theory},
  series       = {London Mathematical Society Lecture Note Series},
  volume       = {332},
  pages        = {49--104},
  publisher    = {Cambridge University Press},
  year         = {2007},
  doi          = {10.1017/CBO9780511735134.005}
}

@article{KL12,
  author       = {Koenig, Steffen and Liu, Yuming},
  title        = {Simple-minded systems in stable module categories},
  journal      = {Q. J. Math.},
  volume       = {63},
  number       = {3},
  pages        = {653--674},
  year         = {2012},
  doi          = {10.1093/qmath/har009}
}

@article{Kra91,
  author       = {Krause, Henning},
  title        = {Maps between Tree and Band Modules},
  journal      = {J. Algebra},
  volume       = {137},
  number       = {1},
  pages        = {186--194},
  year         = {1991}
}

@article{KV87,
  author       = {Keller, Bernhard and Vossieck, Dieter},
  title        = {Sous les cat{\'e}gories d{\'e}riv{\'e}es},
  journal      = {C. R. Acad. Sci. Paris S{\'e}r. I Math.},
  volume       = {305},
  number       = {6},
  pages        = {225--228},
  year         = {1987}
}

@article{KY14,
  author       = {Koenig, Steffen and Yang, Dong},
  title        = {Silting objects, simple-minded collections, {$t$}-structures and co-{$t$}-structures for finite-dimensional algebras},
  journal      = {Doc. Math.},
  volume       = {19},
  pages        = {403--438},
  year         = {2014},
  doi          = {10.4171/DM/451}
}

@article{KY16,
  author       = {Kalck, Martin and Yang, Dong},
  title        = {Relative singularity categories {I}: {Auslander} resolutions},
  journal      = {Adv. Math.},
  volume       = {301},
  pages        = {973--1021},
  year         = {2016}
}

@article{Lin96,
  author       = {Linckelmann, Markus},
  title        = {Stable equivalences of {Morita} type for self-injective algebras and {$p$}-groups},
  journal      = {Math. Z.},
  volume       = {223},
  number       = {1},
  pages        = {87--100},
  year         = {1996},
  doi          = {10.1007/BF02621590}
}

@article{Liu03,
  author       = {Liu, Yuming},
  title        = {On stable equivalences of {Morita} type for finite dimensional algebras},
  journal      = {Proc. Amer. Math. Soc.},
  volume       = {131},
  number       = {9},
  pages        = {2657--2662},
  year         = {2003}
}

@article{Liu06,
  author       = {Liu, Yuming},
  title        = {On stable equivalences induced by exact functors},
  journal      = {Proc. Amer. Math. Soc.},
  volume       = {134},
  number       = {6},
  pages        = {1605--1613},
  year         = {2006},
  doi          = {10.1090/S0002-9939-05-08157-8}
}

@article{LP20,
  author       = {Lekili, Yank{\i} and Polishchuk, Alexander},
  title        = {Derived equivalences of gentle algebras via {Fukaya} categories},
  journal      = {Math. Ann.},
  volume       = {376},
  number       = {1--2},
  pages        = {187--225},
  year         = {2020},
  doi          = {10.1007/s00208-019-01894-5}
}

@article{LX26,
  author       = {Liu, Yuming and Xing, Bohan},
  title        = {Trivial extensions of monomial algebras are symmetric fractional {Brauer} configuration algebras},
  journal      = {Algebr. Represent. Theory},
  volume       = {29},
  pages        = {183--202},
  year         = {2026},
  doi          = {10.1007/s10468-026-10378-5}
}

@incollection{Mar80,
  author       = {Mart{\'i}nez{-}Villa, Roberto},
  title        = {Algebras stably equivalent to {$l$}-hereditary},
  editor       = {Dlab, Vlastimil and Gabriel, Peter},
  booktitle    = {Representation Theory {II}},
  series       = {Lecture Notes in Mathematics},
  volume       = {832},
  pages        = {396--431},
  publisher    = {Springer-Verlag},
  address      = {Berlin, Heidelberg},
  year         = {1980}
}

@article{Mar90,
  author       = {Mart{\'i}nez{-}Villa, Roberto},
  title        = {Properties that are left invariant under stable equivalence},
  journal      = {Comm. Algebra},
  volume       = {18},
  number       = {12},
  pages        = {4141--4169},
  year         = {1990}
}

@article{MM05,
  author       = {Mart{\'i}nez{-}Villa, Roberto and Monta{\~n}o{-}Berm{\'u}dez, Gustavo},
  title        = {Node deletion and stably equivalent {Koszul} algebras},
  journal      = {Comm. Algebra},
  volume       = {33},
  number       = {8},
  pages        = {2569--2585},
  year         = {2005}
}

@article{Ohn02,
  author       = {Ohnuki, Yosuke},
  title        = {Stable equivalence induced by a socle equivalence},
  journal      = {Osaka J. Math.},
  volume       = {39},
  number       = {1},
  pages        = {259--266},
  year         = {2002}
}

@misc{OPS,
  author       = {Opper, Sebastian and Plamondon, Pierre-Guy and Schroll, Sibylle},
  title        = {A geometric model for the derived category of gentle algebras},
  year         = {2018},
  eprint       = {1801.09659},
  archivePrefix= {arXiv},
  primaryClass = {math.RT},
  url          = {https://arxiv.org/abs/1801.09659},
  note         = {arXiv:\allowbreak 1801.09659}
}

@article{Orl04,
  author       = {Orlov, Dmitri},
  title        = {Triangulated categories of singularities and {D}-branes in {Landau--Ginzburg} models},
  journal      = {Proc. Steklov Inst. Math.},
  volume       = {246},
  number       = {3},
  pages        = {227--248},
  year         = {2004}
}

@article{Ric89,
  author       = {Rickard, Jeremy},
  title        = {{Morita} theory for derived categories},
  journal      = {J. Lond. Math. Soc. (2)},
  volume       = {39},
  number       = {3},
  pages        = {436--456},
  year         = {1989}
}

@article{Ric02,
  author       = {Rickard, Jeremy},
  title        = {Equivalences of derived categories for symmetric algebras},
  journal      = {J. Algebra},
  volume       = {257},
  number       = {2},
  pages        = {460--481},
  year         = {2002},
  doi          = {10.1016/S0021-8693(02)00520-3}
}

@article{Rin12,
  author       = {Ringel, Claus Michael},
  title        = {On the representation dimension of {Artin} algebras},
  journal      = {Bull. Inst. Math. Acad. Sin. (N.S.)},
  volume       = {7},
  number       = {1},
  pages        = {33--70},
  year         = {2012}
}

@article{Sch15,
  author       = {Schroll, Sibylle},
  title        = {Trivial extensions of gentle algebras and {Brauer} graph algebras},
  journal      = {J. Algebra},
  volume       = {444},
  pages        = {183--200},
  year         = {2015},
  doi          = {10.1016/j.jalgebra.2015.07.037}
}

@incollection{Sch26,
  author       = {Schroll, Sibylle},
  title        = {On geometric models in representation theory},
  editor       = {Friedlander, Susan and Tschinkel, Yuri},
  booktitle    = {Proceedings of the International Congress of Mathematicians 2026. Vol.~3: Invited Lectures (Sections 1--4)},
  pages        = {250--271},
  publisher    = {Society for Industrial and Applied Mathematics},
  year         = {2026},
  doi          = {10.1137/25M1810647}
}

@article{SW83,
  author       = {Skowro{\'n}ski, Andrzej and Waschb{\"u}sch, Josef},
  title        = {Representation-finite biserial algebras},
  journal      = {J. Reine Angew. Math.},
  volume       = {345},
  pages        = {172--181},
  year         = {1983},
  doi          = {10.1515/crll.1983.345.172}
}

@article{SZ03,
  author       = {Schr{\"o}er, Jan and Zimmermann, Alexander},
  title        = {Stable endomorphism algebras of modules over special biserial algebras},
  journal      = {Math. Z.},
  volume       = {244},
  number       = {3},
  pages        = {515--530},
  year         = {2003}
}

@article{HW84,
  author       = {von H{\"o}hne, Hans-Joachim and Waschb{\"u}sch, Josef},
  title        = {{Die Struktur} {$n$}-reihiger {Algebren}},
  journal      = {Comm. Algebra},
  volume       = {12},
  number       = {9--10},
  pages        = {1187--1206},
  year         = {1984},
  doi          = {10.1080/00927878408823049}
}

@misc{WHL25,
  author       = {Wang, Demin and Huang, Zhaoyong and Liu, Yu-Zhe},
  title        = {Homological dimensions over almost gentle algebras},
  year         = {2025},
  eprint       = {2503.04445},
  archivePrefix= {arXiv},
  primaryClass = {math.RT},
  url          = {https://arxiv.org/abs/2503.04445v1},
  note         = {arXiv:\allowbreak 2503.04445v1}
}

@article{WW85,
  author       = {Wald, Burkhard and Waschb{\"u}sch, Josef},
  title        = {Tame biserial algebras},
  journal      = {J. Algebra},
  volume       = {95},
  number       = {2},
  pages        = {480--500},
  year         = {1985},
  doi          = {10.1016/0021-8693(85)90119-X}
}

@article{Zha24BGA,
  author       = {Zhang, Zhen},
  title        = {On simple-minded systems over domestic {Brauer} graph algebras},
  journal      = {Eur. J. Math.},
  volume       = {10},
  year         = {2024},
  doi          = {10.1007/s40879-023-00717-x},
  note         = {Article no.~1}
}

@misc{Zha25,
  author       = {Zhang, Zhen},
  title        = {A construction of simple-minded systems over domestic {Brauer} graph algebras {I}: the 2-domestic case},
  year         = {2025},
  eprint       = {2509.24184},
  archivePrefix= {arXiv},
  primaryClass = {math.RT},
  url          = {https://arxiv.org/abs/2509.24184},
  note         = {arXiv:\allowbreak 2509.24184}
}

@misc{Zha26BGAII,
  author       = {Zhang, Zhen},
  title        = {A construction of simple-minded systems over domestic {Brauer} graph algebras {II}: the 1-domestic case},
  year         = {2026},
  eprint       = {2606.16085},
  archivePrefix= {arXiv},
  primaryClass = {math.RT},
  url          = {https://arxiv.org/abs/2606.16085},
  note         = {arXiv:\allowbreak 2606.16085}
}

@book{Zim14,
  author       = {Zimmermann, Alexander},
  title        = {Representation Theory: A Homological Algebra Point of View},
  series       = {Algebra and Applications},
  volume       = {19},
  publisher    = {Springer},
  address      = {Cham},
  year         = {2014},
  doi          = {10.1007/978-3-319-07968-4}
}
\endgroup
\end{document}